\documentclass[twoside]{article}
\usepackage{eurosym}
\usepackage[utf8]{inputenc}
\usepackage{graphicx}
\usepackage{latexsym}
\usepackage{amsmath}
\usepackage{amsfonts}
\usepackage{amssymb}
\usepackage{hyperref}
\usepackage[usenames]{color}

\newcommand{\disp}{\displaystyle}
\newcommand{\proof}{\par
	\noindent {\sc Proof}.\quad}

\newcommand{\cA}{{\cal A}}

\newcommand{\rA}{{\rm A}}

\newcommand{\rB}{{\rm B}}

\newcommand{\bB}{{\bf B}}

\newcommand{\cC}{{\cal C}}
\newcommand{\rC}{{\rm C}}
\newcommand{\rCh}{\widehat{\rm C}}
\newcommand{\rd}{{\rm d}}
\newcommand{\rD}{{\rm D}}

\newcommand{\rE}{{\rm E}}

\newcommand{\rFh}{\widehat{\rm F}}
\newcommand{\cF}{{\cal F}}

\newcommand{\rG}{{\rm G}}
\newcommand{\cG}{{\cal G}}

\newcommand{\rH}{{\rm H}}

\newcommand{\bh}{{\bf h}}
\newcommand{\rI}{{\rm I}}
\newcommand{\mI}{\text{\rm \sffamily \bfseries I}}

\newcommand{\cJ}{{\cal J}}
\newcommand{\rK}{{\rm K}}

\newcommand{\rL}{{\rm L}}

\newcommand{\rM}{{\rm M}}
\newcommand{\rN}{{\rm N}}

\newcommand{\bn}{{\bf n}}

\newcommand{\cO}{{\cal O}}

\newcommand{\bp}{{\bf p}}

\newcommand{\bq}{{\bf q}}
\newcommand{\cQ}{{\cal Q}}

\newcommand{\R}{{\mathbb{R}}}
\newcommand{\rR}{{\rm R}}

\newcommand{\cS}{{\cal S}}

\newcommand{\rU}{{\rm U}}
\newcommand{\cU}{{\cal U}}
\newcommand{\rV}{{\rm V}}
\newcommand{\cV}{{\cal V}}
\newcommand{\rW}{{\rm W}}

\newcommand{\cW}{{\cal W}}

\newcommand{\rX}{{\rm X}}

\newcommand{\rY}{{\rm Y}}

\newcommand{\fin}{\hfill\mbox{$\quad{}_{\Box}$}}
\newcommand{\fineq}{\vspace{-.75cm$\fin$}\par\bigskip}
\newcommand{\fineqnum}{\vspace{-.4cm$\fin$}\par\bigskip}
\newcommand{\n}[1]{\disp{\|#1\|}}

\newcommand{\pe}[2]{\disp{\langle #1,#2\rangle}}
\newcommand{\bcero}{{\bf 0}}
\newcommand{\mcero}{\text{\rm \sffamily \bfseries 0}}
\newcommand{\trace}{{\bf trace}~}

\newtheorem{coro}{\bf \sffamily Corollary}[section]
\newtheorem{theo}{\bf \sffamily Theorem}[section]
\newtheorem{prop}{\bf \sffamily Proposition}[section]
\newtheorem{rem}{\bf \sffamily Remark}[section]
\newtheorem{lemma}{\bf \sffamily Lemma}[section]

\newtheorem{exam}{\bf \sffamily Example}[section]

\begin{document}
	
\title{{\Large \bfseries\sffamily Large solutions of elliptic semilinear equations non-degenerate near the boundary}
\author{\bfseries\sffamily G. D\'{\i}az}
\date{}
}	
\maketitle
\begin{flushright}
	{\em Dedicated to Ilde, brother and master}
\end{flushright}  
\begin{abstract}
In this paper we study the so-called large solutions of elliptic semilinear equations with non null sources term, thus solutions blowing up on the boundary of the domain for which reason they are greater than any other solution whenever Weak Maximum Principle holds. The main topic about large solutions is  uniqueness results and their behavior near the boundary. It is much less than being simple. The structure of the semilinear equations considered includes the well known Keller-Osserman integral and an assumption on the ellipticity of the leading part of the differential operator. 
In our study an uniform ellipticity near the boundary is required. We consider source terms in the PDE whose boundary explosion is consistent with the Keller-Osserman condition. Extra Keller-Osserman explosions on the source are also studied, showing in particular that in some cases the PDE only admits large solutions.
\end{abstract}

\section{Introduction}
\label{sec:intro}
 This paper deals with solutions of semilinear elliptic equations with non null source term 
\begin{equation}
-\trace \cA(\cdot)\rD^{2}u+\beta (u)=f\quad \hbox{in $\Omega$},
\label{eq:fullyequation}
\end{equation}
 blowing up on the boundary, $u(x)=+\infty,~x\in\partial \Omega$. 
When a Maximum Principle holds these solutions are called {\em large solution} because they are bigger than any other bounded solution. An ingredient or the existence of large solutions is the regularity on the domains  compatible with a suitable structure of the equation, preferably near the boundary. The main topic is related with the uniquenees of large solutions. It is much less than being simple. There is a wide bibliography on large solutions, mainly for the homogeneous case $f\equiv 0$. A short list surely not the best would contain, for instance, \cite{ADR,B,Bi,DL,GM,K,LL,LeGall1,Li,LMV,O} and the references therein.
\par
The main tool near the boundary in this paper is the Lipschitz continuous distance function $\rd(x)\doteq \hbox{dist} (x,\partial \Omega),~x\in \overline{\Omega},$ defined on an open bounded set $\Omega\subset \R^{\rN},~\rN>1$. When $\partial \Omega \in\cC^{1}$ there exists a positive constant $\delta_{\Omega}$ such that $\rd\in\cC^{1}$ in the subset $\Omega_{\delta_{\Omega}}\doteq \{x\in\Omega:~\rd(x)<\delta_{\Omega}\}$ as it was proved in \cite{GT} where also it is shown the {\em eikonal equation}
\begin{equation}
|\rD \rd(x)|=1,\quad x\in\Omega_{\delta_{\Omega}},
\label{eq:eikonal}
\end{equation}
where we will denote $\bn(x)=-\rD \rd(x)$ even if $x$ is not on the boundary. If $x\in\partial \Omega$, $\bn(x)$ is just the unit outward normal to $\partial \Omega$ at $x$.
Futhermore, when $\partial \Omega \in\cC^{m},~m\ge 2,$ one has $\rd	\in\cC^{m}(\Omega_{\delta_{\Omega}})$ and
$$
\Delta \rd(x)=-(\rN-1)\rH(z_{x})+o(1),\quad x\in\Omega_{\delta_{\Omega}},
$$
where $\rH(z_{x})$ is the mean curvature of $\partial \Omega$ at the boundary point $z_{x}\in\partial \Omega$ such that $\rd(x)=|x-z_{x}|$. 
So that, we emphasize that all reasoning near the boundary $\partial \Omega$ will be considered inside $\Omega_{\delta_{\Omega}}$.
\begin{rem}\rm In some papers the global regularity of $\partial \Omega\in\cC^{2}$ can be replaced by working on regular boundary points in the sense of an interior and an exterior ball condition is fulfilled (see \cite{Ma} of \cite[Theorem 2]{ADLR}). From a different approach  the uniqueness of large solutions where the regularity of $\partial \Omega$ is slightly lowered have been studies for several authors (see, for instance, \cite{La} or \cite{LMV}).$\fin$ 
\end{rem}
\par
In this paper we will suppose that the leading term of the PDE is a Lipschitz continuous function $\cA:~\Omega\rightarrow \cS^{\rN}_{+}\setminus\{\mcero\}$ satisfying the general ellipticity condition
\begin{equation}
0\le \trace \cA (x)\xi\otimes \xi\le \Lambda,\quad x\in\Omega,~\xi\in\R^{\rN},~|\xi|=1
\label{eq:ellipticity}
\end{equation}
that becomes {\em uniform ellipticity} near the boundary  in the sense
\begin{equation}
\lambda\le \trace \cA (x)\bn(x)\otimes \bn(x),\quad x\in\Omega_{\delta_{\Omega}}
\label{eq:ellipticityneraboundaty}
\end{equation}
for some positive constants $0<\lambda \le \Lambda$. Certainly, \eqref{eq:ellipticityneraboundaty} is a kind of eikonal property (see \eqref{eq:eikonal}). So, that, we will assume in this paper that the eikonal-like term
$\rE(x)\doteq \cA (x)\bn(x)\otimes \bn(x),~x\in\Omega_{\delta_{\Omega}}$, takes values on the bounded interval $[\lambda,\Lambda]$.
As usual, $\cS^{\rN}$ is the symmetric $\rN\times\rN$ real matrices set equipped with the  ordering
$$
\rY\in\cS^{\rN}_{+} \quad \Leftrightarrow \quad \hbox{the eigenvalues of $\rY$ are nonnegative}.
$$
Since $\cA\not\equiv \bcero$ one satisfies
\begin{equation}
\cA(x)\in\cS^{\rN}_{+},~x\in\cA	\quad \Rightarrow \quad \Lambda_{\cA}\ge \lambda_{\cA}>0,
	\label{eq:positivetrace}
\end{equation}
where $\lambda_{\cA}\doteq \inf_{x\in\Omega}\trace \cA(x)\rI$ and $\Lambda_{\cA}\doteq \sup_{x\in\Omega}\trace \cA(x)\rI$. We also will use the notation $\rC_{\partial \Omega,\cA}\doteq \n{\trace \cA(\cdot)\rD^{2}\rd}_{ \rL^{\infty}(\Omega_{\delta_{\Omega}})}<+\infty$.
\begin{rem}\rm Involving to stochastic process it is usual the case
$$
\cA(x)=\sigma(x)\sigma^{\tt t}(x),\quad x\in\Omega,
$$
where $\sigma :~\Omega\rightarrow \R^{\rN}\times\R^{p},~p\le \rN$, is a Lipschitz continuous function. Clearly \eqref{eq:ellipticity} holds by construction and then \eqref{eq:ellipticityneraboundaty} only can be available when $p=\rN$. If $\cA\in\rW^{2,\infty}(\Omega)$ then $\sqrt{\cA}$ is locally Lipschitz continuous in $\Omega$, as it was proved in \cite{Fr}.$\fin$
\label{rem:elipticity}
\end{rem}
In Section \ref{sec:SMP} we will prove a Strong Maximum Principle (see Theorem \ref{theo:SMPwithoutcoercivity} below) assuming that for any $\delta>0$ there exists $\lambda(\delta)>0$ such that
\begin{equation}
\trace \cA(x)\bp\otimes \bp\ge \lambda (\delta),\quad x\in\Omega,~\delta <|\bp|<\delta ^{-1}.
\label{eq:lightellipticity}
\end{equation}
as well as that for any $\rR>0$ and $x\in\Omega$ there exist $\alpha >0$  and $\delta>0$ such that for all $\bp(x)\in\bB_{\delta}(\bcero)\setminus \{\bcero\}$  one satisfies
\begin{equation}
\trace \cA(x) \bp(x)\otimes \bp(x)\ge \dfrac{2\Lambda_{\cA}}{\alpha \rR^{2}}|\bp(x)|^{2},\quad x\in\Omega.
\label{eq:structureannulusSMPsub}
\end{equation}
\begin{rem}\rm Under an uniformly elliptic condition
$$
\trace \cA(x)\bp\otimes \bp\ge \lambda>0,\quad |\bp|=1\
$$
all properties \eqref{eq:ellipticity}, \eqref{eq:ellipticityneraboundaty}, \eqref{eq:lightellipticity} and \eqref{eq:structureannulusSMPsub} hold.$\fin$
\end{rem}
\par
We also consider zeroth order terms given by continuous and increasing function $\beta:~\R_{+}\rightarrow \R_{+}$, with $\beta (0)=0$ by simplicity.  Certainly, some additional assumptions on the data must be required in order to prove that the existence of classical large solutions. An almost ``unavoidable'' hypothesis for these goal is the {\em non-degeneracy} of the operator $\cA$. In general, these condition will not be required in this paper, as it ocurrs in some important examples of the applications and consequently the $\cC^{2}$ regularity  even the continuity of the solutions can not be guaranteed (see \cite{Fr} or \cite{OleRa}).  The viscosity solution notion is adequate in order to remove {\em non-degeneracy} hypothesis. We send to the monograph by M.G. Crandall, H. Ishii and P.L. Lions \cite{CIL} to understand how semi-continuous functions can solve \eqref{eq:fullyequation} (see Section \ref{sec:existence}).  The term {\em viscosity solutions} is an artifact motivated by their 
consistence with the {\em vanishing viscosity method}. Since this viscosity notation is the primary ones sometimes we drop the term {\em viscosity} and hereafter simply refer to solution.
\par 
Some other structural assumptions are required when we consider the differential operator
$$
\cF(\rX,t,x)\doteq-\trace \big (\cA (x)\rX\big )+\beta (t),\quad (\rX,x,t)\in\cS^{\rN}\times \Omega\times \R_{+}.
$$
When the boundary of $\Omega$ is smooth enough, the existence of  large solutions is equivalent to the so-called Keller-Osserman condition
\begin{equation}
\int^{+\infty}\dfrac{ds}{\sqrt{\rG(s)}}<+\infty \quad \hbox{where $\disp \rG(s)=\int^{s}_{0}\beta (r)dr$}.
\label{eq:KellerOsserman}
\end{equation}
We notice that under \eqref{eq:KellerOsserman} we may define implicitly the decreasing function
\begin{equation}
\delta =\Phi\big (\phi(\delta )\big)\quad \hbox{where $\disp \Phi(t)=\int^{+\infty}_{t}\dfrac{ds}{\sqrt{2\rG(s)}},\quad t>0$}
\label{eq:profilephi}
\end{equation}
that satisfies
\begin{equation}
\phi(0)=+\infty\quad \hbox{and}\quad \phi''(\delta )=\beta \big (\phi (\delta )\big )\quad\hbox{for small $\delta >0$.}
\label{eq:ODEKO}
\end{equation}
Sometimes we say that $\phi$ provides the Keller-Osserman explosive profile.
\begin{rem}\rm In fact, by strightforward computations one proves that \eqref{eq:KellerOsserman} implies
$$
\int^{+\infty}\dfrac{ds}{\sqrt{\rG(s)-\rG(a)}}<+\infty 
$$
for any $a\ge 0$. Then one defines the function $\phi_{a}(\delta)$ by
$$
\delta =\Phi_{a}\big (\phi_{a}(\delta )\big)\quad \hbox{where $\disp \Phi_{a}(t)=\int^{+\infty}_{t}\dfrac{ds}{\sqrt{2\big (\rG(s)-\rG(a)\big )}},\quad t>0$}
$$
that also satisfies \eqref{eq:ODEKO}. Moreover, one proves
\begin{equation}
\lim_{\delta\searrow 0}\big (\phi_{a_{1}}(\delta)-\phi_{a_{2}}(\delta))=0
\label{eq:uniquenessprofileprimitive}
\end{equation}
(see \cite[Theorem 1]{BM}).$\fin$
\label{rem:primitive}
\end{rem}
As it will be proved later, the function $\phi$ governs boundary profiles of the large solutions independent on the primitive of $\beta$ (see Remark \ref{rem:primitive}). 
\begin{rem}[Early examples]\rm In certain differential equations associated with the Riemannian metrics appear large solutions relative to $\Lambda=\lambda=1,~f\equiv 0$ and  $\beta (t)=t^{\frac{\rN+2}{\rN-2}}$. It was studied in \cite{LN}. In general, for the power like case $\beta _{m}(t)=t^{m}$ the condition \eqref{eq:KellerOsserman} becomes $m>1$ and
$$
\Phi_{m}(t)=\sqrt{\dfrac{m+1}{2}}\dfrac{2}{m-1}t^{-\frac{m-1}{2}}
\quad \hbox{and}\quad
\phi_{m}(\delta)=\left (\dfrac{(m-1)^{2}}{2(m+1)}\right )^{ -\frac{1}{m-1}}\delta^{-\frac{2}{m-1}}.
$$
The study of large solutions trace back to \cite{Bi} for $\Lambda=\lambda=1,~f\equiv 0$ and  $\beta (t)=e^{t}$ that satisfies \eqref{eq:KellerOsserman} for which 
$$
\phi(\delta )=\log \left (\dfrac{2}{\delta^{2}}\right ).
$$
The upcoming choice  $\beta (t)=te^{2t}$ also satisfies \eqref{eq:KellerOsserman}. Then
$$
\phi(\delta)=\sqrt{2}\hbox{ erfc}^{-1}\left(\dfrac{\delta}{\sqrt{\pi}}
\right),
$$
where $\disp \hbox{erfc}(\delta)=1-\hbox{erf}(\delta)=\dfrac{2}
{\sqrt{\pi}}\int^{\infty}_{\delta}e^{-s^{2}}ds$.
\fineq
\label{rem:exemplesbeta}
\end{rem}
\begin{rem}\rm From the inequality
$$
\widehat{\beta }(t)\ge \beta (t)\quad \hbox{for large $t$},
$$
it follows that if $\beta $ verifies  \eqref{eq:KellerOsserman} the same happens with $\widehat{\beta }$. In particular, any function $\widehat{\beta}(t)\ge t q(t)$ such that 
$$
\liminf_{t \rightarrow \infty} \dfrac{q(t)}{t^{\gamma}}\in (0,+\infty]\quad 
\hbox{for some $\gamma >0$}
$$
satisfies \eqref{eq:KellerOsserman}. For instance, we may choose  $q(t)\ge (\log t)^{\gamma},~\gamma\ge 1,$ or $q(t)\ge \log (\log (\cdots \log (t)\cdots )).\fin$
\label{rem:exemplesmasg}
\end{rem}
\par
The plan in this paper is as follows. In Sections \ref{sec:KOgrowth} and  \ref{sec:KOgrowthextra} we obtain a kind of comparison of the boundary behavior by means of the property
\begin{equation}
\limsup_{\rd(x)\rightarrow 0}\dfrac{u^{*}(x)}{v_{*}(x)}\le 1,
\label{eq:UVontheboundary}
\end{equation}
for any subsolution $u$ and any nonnegative large supersolution $v$  of  
$$
-\trace \cA(x)\rD^{2}u+\beta (u)=f\quad \hbox{in $\Omega$}.
$$
Since the large solutions are singular on the boundary, the condition \eqref{eq:UVontheboundary} is very useful when we try to transfer boundary comparison to the interior.  We note that \eqref{eq:UVontheboundary} avoids working with the {\em relaxed limit condition} (see \cite{CIL}). First, in Theorems \ref{theo:uniquenessboundary} and  \ref{theo:absoluteboundarycomparison}, in Section \ref{sec:KOgrowth}, the property \eqref{eq:UVontheboundary} is obtained under sources  terms, $f$, with Keller-Osserman type growth on the boundary, including explosive and bounded boundary behaviour of $f$. 
\par
In Section \ref{sec:KOgrowthextra} we also obtain \eqref{eq:UVontheboundary}
now for any subsolution $u$ and any large supersolution $v$  of  
$$
-\trace \cA(x)\rD^{2}u+\beta (u)=f\quad \hbox{in $\Omega$},
$$
under sources  terms, $f$, with extra Keller-Osserman type growth on the boundary. 
Here, the boundary behaviour of $f$ only can be explosive. It is collected in Theorem \ref{theo:comparisonboundaryextrageneral}.  
It is about a genuine feature of the non homogeneous equation
$$
-\trace \cA(\cdot)\rD^{2}u+\beta(u)=f\quad \hbox{in $\Omega$}.
$$
The key idea of Theorem \ref{theo:comparisonboundaryextrageneral},  announced in \cite{DL} and \cite{ADLR}, is the scheme $ \beta(u)\approx f$ for large values of $u$ near the boundary while the leading part of the equation grows up with a lower order. Certainly the source term $f$ must be more explosive on the boundary than  the studied in the Keller-Osserman type growth case.  Essentially we require a compatiblity condition as
\begin{equation}
f(x)=h\big (\rd(x)\big )+o(\rd(x))\quad \hbox{with}\quad h(\delta)> \beta \big (\phi (\delta) \big )\quad \hbox{for $\rd(x)$ small} 
\label{eq:compatibility}
\end{equation}
provided \eqref{eq:KellerOsserman}. We recall that, as it was proved in \cite{K}, the assumption \eqref{eq:KellerOsserman} is a necessary condition for the existence of large solutions. In Theorem \ref{theo:comparisonboundaryextrageneral} we will obtain the precise  lower order estimate term
$$
u(x)=\beta^{-1} \big (h(\rd(x))\big )+o(\rd(x))
$$
near the boundary $\partial \Omega$. Here the boundedness of the eikonal like term plays a secondary role.  Since the proof of Theorem \ref{theo:comparisonboundaryextrageneral} is tedious, we illustrate the reasoning for the power like case in Theorem~\ref{theo:comparisonboundaryextrapower} in which a striking case appears  
by means of an argument of rescaling: the explositivity of the sources is transferred  to the supersolutions as it is proved in Theorem \ref{theo:explosivitytransferred}. Consequently, in this case the equation
$$
-\trace \cA(\cdot)\rD^{2}u+u^{m}=f\quad \hbox{in $\Omega$}
$$
only can admit large solutions.

We devote Section \ref{sec:existence} to adequate versions of a Weak  Maximum Principle. More precisely, we obtain comparison results  involved large solution as it is proved in Theorem \ref{theo:comparisonprincipleviscosity} under the useful estimate  \eqref{eq:UVontheboundary}.  So, by Theorem \ref{theo:comparisonprincipleviscosity} we  transfer the comparison given by \eqref{eq:UVontheboundary} to the interior 
$$
u^{*}(x)\le v_{*}(x),\quad x\in\Omega,
$$
provided
\begin{equation}
\dfrac{\beta (t)}{t}\mbox{ is increasing for }  t>0.
\label{eq:crecimientog1}
\end{equation}
Therefore
$$
-\trace \cA(\cdot)\rD^{2}u+\beta (u)=f\quad \hbox{in $\Omega$} 
$$
admits at most a unique large solution $u$. Since it implies
$$
u^{*}(x)\le u_{*}(x),\quad x\in\Omega
$$
the solution must be a continuous function. The proof adapts reasoning of the \cite[Theorem 3.3]{CIL} of the Theory of Viscosity Solutions whose notions are briefly recalled. We emphasize that no uniform ellipticity assumptions is required in the proof of Theorem \ref{theo:comparisonprincipleviscosity}.  
\par
The existence of solutions is obtained in Theorems 
\ref{theo:existenceKOsources} and \ref{theo:existenceextraKOsources} by means of the classical Perron's Method (see \cite{CIL} and the bibliography therein). As it is well konwn the main ingredient of this method concerns to weak pass to the limit in equations satisfying the Maximum Principle. In some sense, the stability property of solutions can be seen as the analogue Minty's device for monotone operators. A main feature of the Viscosity Theory is the pass to limit of semi-continuous solutions. Here we deal with it by appropriate universal estimates obtained in Theorems \ref{theo:universalestimateepsilonn} and \ref{theo:universalestimateepsilonn}. They are simple consequence of the more general Theorem  \ref{theo:universalestimate} proved in Section \ref{sec:universalbound}.  We prove in Theorems \ref{theo:existenceKOsources} and \ref{theo:existenceextraKOsources} that Perron function is in fact the unique solution. 
\par
Section \ref{sec:secondorder} is devoted to the boundary explosive expansion of the large solutions. Since  it requires tedious and dark computations we illustrate by considering only the second order term of the expansion. We emphasize with the influence of the geometry provided of the secondary term that in some cases is explosive, bounded or even vanishig on the boundary. While the case of Keller-Osserman type growth on the boundary of the sources is studied in a general way,  for extra Keller-Osserman type growth we only consider the power like case.
\par
Two last sections are included in this article collecting some results previously used. So, in Section \ref{sec:universalbound} we take advantage of the Keller-Osserman assumption \eqref{eq:KellerOsserman} in order to obtain a universal upper bound very useful in locally pass to the limits (see Theorem  \ref{theo:universalestimate}). On the other hand, some results without coercivity term
have been required. In particular,  in the proof of Theorem \ref{theo:absoluteboundarycomparison} we have used a Weak Maximum Principle without coercivity term (see Theorem \ref{theo:withoutcoercivityI}).
Moreover, the uniqueness of Perron function as the unique large solution, is proved by means of a Strong Maximum Principle without coercivity term (see Theorem \ref{theo:SMPwithoutcoercivity}) under the assumptions  \eqref{eq:lightellipticity} and \eqref{eq:structureannulusSMPsub}. Here, the lack of a coercive term is supplied with a kind of Kruzkov change of variable (see \cite{BBus} or \cite{BDD}).
Theorems \ref{theo:withoutcoercivityI} and \ref{theo:SMPwithoutcoercivity} are been obtained in Section \ref{sec:SMP}. 
\par
Some results of the paper admit extension when the function $\beta(r)$ is replaced by suitable terms $\beta(x,r)$ as it is studied in \cite{GM} or \cite{LMV}.

\setcounter{equation}{0}

\section{Keller-Osserman type growth of the sources on the boundary}
\label{sec:KOgrowth}

As it was pointed out,  the Keller-Osserman condition \eqref{eq:KellerOsserman} 
$$
\int^{+\infty}\dfrac{ds}{\sqrt{\rG(s)}}<+\infty \quad \hbox{where $\disp \rG(s)=\int^{s}_{0}\beta (r)dr$}
$$
is equivalent to the existence of  large solutions (see \cite{K}). The uniqueness of large solutions turns out to be delicate. In this section we study a kind of uniqueness of the boundary behaviours of the large solutions. We begin with some technicallities.
\begin{lemma}[Lemma A.1 of \cite{ADR}]\rm Assumption \eqref{eq:KellerOsserman} implies 
\begin{equation}
\lim_{t\rightarrow +\infty}\dfrac{t}{\sqrt{\rG(t)}}=0,
\quad \lim_{t\rightarrow +\infty}\dfrac{t}{\beta (t)}=0
\label{eq:fromKellerOsserman2}
\end{equation}
and
\begin{equation}
\lim_{t\rightarrow  +\infty}\dfrac{\sqrt{\rG(t)}}{\beta (t)}=0.
\label{eq:fromKellerOsserman}
\end{equation}
\fineqnum
\end{lemma}
\par
In some cases we replace \eqref{eq:KellerOsserman}  by 
\begin{equation}
\dfrac{\beta (t)}{t^{\alpha}}\quad	\hbox{is increasing for large $t$},
\label{eq:crecimientog}
\end{equation}
where $\alpha>1$, slightly more restrictive. 
\begin{lemma} Under \eqref{eq:crecimientog} one has the Keller-Osserman condition \eqref{eq:KellerOsserman}. Moreover,  for $\nu>1$ and $\varepsilon >0$ the inequality
\begin{equation}
\varepsilon +\nu \phi (\delta)\ge \phi\big ( \nu^{-\frac{\alpha-1}{2}}\delta\big )\quad \hbox{for $\delta $ small}
\label{eq:casinopotencias}
\end{equation}
holds. 
\label{lemma:KOrestrictive} 
\end{lemma}
\proof  Let us assume 
$$
\dfrac{\beta (t)}{t^{\alpha}}\quad	\hbox{is increasing for $t>a$}.
$$
Then from Remark  \ref{rem:exemplesmasg} it follows that 	\eqref{eq:crecimientog} leads to the Keller-Osserman condition \eqref{eq:KellerOsserman}. Futhermore,  it  implies 
$$
\rG(\nu t)-\rG(a)=\int^{\nu t}_{a}\beta (s)ds=\nu\int ^{t}_{\frac{a}{\nu}}\beta (\nu s)ds \ge \nu^{\alpha +1}\big (\rG(t)-\rG(a)\big )\quad \hbox{ for $\nu>1$ },
$$
hence we deduce
$$
\Phi _{a}(\nu t)=\int^{+\infty}_{\nu t}\dfrac{d\tau}{\sqrt{2\big (\rG(\tau)-\rG(a)\big )}}=\nu
\int^{+\infty}_{t}\dfrac{d\tau}{\sqrt{2\big (\rG(\nu\tau)-\rG(a)\big )}}\le
\nu^{\frac{1-\alpha}{2}}\Phi _{a}(t),\quad \hbox{for large $t$},
$$
thus $\nu t\ge \Phi_{a}^{-1}\big ( \nu^{\frac{1-\alpha}{2}}\Phi_{a}(t)\big )$.  The change of variable $t=\phi_{a}(\delta)$ proves
$$
\nu \phi_{a}(\delta)\ge \phi_{a}\big ( \nu^{-\frac{\alpha-1}{2}}\delta\big )\quad \hbox{for $\delta $ small}.
$$
Finally, \eqref{eq:casinopotencias} follows from \eqref{eq:uniquenessprofileprimitive}.$\fin$
\begin{rem}\rm The above reasonings prove
for $\nu<1$ and $\varepsilon >0$ the inequality
\begin{equation}
-\varepsilon +\nu \phi (\delta)\le \phi\big ( \nu^{-\frac{\alpha-1}{2}}\delta\big )\quad \hbox{for $\delta $ small}.
\label{eq:casinopotenciasmenor}
\end{equation}
Clearly \eqref{eq:crecimientog} is verified for the examples of Remark \ref{rem:exemplesbeta}. $\fin$
\end{rem}
\par
\medskip
The key idea of the boundary explosive estimate is the simple  equality 
$$
\phi ''(\delta)=\beta \big (\phi (\delta)\big )=\mu \beta \big (\phi (\delta)\big )-(\mu-1)\beta \big (\phi (\delta)\big )
$$
(see \eqref{eq:ODEKO}). The proof is based on following boundary behaviour
\begin{prop}[Boundary maximal universal upper bound] Asume $\partial \Omega \in \cC^{2}$. Suppose \eqref{eq:ellipticity}, \eqref{eq:ellipticityneraboundaty} and \eqref{eq:crecimientog}. Then any solution $u$ 
$$
-\trace \cA(x)\rD^{2}u+ \beta (u)\le f\quad \hbox{in $\Omega$}
$$
admits the maximal boundary behaviour given by
\begin{equation}
\limsup_{\rd(x)\searrow 0}\dfrac{u^{*}(x)}{\phi \left  (\sqrt{\dfrac{1-\ell_{u}}{\Lambda}}\rd(x)\right )}\le 1,
\label{eq:maximalboundaryprofile}
\end{equation}
provided 
\begin{equation}
\limsup_{\rd(x)\searrow 0}\dfrac{f(x)}{\beta\left  (\phi \left  (\sqrt{\dfrac{1-\ell_{u}}{\Lambda}}\rd(x)\right )\right )}\le \ell_{u} < 1,
\label{eq:criticalgrowthd}
\end{equation}
for some positive constant  $\ell_{u}$. Here $\phi$ is the function defined in \eqref{eq:profilephi}.
\label{prop:maximalboundaryprofile}
\end{prop}
\proof The idea is to construct a suitable classical large supersolution.  Let  $\varepsilon>0$. So,  for $0<\delta <\delta_{1}<\delta_{\Omega}$, we consider
$$
\rU_{0}(x)=\phi  _{a}\left (\sqrt{\dfrac{(1-\varepsilon)(1-\ell_{u} )}{\Lambda}}(\rd(x)-\delta)\right  ), \quad x\in \Omega^{\delta}_{\delta_{1}}\doteq \{x\in\Omega:~\delta <\rd(x)<\delta_{1}\}.
$$ 
Here $\phi$ is the decreasing function defined implicitly in \eqref{eq:profilephi}. In the pioneer \cite{LN} this kind of method was used. By \eqref{eq:ODEKO}, straightforward computations as in \cite{ADR} or \cite{DL}  prove
$$
\begin{array}{ll}
\trace \cA(x)\rD ^{2}\rU_{0}(x)&\hspace*{-.2cm}=\dfrac{(1-\varepsilon)(1-\ell_{u} )}{\Lambda}\phi ''(\sigma)\trace \cA(x)\rD \rd(x)\otimes \rD \rd(x)\\ [.25cm]
&+\sqrt{\dfrac{(1-\varepsilon)(1-\ell_{u} )}{\Lambda}}\phi '(\sigma)\trace \cA(x)\rD^{2}\rd(x)\\ [.25cm]
&\hspace*{-.2cm}\le (1-\varepsilon)(1-\ell_{u} )g\big (\phi(\sigma)\big ) +\rC_{\partial \Omega,\cA}\sqrt{\dfrac{2(1-\varepsilon)(1-\ell _{u} )}{\Lambda}\rG\big (\phi(\sigma )\big )},\quad x\in \Omega^{\delta}_{\delta_{1}},
\end{array}
$$
where we are denoting $\sigma=\sqrt{\dfrac{(1-\varepsilon)(1-\ell _{u} )}{\Lambda}}(\rd(x)-\delta)$.
Then 
$$
\begin{array}{ll}
-\trace \cA(x)\rD ^{2}\rU_{0}(x)+\beta \big (\rU_{0}(x)\big )&\hspace*{-.2cm }\ge \beta \big (\phi(\sigma )\big )\bigg [\varepsilon (1-\ell _{u})+\ell _{u} \\
& 
\left .-\rC_{\partial \Omega,\cA}\sqrt{\dfrac{2(1-\varepsilon)(1-\ell _{u} )}{\Lambda}}\dfrac{\sqrt{\rG\big (\phi(\sigma)\big )}}{\beta \big (\phi(\sigma)\big )}\right ], \quad x\in \Omega^{\delta}_{\delta_{1}}. 
\end{array}
$$
On the other hand, assumption \eqref{eq:criticalgrowthd} implies
$$
\ell _{u}\beta \big (\phi (\sigma)\big )\ge  \ell_{u} \beta\left (\phi\left (\sqrt{\dfrac{1-\ell_{u}}{\Lambda}}\rd(x)\right )\right ) \ge f(x),\quad  x\in \Omega^{\delta}_{\delta_{1}},
$$
whence
\begin{equation}
\begin{array}{ll}
-\trace \cA(x)\rD ^{2}\rU_{0}(x)+\beta \big (\rU_{0}(x)\big )-f(x)&\hspace*{-.2cm }\ge \beta \big (\phi(\sigma )\big )\bigg [\varepsilon (1-\ell_{u})\\
&\left . 
-\rC_{\partial \Omega,\cA}\sqrt{\dfrac{2(1-\varepsilon)(1-\ell _{u} )}{\Lambda}}\dfrac{\sqrt{\rG\big (\phi(\sigma)\big )}}{\beta\big (\phi(\sigma)\big )}\right ],
\end{array}
\label{eq:stepfirst}
\end{equation}
for $x\in \Omega^{\delta}_{\delta_{1}}$. Since $0<\sigma<\sqrt{\dfrac{1-\ell_{u}}{\Lambda}}\delta_{1}$, 
taking account into \eqref{eq:fromKellerOsserman} we obtain
$$
\rC_{\partial \Omega,\cA}\sqrt{\dfrac{2(1-\varepsilon)(1-\ell_{u})}{\Lambda}}\dfrac{\sqrt{\rG\big (\phi(\sigma)\big )}}{\beta \big (\phi(\sigma)\big )}<\varepsilon(1-\ell_{u}),
$$
for some $\delta_{1}$ small enough. Because $\beta \big (\phi(\sigma )\big )>0$ we deduce
$$
-\trace \cA(x)\rD ^{2}\rU_{0}(x)+\beta \big (\rU_{0}(x)\big )\ge f(x)\quad x\in \Omega^{\delta}_{\delta_{1}}.
$$
Clearly, $\rU_{0}$ is a classical positive supersolution in $\Omega^{\delta}_{\delta_{1}}$ such that   
$$
u^{*}(x)<+\infty=\rU_{0}(x)\quad \hbox{if $\rd(x)=\delta$}
$$ 
for any subsolution $u$. Moreover, 
$$
\left \{
\begin{array}{l}
-\trace \cA(x)\rD ^{2}\big (\rU_{0}(x)+\rM\big )+\beta \big (\rU_{0}(x)+\rM\big )\ge f(x),\quad x\in \Omega^{\delta}_{\delta_{1}}, \\ [.175cm]
\rU_{0}(x)+\rM=+\infty>u^{*}(x)\quad \hbox{if $\rd(x)=\delta$},
\end{array}
\right .
$$
for any positive constant $\rM$.  In order to conclude \eqref{eq:maximalboundaryprofile}  we need an adequate inequality on the boundary $\rd(x)=\delta_{1}$. In particular
$$
\rU_{0}(x)+\rM\ge \rM \ge \sup_{\rd(x)=\delta_{1}}u^{*}\ge u^{*}(x),\quad \rd(x)=\delta_{1} 
$$
for $\rM\ge \disp \sup_{\rd(x)=\delta_{1}}u^{*}$ (we recall that $u^{*}$ is an upper semi-continuous function). Remark \ref{rem:supverification}  implies
$$
u^{*}(x)\le \rU_{0}(x)+\rM,\quad x\in \Omega_{\delta_{1}}^{\delta},
$$
for all $\delta<\delta_{\Omega}$, whence
\begin{equation}
u^{*}(x)\le \phi  \left (\sqrt{\dfrac{(1-\varepsilon)(1-\ell_{u})}{\Lambda}}\rd(x)\right  ) +\rM,\quad x\in \Omega_{\delta_{1}}.
\label{eq:boundarybehaviouru}
\end{equation}
This maximal bound \eqref{eq:boundarybehaviouru} can be improved by means of  \eqref{eq:casinopotencias}. Indeed, let us choose the constant $\nu=(1-\varepsilon)^{-\frac{1}{\alpha-1}}>1$ for which
$$
u^{*}(x)\le \nu \phi  \left (\sqrt{\dfrac{1-\ell_{u}}{\Lambda}}\rd(x)\right  ) +\rM+\varepsilon,\quad x\in \Omega_{\delta_{1}}
$$
(see \eqref{eq:casinopotencias})
and 
$$
\dfrac{u^{*}(x)}{\phi  \left (\sqrt{\dfrac{1-\ell_{u}}{\Lambda}}\rd(x)\right  )}
\le \nu 
+\dfrac{\rM+\varepsilon}{\phi  \left (\sqrt{\dfrac{1-\ell_{u}}{\Lambda}}\rd(x)\right  )}
$$
hold for $x\in \Omega_{\delta_{1}}$. Then
$$
\limsup_{\rd(x)\searrow 0}\dfrac{u^{*}(x)}{\phi  \left (\sqrt{\dfrac{1-\ell_{u}}{\Lambda}}\rd(x)\right  )} \le \nu.
$$
Taking $\varepsilon \searrow 0$ one concludes \eqref{eq:maximalboundaryprofile}.$\fin$
\par
\medskip
\noindent We emphasize that in Theorem \ref{theo:absoluteboundarycomparison} we do not require any information of $u$ on the $\partial \Omega$.
\begin{prop}[Boundary minimal universal lower bound]  Asume $\partial \Omega \in \cC^{2}$. Suppose \eqref{eq:ellipticity}, \eqref{eq:ellipticityneraboundaty} and \eqref{eq:crecimientog}.  Then any  nonnegative large solution $v$ of
$$
-\trace \cA(x)\rD^{2}v+\beta (v)\ge f\quad \hbox{in $\Omega$}
$$
admits the minimal behaviour given by
\begin{equation}
1\le \liminf_{\rd(x)\searrow 0}\dfrac{v_{*}(x)}{\phi \left (\sqrt{\dfrac{1-\ell_{d}}{\lambda}}\rd(x)\right )}
\label{eq:minimalboundaryprofileg}
\end{equation}
provided 
\begin{equation}
\liminf_{\rd(x)\searrow 0}\dfrac{f(x)}{\beta\left (\phi \left (\sqrt{\dfrac{1-\ell_{d}}{\lambda}}\rd(x)\right )\right )}\ge \ell_{d} \quad\hbox{with $0\le \ell_{d} <1$},
\label{eq:criticalgrowthu}
\end{equation}
for some positive constant $\ell_{d}$. Here $\phi$ is the function defined in \eqref{eq:profilephi}.
\label{prop:minimalboundaryprofile}
\end{prop}
\proof For $\varepsilon >0,~\delta >0$ and $0<\delta_{0}<\delta_{\Omega}$, we construct
$$
\rV_{0}(x)=\phi  \left (\sqrt{\dfrac{(1+\varepsilon)(1-\ell_{d})}{\lambda}}\big (\rd(x)+\delta\big )\right  ), \quad x\in \Omega_{\delta_{0}}.
$$ 
Again $\phi$ is the decreasing function defined implicitly in \eqref{eq:profilephi}. Now, straightforward computations prove
$$
\begin{array}{ll}
\trace \cA(x)\rD ^{2}\rV_{0}(x)&\hspace*{-.2cm}=\dfrac{(1+\varepsilon)(1-\ell_{d})}{\lambda}\phi ''(\sigma)\trace \cA(x)\rD \rd(x)\otimes \rD \rd(x)\\ [.25cm]
& +\sqrt{\dfrac{(1+\varepsilon)(1-\ell_{d})}{\lambda}}\phi'(\sigma)\trace \cA(x)\rD^{2}\rd(x)\\ [.25cm]
&\hspace*{-.2cm}\ge (1+\varepsilon)(1-\ell_{d})g\big (\phi(\sigma)\big ) -\rC_{\partial \Omega,\cA}\sqrt{\dfrac{2(1+\varepsilon)(1-\ell_{d})}{\lambda}\rG\big (\phi(\sigma )\big )},\quad x\in \Omega_{\delta_{0}},
\end{array}
$$
where we are denoting $\sigma=\sqrt{\dfrac{(1+\varepsilon)(1-\ell_{d})}{\lambda}}\big (\rd(x)+\delta)$.
Then 
$$
\begin{array}{ll}
-\trace \cA(x)\rD ^{2}\rV_{0}(x)+\beta \big (\rV_{0}(x)\big )&\hspace*{-.2cm }\le \beta \big (\phi(\sigma )\big )\bigg [-\varepsilon(1-\ell_{d} )+\ell _{d}\\
&\left . 
+\rC_{\partial \Omega,\cA}\sqrt{\dfrac{2(1+\varepsilon)(1-\ell_{d})}{\lambda}}\dfrac{\sqrt{\rG\big (\phi(\sigma)\big )}}{\beta\big (\phi(\sigma)\big )}\right ],\quad x\in \Omega_{\delta_{0}}.
\end{array}
$$
On the other hand, assumption \eqref{eq:criticalgrowthu} implies
$$
\ell _{d}\beta \big (\phi(\sigma)\big )\le  \ell _{d}\beta \left (\phi\left (\sqrt{\dfrac{(1+\varepsilon)(1-\ell_{d})}{\lambda}}\rd(x) \right )\right ) \le f(x),
\quad x\in \Omega_{\delta_{0}}
$$
whence
$$
\begin{array}{ll}
-\trace \cA(x)\rD ^{2}\rV_{0}(x)+\beta \big (\rV_{0}(x)\big )-f(x)&\hspace*{-.2cm}\le \beta \big (\phi(\sigma )\big )\bigg [-\varepsilon(1-\ell _{d})\\
& \left .
+\rC_{\partial \Omega,\cA}\sqrt{\dfrac{2(1+\varepsilon)(1-\ell_{d})}{\lambda}}\dfrac{\sqrt{\rG\big (\phi(\sigma)\big )}}{\beta\big (\phi(\sigma)\big )}\right ],\quad x\in \Omega_{\delta_{0}}.
\end{array}
$$
 We also will use the property \eqref{eq:fromKellerOsserman}. 
So, fixed $\varepsilon _{0}>0$  and $\varepsilon <\varepsilon _{0} $ we choose $0<\delta\le \delta_{0}$ for which
$0<\sigma\le 2\sqrt{\dfrac{(1+\varepsilon_{0})(1-\ell_{d})}{\lambda}}\delta_{0}$ and the inequality 
$$
\rC_{\partial \Omega,\cA}\sqrt{\dfrac{2(1+\varepsilon_{0})(1-\ell_{d})}{\lambda}}\dfrac{\sqrt{\rG\big (\phi(\sigma)\big )}}{\beta\big (\phi(\sigma)\big )}<\varepsilon
$$
holds for $\delta_{0}$ small enough. It implies that the function
$$
\rV_{0}(x)=\phi  \left (\sqrt{\dfrac{(1+\varepsilon)(1-\ell_{d})}{\lambda}}\big (\rd(x)+\delta\big )\right  ), \quad x\in \Omega_{\delta_{0}},
$$ 
satisfies  
$$
-\trace \cA(x)\rD ^{2}\rV_{0}(x)+\beta \big (\rV_{0}(x)\big )\le f(x)\quad x\in \Omega _{\delta_{0}}.
$$
Clearly, $\rV_{0}$ is a classical bounded subsolution in $\Omega _{\delta_{0}}$ such that   
$$
\rV_{0}(x)<+\infty =v_{*}(x)\quad \hbox{if $\rd(x)=0$}
$$ 
for any large supersolution $v$.  So, for any nonnegative constant $\rM$ one has
$$
\left \{
\begin{array}{l}
-\trace \cA(x)\rD ^{2}\big (\rV_{0}(x)-\rM \big )+\beta \big (\rV_{0}(x)-\rM \big )\le f(x),\quad x\in \Omega _{\delta_{0}}, \\ [.175cm]
\rV_{0}(x)-\rM <+\infty=v_{*}(x)\quad \hbox{if $\rd(x)=0$}.
\end{array}
\right .
$$
In order to conclude \eqref{eq:minimalboundaryprofileg} 	 we need an adequate inequality on the boundary $\rd(x)=\delta_{0}$. In particular
$$
\rV_{0}(x)-\rM\le \phi \left (2\delta_{0}\sqrt{\dfrac{(1+\varepsilon)(1-\ell_{d})}{\lambda}}\right )-\rM \le 0 \le v_{*}(x),\quad \rd(x)=\delta_{0} 
$$
for $\rM\ge  \phi \left (2\delta_{0}\sqrt{\dfrac{(1+\varepsilon)(1-\ell_{d})}{\lambda}}\right )$, when we are taking $\delta=\delta_{0}$. Remark \ref{rem:subverification} implies
$$
\rV_{0}(x)-\rM \le v_{*}(x),\quad x\in \Omega_{\delta_{0}},
$$
thus
\begin{equation}
\phi \left (\sqrt{\dfrac{(1+\varepsilon)(1-\ell_{d})}{\lambda}}\big ( \rd(x)+\delta_{0}\big )\right ) - \rM\le v_{*}(x),\quad x\in \Omega_{\delta_{0}}
\label{eq:boundarybehaviourv}
\end{equation}
holds. Again we may improve the minimal bound \eqref{eq:boundarybehaviourv} by means of  \eqref{eq:casinopotenciasmenor} here by choosing
$\nu=(1+\varepsilon)^{-\frac{1}{\alpha-1}}<1$.
Indeed,  \eqref{eq:boundarybehaviourv} becomes
$$
\nu \phi \left (\sqrt{\dfrac{1-\ell_{d}}{\lambda}}\big ( \rd(x)+\delta_{0}\big )\right ) - (\rM+\varepsilon)\le v_{*}(x),\quad x\in \Omega_{\delta_{0}}.
$$
It leads to
$$
\nu-\dfrac{\rM+\varepsilon}{\phi  \left (\sqrt{\dfrac{1-\ell_{d}}{\lambda}}\rd(x)\big )\right )}\le \dfrac{v_{*}(x)}{\phi  \left (\sqrt{\dfrac{1-\ell_{d}}{\lambda}}\rd(x)\big )\right )},\quad x\in \Omega_{\delta_{0}}
$$
and
$$
\nu\le \liminf_{\rd(x)\searrow 0}\dfrac{v_{*}(x)}{\phi  \left (\sqrt{\dfrac{1-\ell_{d}}{\lambda}}\rd(x)\big )\right )}.
$$
Taking $\varepsilon \searrow 0$ one concludes \eqref{eq:minimalboundaryprofileg}.$\fin$
\begin{rem}\rm
The property \eqref{eq:casinopotencias} avoids to consider the assumption as
$$
\limsup_{t\nearrow +\infty}\dfrac{\Phi(\eta t)}{\Phi(t)}<1\quad \hbox{for any $\eta >1$}
$$
or the more sharp {\em borderline case} given by
$$
\limsup_{t\rightarrow\infty} \frac{\Phi(\eta_0 t)}{\Phi (t)}=1\quad \mbox{for some } \eta_0>1
$$
in the proofs of Propositions \ref{prop:maximalboundaryprofile}  and \ref{prop:minimalboundaryprofile} as they  were used in \cite{ADR}.$\fin$
\label{rem:nopotencias} 
\end{rem}
\begin{rem}\rm Obviously, any eventual nonnegative large solution of
$$
-\trace \cA (\cdot)\rD^{2}v+\beta(v)\ge f\quad \hbox{in $\Omega$}
$$
admits a minimal behaviour given by \eqref{eq:minimalboundaryprofileg}, provided \eqref{eq:ellipticity}, \eqref{eq:ellipticityneraboundaty},  \eqref{eq:crecimientog} and $f\ge 0$, without the assumption \eqref{eq:criticalgrowthu}.$\fin$
\label{rem:supersolutionfpositive}
\end{rem}
\begin{rem}\rm As it was pointed out, the assumption  \eqref{eq:crecimientog}  is fulfilled for the choices of Remark \ref{rem:exemplesbeta}. Moreover, for the power like case $\beta _{m}(t)=t^{m}$ and $m>1$ the assumptions \eqref{eq:criticalgrowthd} and \eqref{eq:criticalgrowthu}	correspond with 
$$
\begin{array}{l}
\hspace*{-1cm}\disp \left(\dfrac{2(m+1)}{(m-1)^{2}}\right)  ^{\frac{m}{m-1}}\left (\dfrac{\lambda}{1-\ell_{d}}\right )^{\frac{m}{m-1}}\ell_{d}  \le \liminf_{\rd(x)\searrow 0} f(x)\rd(x)^{\frac{2m}{m-1}}\\ [.2cm]
\hspace*{2.5cm}
\disp \le  \limsup_{\rd(x)\searrow 0} f(x)\rd(x)^{\frac{2m}{m-1}}\le 
\left(\dfrac{2(m+1)}{(m-1)^{2}}\right)  ^{\frac{m}{m-1}}\left (\dfrac{\Lambda}{1-\ell_{u}}\right )^{\frac{m}{m-1}}\ell_{u} .
\end{array}
$$
For $\beta (t)=e^{t}$ the assumptions \eqref{eq:criticalgrowthd} and \eqref{eq:criticalgrowthu}  correspond with 
$$
\hspace*{-2cm}\disp 2\ell_{d} \dfrac{\lambda}{1-\ell_{d}}\le \liminf_{\rd(x)\searrow 0} f(x)\rd(x)^{2}\le \limsup_{\rd(x)\searrow 0} f(x)\rd(x)^{2}\le 
2\ell_{u} \dfrac{\Lambda}{1-\ell_{u}}.
$$
\fineqnum
\label{rem:exemplesg1}
\end{rem}
\begin{rem}\rm We note that under the assumptions of Propositions \ref{prop:maximalboundaryprofile}  and \ref{prop:minimalboundaryprofile} we have
$$
\limsup_{\rd(x)\searrow 0}\dfrac{u^{*}(x)}{\phi \left (\sqrt{\dfrac{1-\ell_{u}}{\Lambda}}\rd(x)\right )}	\le 1\le \liminf_{\rd(x)\searrow 0}\dfrac{v_{*}(x)}{\phi \left (\sqrt{\dfrac{1-\ell_{d}}{\lambda}}\rd(x)\right )}
$$
for any solution of
$$
-\trace \cA(\cdot)\rD^{2}u+\beta (u)\le f\quad \hbox{in $\Omega$},
$$
and any nonnegative large solution of
$$
-\trace \cA(\cdot)\rD^{2}v+\beta (v)\ge f\quad \hbox{in $\Omega$}.
$$
\fineqnum
\end{rem}
From the above results we deduce
\begin{theo}\rm Under the assumptions of Propositions \ref{prop:maximalboundaryprofile}  and \ref{prop:minimalboundaryprofile}
for any solution of
$$
-\trace \cA(\cdot)\rD^{2}u+\beta (u)\le f\quad \hbox{in $\Omega$}
$$
and any non negative large solution of
$$
-\trace \cA(\cdot)\rD^{2}v+\beta (v)\ge f\quad \hbox{in $\Omega$}
$$
one has
$$
\limsup_{\rd(x)\searrow 0}\dfrac{u^{*}(x)}{v_{*}(x)}\le \liminf_{\rd(x)\searrow 0}\dfrac{\phi  \left (\sqrt{\dfrac{1-\ell_{u}}{\Lambda}}\rd(x)\right )}{\phi  \left (\sqrt{\dfrac{1-\ell_{d}}{\lambda}}\rd(x)\right )}.
$$
Then \eqref{eq:UVontheboundary} is fulfilled whenever
\begin{equation}
\dfrac{1-\ell_{u}}{1-\ell _{d}}\ge \dfrac{\Lambda}{\lambda}
\label{eq:uniquenessboundary}
\end{equation}
holds. In particular if
\begin{equation}
\sqrt{\dfrac{1-\ell_{u}}{\Lambda}}=\rC=\sqrt{\dfrac{1-\ell_{d}}{\lambda}}
\quad \Leftrightarrow \quad \dfrac{1-\ell_{u}}{1-\ell _{d}}=\dfrac{\Lambda}{\lambda},
\label{eq:uniquenessbehaviourboundary}
\end{equation}
 the precise boundary behaviour of any nonnegative large solution $u$ of 
$$
-\trace \cA(\cdot)\rD^{2}u+\beta (u)=f\quad \hbox{in $\Omega$},
$$
is given by
$$
\limsup_{\rd(x)\searrow 0}\dfrac{u^{*}(x)}{\phi \big (\rC\rd(x)\big )}	\le 1\le \liminf_{\rd(x)\searrow 0}\dfrac{u_{*}(x)}{\phi \big (\rC\rd(x)\big)},
$$
thus
$$
\lim_{\rd(x)\searrow 0}\dfrac{u(x)}{\phi \big (\rC\rd(x)\big )}=1.
$$
\fineq
\label{theo:uniquenessboundary}
\end{theo}

\par
In the next result we complete the  obtainement of the inequality \eqref{eq:uniquenessboundary} by means of a sharp iterative device by M. Safonov (see \cite{DKS} and \cite{GM}).
\begin{theo}
Asume $\partial \Omega \in\cC^{2}$. Let us suppose \eqref{eq:ellipticity}, \eqref{eq:ellipticityneraboundaty}, \eqref{eq:crecimientog}, \eqref{eq:cA} and 
\begin{equation}
	\dfrac{1-\ell_{u}}{1-\ell_{d}}<\dfrac{\Lambda}{\lambda}.
	\label{eq:laotra}
\end{equation}
Then any subsolution $u$ and any non negative large solution $v$ 
of \eqref{eq:fullyequation}
satisfy the property \eqref{eq:UVontheboundary} provided \eqref{eq:criticalgrowthd} and \eqref{eq:criticalgrowthu} with $0\le \ell_{u},\ell_{d} <1$.
\label{theo:absoluteboundarycomparison}
\end{theo}
\proof  We claim that in this case the inequality
$$
\limsup_{\rd(x)\rightarrow 0}\dfrac{u^{*}(x)}{v_{*}(x)}\le 1
$$
also holds. Indeed,  \eqref{eq:laotra} implies that for $\varepsilon >0$ the inequality
$$
\varepsilon +\nu \phi  \left (\sqrt{\dfrac{1-\ell_{d}}{\lambda}}\rd(x)\big )\right )\ge \phi  \left (\sqrt{\dfrac{1-\ell_{u}}{\Lambda}}\rd(x)\big )\right )\big )\quad \hbox{for $\rd(x) $ small}
$$
(see \eqref{eq:casinopotencias}) holds, where
$$
\nu=\left (\dfrac{(1-\ell_{d})\Lambda}{(1-\ell_{u})\lambda}\right )^{\frac{1}{\alpha-1}}>1.
$$
It leads to
$$
0<\limsup_{\rd(x)\searrow 0}\dfrac{u^{*}(x)}{v_{*}(x)}=\rL\le\left (\dfrac{(1-\ell_{d})\Lambda}{(1-\ell_{u})\lambda}\right )^{\frac{1}{\alpha-1}}. 
$$
Assume $\rL>1$, otherwise the proof ends.  Then for every constant $k$ near $\rL$ such that $1<k<\rL$ there exists $x_{0}\in\Omega_{\delta_{\Omega}}$ satisfying 
$$
u^{*}(x_{0})>kv_{*}(x_{0})\quad \hbox{for $\rd(x_{0})$ small.} 
$$
We will use the notation  $[u^{*}>kv_{*}]\doteq\big \{x\in\Omega_{\rd(x_{0})}:u^{*}(x)>kv_{*}(x)\big \}$. It is clear, that in the open set $[u^{*}>kv_{*}]$ the assumption  \eqref{eq:crecimientog} implies
$$
\beta (u^{*})> \beta (kv_{*})>k^{\alpha}\beta(v_{*}).
$$
Let us consider $r=\dfrac{\rd(x_{0})}{2}$, for which $\rB_{r}(x_{0})\subset\Omega_{\delta_{\Omega}}$, and define the smooth function
$$
w(x)=\rC\big (r^{2}-|x-x_{0}|^{2}\big ),\quad x\in \rB_{r}(x_{0})
$$
that satisfies
$$
\trace \cA(x)\rD^{2}w(x)\ge -2\rC \Lambda,\quad x\in\rB_{r}(x_{0}).
$$
Then 
$$
\left \{
\begin{array}{l}
-\trace \cA(\cdot)\rD^{2}u^{*}\le f-k^{\alpha}\beta(v_{*}),\\ [.15cm]
-\trace \cA(\cdot)\rD^{2}\big (kv_{*}-w\big )\ge k\big (f-\beta(v_{*})\big )
-2\rC \Lambda
\end{array}
\right .
\quad \hbox{in  $[u^{*}>kv_{*}]\cap \rB_{r}(x_{0})$ }
$$
in viscosity sense. Since $r<\rd(x)<3r$, one proves
$$
\phi  \left (\sqrt{\dfrac{1-\ell_{d}}{\lambda}}\rd(x)\big )\right )\ge \phi  \left (\sqrt{\dfrac{1-\ell_{d}}{\lambda}}3r\big )\right ),
$$
hence Proposition \ref{prop:minimalboundaryprofile} implies
$$
k(k^{\alpha-1}-1)\beta\big (v_{*}(x)\big )\ge
k(k_{0}^{\alpha}-1)\beta \left ((1+\varepsilon)\phi  \left (\sqrt{\dfrac{1-\ell_{d}}{\lambda}}3r\big )\right )\right )
$$ 
for $\varepsilon>0$ and $\rd(x_{0})$ small, where $k\ge k_{0}>1$. So that,
$$
f-k^{\alpha}\beta(v_{*})-k\big (f-\beta(v_{*})\big )
+2\rC \Lambda\le 0\quad \hbox{in $[u^{*}>kv_{*}]\cap \rB_{r}(x_{0})$},
$$
provided
$$
\rC=\dfrac{k(k_{0}^{\alpha}-1)}{2\Lambda}\beta \left ((1+\varepsilon)\phi  \left (\sqrt{\dfrac{1-\ell_{d}}{\lambda}}3r\big )\right )\right ).
$$
By means of the Weak Maximum Principle without coercivity term (see Theorem \ref{theo:withoutcoercivityI}) one proves that there exists $x_{1}\in\partial \big ([u^{*}>kv_{*}]\cap \rB_{r}(x_{0})\big )$ such that
\begin{equation}
u^{*}(x)-k v_{*}(x)+w(x)\le u^{*}(x_{1})-kv_{*}(x_{1})+w(x_{1}),\quad x\in [u^{*}>kv_{*}]\cap \rB_{r}(x_{0}).
\label{eq:usingWMP}
\end{equation}
In fact, $x_{1}\in\partial \rB_{r}(x_{0})$ because   $x_{1}\in\rB_{r}(x_{0})$ implies $x_{1}\in\partial [u^{*}>kv]$, thus $u^{*}(x_{1})=kv(x_{1})$ whence the contradiction
$$
\rC r^{2}=w(x_{0})\le w(x_{1})<\rC r^{2}
$$
follows. So that, $x_{1}\in\partial \rB_{r}(x_{0})$ and $w(x_{1})=0$ as well as 
$$
u^{*}(x_{1})-kv_{*}(x_{1})=u^{*}(x_{1})-kv_{*}(x_{1})+w(x_{1})\ge u^{*}(x_{0})-kv_{*}(x_{0})+w(x_{0})\ge w(x_{0})=\rC r^{2}.
$$
Thus  
$$
u^{*}(x_{1})-kv_{*}(x_{1})\ge k\dfrac{\big (k_{0}^{\alpha-1}-1\big )r^{2}}{2\Lambda}\beta \left ((1+\varepsilon)\phi\left (\sqrt{\dfrac{1-\ell_{d}}{\lambda }}3r\right )\right ).
$$
Now the property  \eqref{eq:fromKellerOsserman2} implies
$$
u^{*}(x_{1})-kv_{*}(x_{1})
\ge k \dfrac{(1+\varepsilon)
\big (k_{0}^{\alpha-1}-1\big )r^{2}}{2\Lambda \varepsilon_{0}}
\phi\left (\sqrt{\dfrac{1-\ell_{d}}{\lambda }}3r\right ),
$$
for $0<\varepsilon_{0}$ and $\rd(x_{0})$ small. On the other hand, from \eqref{eq:casinopotencias} we deduce
$$
\varepsilon +\nu\phi\left (\sqrt{\dfrac{1-\ell_{d}}{\lambda }}3r\right )\ge
\phi\left (\sqrt{\dfrac{1-\ell_{u}}{\Lambda }}r\right )
$$
where
$$
\nu\doteq \left (9\dfrac{(1-\ell_{d})\Lambda}{(1-\ell_{u})\lambda}\right )^{\frac{1}{\alpha-1}}>1
$$
(see  \eqref{eq:laotra}). Since $v$ is a solution of \eqref{eq:fullyequation} we lead to 
$$
\phi\left (\sqrt{\dfrac{1-\ell_{d}}{\lambda }}3r\right )\ge\dfrac{1}{\nu}\phi\left (\sqrt{\dfrac{1-\ell_{u}}{\Lambda }}r\right )-
\dfrac{\varepsilon}{\nu}\ge \dfrac{1}{\nu}v^{*}(x_{1})-
\dfrac{\varepsilon}{\nu}
$$
(see Proposition \ref{prop:maximalboundaryprofile}).   So that, 
$$
u^{*}(x_{1})\ge k\left (1+ \dfrac{(1+\varepsilon)
\big (k_{0}^{\alpha-1}-1\big )r^{2}}{2\Lambda\nu\varepsilon_{0}(1-\varepsilon)} \right )v_{*}(x_{1})-
\dfrac{(1+\varepsilon)\varepsilon\big (k_{0}^{\alpha-1}-1\big )r^{2}}{2\Lambda\nu\varepsilon_{0}}.
$$
Therefore, from
$$
0<\limsup_{\rd(x)\searrow 0}\dfrac{u^{*}(x)}{v_{*}(x)}=\rL <+\infty 
$$
we deduce
$$
\rL +\varepsilon \ge k\left (1+ \dfrac{(1+\varepsilon)
	\big (k_{0}^{\alpha-1}-1\big )r^{2}}{2\Lambda\nu\varepsilon_{0}} \right ).
$$
So, the choice $k=\rL -\varepsilon >k_{0}>1$, thus $\varepsilon <\rL -k_{0}$, implies
$$
\rL+\varepsilon \ge (\rL -\varepsilon)\left (1+ \dfrac{(1+\varepsilon)
	\big (k_{0}^{\alpha-1}-1\big )r^{2}}{2\Lambda\nu\varepsilon_{0}(1-\varepsilon)} \right ).
$$
Therefore one arrives to a contradiction by letting $\varepsilon \searrow 0.\fin$
\begin{rem}\rm We also may obtain \eqref{eq:usingWMP} by using  Theorem~\ref{theo:withoutcoercivityII}
provided \eqref{eq:lightellipticity} and 
\eqref{eq:structureannulusSMPsub}.$\fin$
\end{rem}
\begin{rem}\rm Property \eqref{eq:uniquenessbehaviourboundary} requires $\ell_{d}\le \ell_{u}$ and
$$
\ell_{d}\le \liminf_{\rd(x)\searrow 0}\dfrac{f(x)}{\beta\left  (\phi \big  (\rC\rd(x)\big )\right )}\le \limsup_{\rd(x)\searrow 0}\dfrac{f(x)}{\beta\left  (\phi \big  (\rC\rd(x)\big )\right )}\le \ell_{u} < 1.
$$
When $\Lambda =\lambda$ equality \eqref{eq:uniquenessbehaviourboundary} follows if $0\le \ell_{u}=\ell_{d}<1$.  It extends the result obtained in \cite{DL} for the unit operator $\cA(\cdot)\equiv \rI.\fin$
\label{rem:uniquenessboundary}
\end{rem}
\setcounter{equation}{0}
\section{Extra Keller-Osserman type growth of the sources on the boundary}
\label{sec:KOgrowthextra}

As it was pointed out, first we illustrate the power like case corresponding to
$$
-\trace \cA(\cdot)\rD^{2}u+u^{m}=d^{-q}\quad \hbox{in $\Omega$},
$$
where we require 
$$
q>\dfrac{2m}{m-1},~m>1,
$$
in order to consider an extra Keller-Osserman type growth at the boundary of solutions for which we want to obtain the property \eqref{eq:UVontheboundary}
$$
\limsup_{\rd(x)\rightarrow 0}\dfrac{u^{*}(x)}{v_{*}(x)}\le 1.
$$
\begin{theo}[Comparison on the boundary (power-like case)] Asume $\partial \Omega \in \cC^{2}$. Suppose \eqref{eq:ellipticity}, \eqref{eq:ellipticityneraboundaty} and
\begin{equation}
\limsup_{\rd(x)\searrow 0} f(x)\rd(x)^{q}\le\ell_{u} \quad (\ell_{u}>0),
\label{eq:criticalgrowthextrapoweru}
\end{equation}
for $q>\dfrac{2m}{m-1},~m>1$.  Then any solution $u$ 
of 
$$
-\trace \cA(\cdot)\rD^{2}u+u^{m}\le f\quad \hbox{in $\Omega$}
$$
admits a maximal behaviour given by
\begin{equation}
\limsup_{\rd(x)\searrow 0}u^{*}(x)\rd(x)^{\frac{q}{m}}\le \ell_{u}^{\frac{1}{m}}.
\label{eq:maximalboundaryprofileextrapower}
\end{equation}
On the other hand, under 
\begin{equation}
\liminf_{\rd(x)\searrow 0} f(x)\rd(x)^{q}\ge \ell_{d} >0,
\label{eq:criticalgrowthextrapowerd}
\end{equation}
any large non negative supersolution $v$ of 
$$
-\trace \cA(\cdot)\rD^{2}v+v^{m}\ge f\quad \hbox{in $\Omega$}
$$
admits the  behaviour given by
\begin{equation}
\liminf_{\rd(x)\searrow 0}v_{*}(x)\rd(x)^{\frac{q}{m}}\ge \ell_{d}^{\frac{1}{m}}.
\label{eq:minimalboundaryprofileextrapower}
\end{equation}
Therefore, one satisfies \eqref{eq:UVontheboundary}.
\label{theo:comparisonboundaryextrapower}
\end{theo}
\proof  Once more the idea is to construct suitable classical large supersolutions. So,  for $0<\delta <\delta_{\Omega}$, we consider
$$
\rU_{0}(x)=\rC_{0} \sigma^{-\frac{q}{m}}, \quad x\in \Omega^{\delta}_{\delta_{\Omega}}\doteq \{x\in\Omega:~\delta <\rd(x)<\delta_{\Omega}\},
$$ 
where $\rC_{0}>0$ is independent on $\delta$. Here $\sigma \doteq \rd(x)-\delta$. Straightforward computations show
$$
-\trace \cA(x)\rD^{2}\rU_{0}(x)=\rC_{0} \dfrac{q}{m}\sigma ^{-\frac{q+2m}{m}}\bigg (-\frac{q+m}{m}\rE(x)+\sigma\trace \cA(x)\rD^{2}\rd(x)\bigg ),\quad x\in \Omega^{\delta}_{\delta_{1}},
$$
where $0<\delta_{1}<\delta_{\Omega}$. On the other hand, let $0<\varepsilon<
\rC_{0}^{m}-\ell_{u} $, assumed $\rC_{0}^{m}>\ell_{u}$. Then, 
the condition \eqref{eq:criticalgrowthextrapoweru} implies
$$
-f(x)\ge -(\ell_{u} +\varepsilon)\rd(x)^{-q}\ge  -(\ell_{u}+\varepsilon)\sigma^{-q},\quad x\in \Omega^{\delta}_{\delta_{1}}.
$$
So that
$$
\begin{array}{ll}
-\trace \cA(\cdot)\rD^{2}\rU_{0}+\rU_{0}^{m}-f&\hspace*{-.2cm}\ge \sigma^{-q}
\bigg  [-\rC_{0}\dfrac{q}{m}\sigma ^{\frac{q(m-1)-2m}{m}}\bigg (\dfrac{q	+m}{m}\rE(x)+\sigma\trace \cA(x)\rD^{2}\rd(x)\bigg )\\
&+\rC_{0}^{m}-(\ell_{u}+\varepsilon)\bigg  ],
\end{array}
$$
for $ x\in\Omega^{\delta}_{\delta_{1}}$. We note the inequality
$$
\dfrac{q(m-1)-m}{m}>\dfrac{q(m-1)-2m}{m}>0,
$$
whence
$$
-\trace \cA\rD^{2}\rU_{0}+\rU_{0}^{m}-f\ge \sigma^{-q}
\left [-\rC_{0}\dfrac{q}{m}\dfrac{q+m}{m}\Lambda\delta_{1} ^{\frac{q(m-1)-2m}{m}}-\rC_{0} \dfrac{q}{m}\rC_{\partial \Omega, \cA}\delta_{1} ^{\frac{q(m-1)-m}{m}}+\rC_{0}^{m}-(\ell_{u}+\varepsilon)\right ],
$$
holds in $ \Omega^{\delta}_{\delta_{1}}$. Since
$$
\rC_{0}\left [\dfrac{q}{m}\delta_{1}^{\frac{q(m-1)-2m}{m}}\left ( \frac{q+m}{m}\Lambda -\rC_{\partial \Omega,\cA}\delta_{1}\right )\right ]\le \rC_{0}^{m}-(\ell_{u}+\varepsilon)
$$
is fulfilled  for $\delta_{1}$ small enough uniformly on $\rC_{0}$, we obtain
$$
-\trace \cA(x)\rD^{2}\rU_{0}(x)+(\rU_{0}(x))^{m}\ge f(x),\quad x\in  \Omega^{\delta}_{\delta_{1}}.
$$
By means of arguments used in Proposition \ref{prop:maximalboundaryprofile} one deduces 
$$
\limsup_{\rd(x)\nearrow 0}u^{*}(x)\rd(x)^{\frac{q}{m}}\le \rC_{0}
$$
and then \eqref{eq:maximalboundaryprofileextrapower} holds by letting $\rC_{0} \searrow \ell_{u}^{\frac{1}{m}}$. 
Next for $0<\delta$ and $0<\delta_{1}\le \delta_{1}$, we construct
$$
\rV_{0}(x)=\rC_{0} \sigma^{-\frac{q}{m}}, \quad x\in \Omega _{\delta_{1}}.
$$ 
for $\sigma=\rd(x)+\delta$.  As above, straightforward computations show
$$
-\trace \cA(x)\rD^{2}\rV_{0}(x)=\rC_{0} \dfrac{q}{m}\sigma ^{-\frac{q+2m}{m}}\bigg (-\frac{q+m}{m}\rE(x)+\sigma\trace \cA(x)\rD^{2}\rd(x)\bigg ),\quad x\in \Omega_{\delta_{1}}.
$$
On the other hand, let  $0<\varepsilon<
\ell_{d} -\rC_{0}^{m}$, assumed $\ell_{d}>\rC_{0}^{m}$. Then, 
the condition \eqref{eq:criticalgrowthextrapowerd} implies
$$
-f(x)\le -(\ell_{d} -\varepsilon)\rd(x)^{-q}\ge  (\ell_{d} -\varepsilon)\big (\rd(x)+\delta \big )^{-q},\quad x\in \Omega _{\delta_{1}}.
$$
So that
$$
\begin{array}{ll}
-\trace \cA\rD^{2}\rV_{0}(x)+\big (\rV_{0}(x)\big )^{m}-f(x)&\hspace*{-.2cm}\le \sigma^{-q}\bigg  [-\rC_{0} \dfrac{q}{m}\dfrac{q+m}{m}\rE(x)\sigma^{\frac{q(m-1)-2m}{m}} \\ [.2cm]
&+\rC_{0} \dfrac{q}{m}\sigma^{\frac{q(m-1)-m}{m}}\trace\cA(x)\rD^{2}\rd(x)+\rC_{0}^{m} -(\ell_{d} -\varepsilon) \bigg ]\\ [.2cm]
&\hspace*{-.2cm}\le \sigma^{-q}\bigg  [\rC_{0} \dfrac{q}{m}\rC_{\partial\Omega,\cA}\sigma^{\frac{q(m-1)-m}{m}}+\rC_{0}^{m} -(\ell_{d} -\varepsilon) \bigg ],
\end{array}
$$
for $x\in \Omega_{\delta_{1}}$.
We note $\ell_{d} -\rC_{0}^{m} -\varepsilon >0$. 
Again we use the inequality
$$
\dfrac{q(m-1)-m}{m}>\dfrac{q(m-1)-2m}{m}\ge 0.
$$
Therefore, taking $\delta =\delta_{1}<\delta _{\Omega}$
$$
\rC_{0}\dfrac{q}{m}(2\delta_{1})^{\frac{q(m-1)-m}{m}}\rC_{\partial\Omega,\cA}\le (\ell_{d} -\varepsilon) -\rC_{0}^{m} 
$$
holds for some $\delta_{1}$ small enough uniformly on $\rC_{0}$. Thus
$$
-\trace \cA(x)\rD^{2}\rV_{0}(x)+\big (\rV_{0}(x)\big )^{m}\le f(x),\quad x\in \Omega_{\delta_{1}}.  
$$
So, by means of arguments used in Proposition \ref{prop:minimalboundaryprofile} we deduce
$$
\liminf_{\rd(x)\nearrow 0}v_{*}(x)\rd(x)^{\frac{q}{m}}\ge \rC_{0}
$$
and then \eqref{eq:minimalboundaryprofileextrapower} holds by letting $\rC_{0} \nearrow \ell_{d}^{\frac{1}{m}}$. Since $\ell_{d}\le \ell_{u}$ the property \eqref{eq:UVontheboundary} hods.$\fin$ 
\begin{rem}\rm Since the diffusive term of the PDE is negligible with respect to the other one  the ellipticity constans $\lambda$ and $\Lambda$ do not appear in the estimations \eqref{eq:maximalboundaryprofileextrapower} and \eqref{eq:minimalboundaryprofileextrapower}. In Theorem \ref{theo:secondarytermextraKOgrowth} we obtain a more general result.$\fin$
\end{rem}
\par
\medskip
With the above simple ideas we get the general result 
\begin{theo}[Comparison on the boundary (general case)] Asume $\partial \Omega \in \cC^{2}$ and $\beta\in\cC^{2}$. Suppose \eqref{eq:ellipticity}, \eqref{eq:ellipticityneraboundaty} and
\begin{equation}
\liminf_{\rd(x)\searrow 0} \dfrac{f(x)}{h\big (\rd(x)\big )}\le 1
\label{eq:criticalgrowthextrau}
\end{equation}
for some $\cC^{2}$ decreasing function $h:~\R_{+}\rightarrow \R_{+} $ such that $h(0)=\infty$ and $h'(0)=0$ satisfying the compatiblity condition \eqref{eq:compatibility}.
Futhermore, we assume that 
\begin{equation}
\left \{
\begin{array}{l}
\rH_{1}(\sigma)\doteq \dfrac{h(\sigma)}{\beta'\big ( \beta^{-1}\big (h(\sigma)\big )}\left (h''(\sigma)-\dfrac{\big (h'(\sigma )\big )^{2}\beta ''\big ( \beta ^{-1}\big (h(\sigma)\big )\big )}{\big (\beta '\big (\beta ^{-1}\big (h(\sigma)\big )\big )^{2}}\right ),\\ [.2cm]
\rH_{2}(\delta)\doteq -\dfrac{h'(\sigma)}{h(\sigma)\beta '\big (\beta^{-1}\big (h(\sigma)\big )\big )}
\end{array}
\right .
\label{eq:criticalgrowthextraexpoademas}
\end{equation}
are increasing functions  vanishing  at the origin. Then any subsolution $u$ 
of 
$$
-\trace \cA(\cdot)\rD^{2}u+\beta (u)=f\quad \hbox{in $\Omega$}
$$
admits a maximal behaviour given by
\begin{equation}
\limsup_{\rd(x)\searrow 0}\dfrac{u^{*}(x)}{\beta^{-1}\big (h(\rd(x))\big )}\le 1.
\label{eq:maximalboundaryprofileextrau}
\end{equation}
On the other hand, under
\begin{equation}
\liminf_{\rd(x)\searrow 0} \dfrac{f(x)}{h\big (\rd(x)\big )}\ge 1 
\label{eq:criticalgrowthextrad}
\end{equation}
any large nononegative supersolution $v$  admits the minimal behaviour given by
\begin{equation}
\liminf_{\rd(x)\searrow 0}\dfrac{v_{*}(x)}{\beta ^{-1}\big (h(\rd(x))\big )}\ge 1. 
\label{eq:minimalboundaryprofileextrad}
\end{equation}  
Hence assuming 
\begin{equation}
\lim_{\rd(x)\searrow 0} \dfrac{f(x)}{h\big (\rd(x)\big )}=1,
\label{eq:criticalgrowthextra}
\end{equation}
one has \eqref{eq:UVontheboundary}.
\label{theo:comparisonboundaryextrageneral}
\end{theo}
\proof  The idea is to construct a classical large subsolution. Let
$$
\rU_{0}(x)=\varrho \beta ^{-1}\big (h(\rd(x)-\delta)\big ),\quad x\in\Omega^{\delta}_{\delta_{1}},
$$
where $0<\delta<\delta_{1}<\delta_{\Omega}$ and $\varrho >1$. Denoting $\sigma =\rd(x)-\delta$, straightforward computations show
$$
\begin{array}{ll}
-\trace \cA(x)\rD^{2}\rU_{0}(x) &\hspace*{-.2cm}\ge  \Lambda \left (-\dfrac{h''(\sigma)}{\beta '\big (\beta ^{-1}\big (h(\sigma)\big )\big )} +\dfrac{\big (h'(\sigma )\big )^{2}g\beta ''\big (\beta ^{-1}\big (h(\sigma)\big )\big )}{\big (\beta '\big (\beta ^{-1}\big (h(\sigma)\big )\big )^{3}}\right )\\
& +\dfrac{h'(\sigma)}{\beta '\big (\beta ^{-1}\big (h(\sigma)\big )\big )}\cC_{\partial \Omega,\cA},\quad x\in\Omega^{\delta}_{\delta_{1}},\quad x\in\Omega^{\delta}_{\delta_{1}}
\end{array}
$$
(see \eqref{eq:criticalgrowthextraexpoademas}). On the other hand, let $0<\varepsilon<\varrho-1$. Then, 
the condition \eqref{eq:criticalgrowthextrau} implies
$$
-f(x)\ge -(1 +\varepsilon)h(\rd(x))\ge -(1+\varepsilon)h(\rd(x)-\delta),\quad x\in \Omega^{\delta}_{\delta_{1}}.
$$
Consequently
$$
-\trace \cA(x)\rD^{2}\rU_{0}(x)+\beta \big (\rU_{0}(x)\big )-f(x)\ge h(\delta_{1})\left [\varrho \bigg 
(-\Lambda \rH_{1}(\delta_{1})-\rH_{2}(\delta_{1})\cC_{\partial \Omega,\cA}\bigg )+\varrho-(1+\varepsilon)\right ],
$$
for $x\in \Omega^{\delta}_{\delta_{1}} $. Since
$$
\Lambda \rH_{1}(\delta_{1}) +\rH_{2}(\delta_{1})\le \dfrac{\varrho-(1+\varepsilon)}{\varrho}
$$
is fulfilled  for $\delta_{1}$ small enough (see \eqref{eq:criticalgrowthextraexpoademas}), we obtain
$$
-\trace \cA(x)\rD^{2}\rU_{0}(x)+\beta\big (\rU_{0}(x)\big )\ge f(x)\quad \hbox{for $x\in  \Omega^{\delta}_{\delta_{1}}$}.
$$
Finally one deduces 
$$
\limsup_{\rd(x)\nearrow 0}\dfrac{u^{*}(x)}{\beta^{-1}\big (h(\rd(x))\big )}\le \varrho 
$$
and then \eqref{eq:maximalboundaryprofileextrau} holds by letting $\varrho \searrow 1$. 
\par
\medskip
Next for $0<\delta$ and $0<\delta_{1}\le \delta_{\Omega}$, we construct
$$
\rV_{0}(x)=\varrho \beta ^{-1}\big (h(\rd(x)+\delta)\big ),\quad x\in\Omega_{\delta_{1}},
$$
where $0<\varrho<1 $. Let $0<\varepsilon<1-\varrho $. Then, 
the condition  \eqref{eq:criticalgrowthextrad} implies
$$
-f(x)\le -(1 -\varepsilon) h(\rd(x))\le -(1 -\varepsilon) h(\rd(x)+\delta),\quad x\in \Omega_{\delta_{1}}.
$$
As above,  by straightforward computations we have
$$
-\trace \cA(x)\rD^{2}\rV_{0}(x)+\beta\big (\rV_{0}(x)\big )-f(x)\le h(\delta_{1})\left [\varrho \bigg 
(-\lambda \rH_{1}(2\delta_{\Omega})+\rH_{2}(2\delta_{1})\cC_{\partial \Omega,\cA}\bigg )+\varrho-(1-\varepsilon)\right ],
$$
for $x\in \Omega _{\delta_{1}}$ (see  \eqref{eq:criticalgrowthextraexpoademas}).  Here by simplicity we choose $\delta=\delta_{1}<\delta_{\Omega}$. Since
$$
\rH_{2}(2\delta_{1})\le \dfrac{(1-\varepsilon)-\varrho}{\varrho}
$$
is fulfilled  for $\delta_{1}$ small enough, and
$$
-\trace \cA(x)\rD^{2}\rV_{0}(x)+\beta \big (\rV_{0}(x)\big )\le f(x)\quad \hbox{for $x\in  \Omega^{\delta}_{\delta_{1}}$}
$$
holds. Finally one deduces 
$$
\liminf_{\rd(x)\nearrow 0}\dfrac{v_{*}(x)}{\beta^{-1}\big (h(\rd(x))\big )}\ge \varrho.
$$
Now \eqref{eq:minimalboundaryprofileextrad} follows by letting $\varrho \nearrow 1 .\fin$ 

\begin{rem}\rm For the power-like case $\beta_{m}(t)= t^{m},~m>1$, and $h_{q}(\delta)=\delta^{-q}$ one has
$$
\left \{
\begin{array}{l}
	\rH_{1}(\sigma)=\dfrac{q}{m}\dfrac{q+m}{m}\sigma^{-\frac{q(m-1)-2m}{m}},\\ [.275cm]
	\rH_{2}(\sigma)=\dfrac{q}{m}\sigma^{\frac{q(m-1)-m}{m}}.
\end{array}
\right .
$$
Certainly, both functions are increasing functions vanishing at the origin provided $q>\dfrac{2m}{m-1}.\fin$
\label{rem:examplepower}
\end{rem}
\begin{rem}[Bieberbach choice]\rm For the exponential case $\beta (t)=e^{t}$ the Keller-Osserman boundary profile is governed by 
$$
\phi(\delta )=\log \left (\dfrac{2}{\delta^{2}}\right )
$$
(see Remark \ref{rem:exemplesbeta}). Here, available extra Keller-Osserman boundary profiles are governed by functions $h(\delta)$ satisfying
$$
h(\delta)>\dfrac{2}{\delta^{2}}\quad \hbox{for $\delta$ small}
$$
(see \eqref{eq:compatibility}) such that  
$$
\left \{
\begin{array}{l}
	\rH_{1}(\sigma)=h''(\sigma)-\dfrac{\big (h'(\sigma) \big)^{2}}{h(\sigma)},\\ [.275cm]
	\rH_{2}(\sigma)=-\dfrac{h'(\sigma)}{\big (h(\sigma)\big )^{2}}
\end{array}
\right .
$$
are increasing functions vanishing at the origin. In order to avoid tedious computations we only give the simple example $h(\delta)=\dfrac{\mu}{\delta^{2}}$ for $\mu >2.\fin$
\label{rem:exampleexpo}
\end{rem}
\par
\noindent 
\par
\medskip
In this case of extra Keller-Osserman type growth by means of an argument of rescaling the  explositivity of the sources is transferred  to the supersolutions. 
\begin{theo}[Transferring the explosivity] Suppose $q\ge \dfrac{2m}{m-1}$. Let $f$ be a source for which
\begin{equation}
\liminf_{\rd(x)\searrow 0}f(x)\rd(x)^{q}=\ell_{d}>0.
\label{eq:sourcestrasferringpower}
\end{equation}
The any nononegative supersolution $v$ of
$$
-\trace \cA(\cdot)\rD^{2}v+v^{m}= f\quad \hbox{in $\Omega$}
$$
is a large supersolution. 
\label{theo:explosivitytransferred}
\end{theo}
\proof Let $x_{0}\in\Omega_{\delta_{\Omega}}$ and $2r=\rd(x_{0})$ with $r<1$. Then the supersolution $v_{*}$ solves, in the viscosity sense,
$$
\left \{
\begin{array}{ll}
-\trace \cA(\cdot)\rD^{2}v_{*}+v_{*}^{m}\ge \widehat{\rC}r^{-q},&\quad x\in \bB_{r}(x_{0}),\\ [.15cm]
v_{*}(x)\ge 0 ,&\quad x\in \partial \bB_{r}(x_{0}),
\end{array}
\right .
$$
for $\widehat{\rC}=\ell_{d} 2^{-q}$. Since $\cA$ is uniform elliptic in $\Omega_{\delta_{\Omega}}$ there exists a classical solution $\rW$ of 
$$
\left \{
\begin{array}{ll}
-\trace \cA(x)\rD^{2}\rW(x)+\rW(x)^{m}=\widehat{\rC}r^{-q},&\quad x\in \bB_{r}(x_{0}),\\ [.15cm]
\rW(x)=0 ,&\quad x\in \partial \bB_{r}(x_{0}).
\end{array}
\right .
$$
Then $v_{*}(x)\ge \rW(x),~x\in \overline{\bB}_{r}(x_{0})$. Next, we construct $\widehat{\rW}(x)=\rW(x-x_{0})$ that solves
$$
\left \{
\begin{array}{ll}
-\trace \cA(x-x_{0})\rD^{2}\widehat{\rW}(x)+\widehat{\rW}(x)^{m}=\widehat{\rC}r^{-q},&\quad x\in \bB_{r}(0),\\ [.15cm]
\widehat{\rW}(x)=0 ,&\quad x\in \partial \bB_{r}(0),
\end{array}
\right .
$$
hence
$$
-\trace \cA(x-x_{0})\rD^{2}\widehat{\rW}(x)+\widehat{\rW}(x)^{m}\ge \widehat{\rC}r^{-\frac{2m}{m-1}},\quad x\in \bB_{r}(0)
$$
or
$$
r^{\frac{2m}{m-1}}\trace \cA(x-x_{0})\rD^{2}\widehat{\rW}(x)+r^{\frac{2m}{m-1}}\widehat{\rW}(x)^{m}\ge \widehat{\rC},\quad x\in \bB_{r}(0)
$$
for $r$ small. Since $\dfrac{2}{m-1}+2=\dfrac{2m}{m-1}$, the function
$\rW_{r}(x)=r^{\frac{2}{m-1}}\widehat{\rW}(rx)$ satisfies
$$
-\trace \cA(x-x_{0})\rD^{2}\rW_{r}(x)+\rW_{r}(x)^{m}\ge 
\widehat{\rC},\quad x\in\bB_{1}(0),
$$
and we get the problem
$$
\left \{
\begin{array}{ll}
-\trace \cA(x-x_{0})\rD^{2}\rV(x)+\rV(x)^{m}=\widehat{\rC},&\quad x\in \bB_{1}(0),\\ [.15cm]
\rV(x)=0 ,&\quad x\in \partial \bB_{1}(0).
\end{array}
\right .
$$
A simple adaptation of Strong Maximum Principle (see \cite{V1}) shows $\rV>0$ in $\bB_{1}(0)$ whence
$$
0<\rV(0)\le \left (\dfrac{\rd(x_{0})}{2}\right )^{\frac{2}{m-1}}v_{*}(x_{0})
$$
concludes the result. $\fin$

\setcounter{equation}{0}

\section{Comparison, uniqueness and existence of solutions}
\label{sec:existence}

As it was pointed out in Section \ref{sec:intro} the eventual degeneracy inside the domain of the leading part of the equation prevents, in general, to consider $\cC^{2}$ solutions. They are replaced by viscosity solutions that solve the PDE more intrinsically. We send to the mongraphy \cite{CIL} for a full introduction of the Viscosity Solution Theory. In short, we recall that for any function $u:\cO\rightarrow \R,~\cO\subseteq\R^{\rN}$, locally bounded from above, one defines the upper semi-continuous envelope
$$
u^{*}(x)=\limsup_{y\rightarrow x\atop y\in\cO}u(y),\quad x\in\cO
$$ 
and the superjets
$$
\cJ^{2,+}_{\cO}u^{*}(x)=\big \{(\bp,\rX)\in\R^{\rN}\times \cS^{\rN}:~u^{*}(y)\le u^{*}(x)+\pe{\bp}{y-x}+\dfrac{1}{2}\pe{\rX(y-x)}{y-x}+o(|y-x|^{2}) \big \}.
$$
In particular, a locally bounded from above function $u:\Omega\rightarrow \R$ is a (viscosity) subsolution of 
$$
-\trace \cA(\cdot)\rD^{2}u +\beta (u)=f\quad \hbox{in $\Omega$}
$$
if 
$$
-\trace \cA(x)\rX+\beta \big (u^{*}(x)\big )\le f(x)
$$
for all $(\bp,\rX)\in \cJ^{2,+}_{\Omega}u^{*}(x),~x\in\Omega$. By simplicity, sometimes we say that $u$ is a solution of
$$
-\trace \cA(\cdot)\rD^{2}u +\beta (u)\le f\quad \hbox{in $\Omega$}.
$$
\begin{rem}\rm From the above definition we deduce a property near the well known verification results of the Optimal Control Problems. More precisely any (viscosity) subsolution of 
$$	-\trace \cA(\cdot)\rD^{2}u +\beta (u)=f\quad \hbox{in $\Omega$}
$$
satisfies
\begin{equation}
u^{*}(x)\le \inf_{\rV\in\cV_{+}^{u}}\rV(x):,\quad x\in \overline{\Omega},
\label{eq:supverification}
\end{equation}
where
$$
\cV_{+}^{u}=\sup \big \{\rV\in\cC^{2} \hbox{ supersolution in $\Omega$ such that $u^{*}\le \rV$ on $\partial \Omega$} \big \}.
$$
Indeed, given a such function $\rV$, if there exists a positive interior maximum point of $u^{*}-\rV$ at some $x_{0}\in \Omega$ one has
$$
u^{*}(x)-u^{*}(x_{0})\le \rV(x)-\rV(x_{0})=\pe{\rD \rV(x_{0})}{x-x_{0}}+\dfrac{1}{2}\pe{\rD^{2}\rV(x_{0})}{x-x_{0}}+o(|x-x_{0}|^{2})^,
$$
thus $\big (\rD \rV(x_{0}),\rD ^{2}\rV(x_{0})\big )\in \cJ^{2,+}_{\Omega}u^{*}(x_{0})$ for which  we derive the contradiction
$$
-\trace \cA(x_{0})\rD^{2}\rV(x_{0})+\beta \big (u^{*}(x_{0})\big )\le 
-\trace \cA(x_{0})\rD^{2}\rV(x_{0})+\beta \big (\rV(x_{0})\big ).
$$
Therefore 
$$
u^{*}(x)-\rV(x)\le u^{*}(x_{0})-\rV(x_{0})\le 0,\quad x\in\overline{\Omega}
$$
and \eqref{eq:supverification} holds.$\fin$
\label{rem:supverification}
\end{rem}

Analogously, for $v:\cO\rightarrow \R,~\cO\subseteq\R^{\rN}$, locally bounded from below, one defines the lower semi-continuous envelope
$$
v^{*}(x)=\limsup_{y\rightarrow x\atop y\in\cO}v(y),\quad x\in\cO
$$ 
and the subjets
$$
\cJ^{2,-}_{\cO}u^{*}(x)=\big \{(\bp,\rX)\in\R^{\rN}\times \cS^{\rN}:~v_{*}(y)\ge v_{*}(x)+\pe{\bp}{y-x}+\dfrac{1}{2}\pe{\rX(y-x)}{y-x}+o(|y-x|^{2}) \big \}.
$$
In particular, a locally bounded from below function $v:\Omega\rightarrow \R$ is a (viscosity) supersolution of 
$$
-\trace \cA(\cdot)\rD^{2}v +\beta (v)=f\quad \hbox{in $\Omega$}
$$
if 
$$
-\trace \cA(x)\rX+\beta \big (v_{*}(x)\big )\ge f(x)
$$
for all $(\bp,\rX)\in \cJ^{2,-}_{\Omega}v_{*}(x),~x\in\Omega$. Sometimes we say that $v$ is a solution of
$$
-\trace \cA(\cdot)\rD^{2}v +\beta (u)\ge f\quad \hbox{in $\Omega$}.
$$
\begin{rem}\rm Arguing as in Remark	\ref{rem:subverification} one proves that any  (viscosity) supersolution of 
$$	
-\trace \cA(\cdot)\rD^{2}v +\beta (v)=f\quad \hbox{in $\Omega$}
$$
satisfies
\begin{equation}
\sup_{\rU\in\cV_{-}^{v}} \rU(x)\le v_{*}(x),\quad x\in \overline{\Omega},
\label{eq:subverification}
\end{equation}
where
$$
\cV_{-}^{v}=\big \{\rU\in \cC^{2} \hbox{ is a subsolution in $\Omega$ such that $\rU\le v_{*}$ on $\partial \Omega$} \big \}.
$$
\fineq
\label{rem:subverification}
\end{rem}

Clearly, a locally bounded function $u:~\Omega\rightarrow \R$ is a (viscosity) solution if its upper semi-continuous envelope $u^{*}$ and  its lower semi-continuous envelope $u_{*}$
are sub and supersolution, respectively. We send to the monography \cite{CIL} for the theory of (viscosity) solutions.
\par
Next, we study a kind of  superlinearity  property 
\begin{equation}
q(t+s)\ge q(s)+q(t),\quad r,s\ge 0,
\label{eq:superlineal}
\end{equation}
which  can be transferred under a suitable condition
\begin{lemma}[Superlinearity transferred] Let $q$ and $\beta $ two nonnegative functions on $\R_{+}$ such that
$$
\dfrac{\beta (t)}{q(t)}\quad \hbox{is nondecreasing}.
$$
Then if $q(t)$ verifies {\rm \eqref{eq:superlineal}} the same property holds for $\beta (t)$.
\label{lemma:trsaferingsuoerlinearity}
\end{lemma}
\proof It follows from
$$
\beta (t)+\beta (s)=q(t)\dfrac{\beta (t)}{q(t)}+q(s)\dfrac{\beta (s)}{q(s)}\le \big (q(t)+q(s)\big )\dfrac{\beta(t+s)}{q(t+s)}\le \beta (t+s).
$$
\fineq
\begin{rem}\rm Since the identity is a function satisfying \eqref{eq:superlineal} and $t^{\alpha-1}$ is a nondecreasing function if $\alpha\ge 1$ we deduce that $q(t)=t^{\alpha},~\alpha \ge 1$ satisfies  \eqref{eq:superlineal}.$\fin$
\end{rem}

We give some details on the result of uniqueness. First we obtain the uniqueness of eventual classical solutions of \eqref{eq:fullyequation} in a direct way when
$$
\limsup_{\rd(x)\rightarrow 0}\dfrac{u(x)}{v(x)}\le 1
$$
is assumed (see \eqref{eq:UVontheboundary}). Since we only intend to give an idea of the reasoning we omit additional asumptions on the ellipticity of the $\cA(x)$ for which \eqref{eq:fullyequation} admits classical solutions in this motivation (see \cite{Ev} or \cite{Sa}). 

\begin{prop}[Classical Maximum Principle] Suppose \eqref{eq:crecimientog1}.
Let $f$ and $g$ be two continuous functions. If $u,v\in \mathcal{C}^2(\Omega)$ are two nonnegative functions verifying
$$
\left \{
\begin{array}{l}
-\trace \cA(x) \rD^{2}u(x) + \beta \big (u(x)\big )\le  f(x)\\ [.2cm]
-\trace \cA(x)\rD^{2}v(x)+\beta\big (v(x)\big )\ge g(x),
\end{array}
\quad x\in \Omega,
\right .
$$
and
\begin{equation}
u(x)\le v(x),\quad 	x\in\partial \Omega
\label{eq:desionboundarybounded}
\end{equation}
then 
\begin{equation}
u(x)\le v(x)+\beta^{-1}\left (\n{(f-g)_{+}}_{\infty}\right ),\quad x\in\Omega.
\label{eq:comparison}
\end{equation}
When we replace \eqref{eq:desionboundarybounded} by
\begin{equation}
	\limsup_{x\rightarrow \partial\Omega}\frac{u(x)}{v(x)}\le 1\quad \hbox{and $v>0$ near $\partial \Omega$}
	\label{eq:desionboundary}
\end{equation}
 the inequality  \eqref{eq:comparison} also holds.
Here $r_{+}=\max \{r,0\}$.
\label{prop:classicalcomparison}
\end{prop}
\proof If the maximum of the continuous function $u-v$ is attained on $\partial \Omega$ or it is a nonegative value the inequality \eqref{eq:comparison} follows.
On the contrary, let us assume that $u-v$ admits some point $x_{0}\in \Omega$ such that
$$
(u-v)(x_{0})=\max_{\Omega}~(u-v)>0.
$$
Since
\begin{equation}
0\ge \trace \cA(x_{0})\rD^{2}(u-v)(x_{0})
\label{eq:classicalmaxomum}
\end{equation}
we deduce the inequality
\begin{equation}
(f-g)(x_{0})+\trace \cA(x_{0})\rD^{2}(u-v)(x_{0})\ge \beta \big (u(x_{0})\big )-\beta \big (v(x_{0})\big )>0.
\label{eq:contradiction}
\end{equation}
Then from Lemma \ref{lemma:trsaferingsuoerlinearity} we get
$$
\beta \big ((u-v)(x_{0})\big )\le (f-g)_{+}(x_{0})
$$
and \eqref{eq:comparison} follows. On the other hand, it is clear that 
\begin{equation}
\limsup_{x\rightarrow \partial\Omega}\frac{u(x)}{v(x)}< 1\quad \hbox{and $v>0$ near $\partial \Omega$}
\label{eq:unicityonboundarystrict}
\end{equation}
lead to $u\le v$ near $\partial \Omega$ and the above reasoning is true. The more general assumption \eqref{eq:desionboundary} implies
$$
\limsup_{x\rightarrow \partial\Omega}\frac{u(x)}{(1+\varepsilon)v(x)}<1,
$$
for all $\varepsilon >0$ and then $u\le (1+\varepsilon)v$ near $\partial \Omega$. Moreover, from \eqref{eq:crecimientog1} we have
$$
-(1+\varepsilon)\trace \cA(\cdot)\rD^{2}v+\beta \big ((1+\varepsilon)v\big )\ge (1+\varepsilon)\left (-\trace \cA(\cdot)\rD^{2}v+\beta (v)\right )\ge (1+\varepsilon)g \quad \hbox{in $\Omega$}.
$$
The above reasoning leads to
$$
\beta (u-(1+\varepsilon)v)(x_{0})\le (f-g)(x_{0})-\varepsilon g(x_{0})\le
\n{(f-g)_{+}}_{\infty}+\varepsilon \n{g}_{\infty}
$$
thus
$$
u(x)-(1+\varepsilon)v(x)\le \beta^{-1}\left (\n{(f-g)_{+}}_{\infty}+\varepsilon \n{g}_{\infty}\right ),\quad x\in\Omega,
$$
and the result follows by letting $\varepsilon \rightarrow  0.\fin$
\begin{rem}\rm For a direct Comparison Principle as 
$$
\hbox{$u\le v$ on $\partial \Omega$  and $f\le g$ in $\Omega\quad \Rightarrow \quad u\le v$ in $\Omega$}
$$ 
the assumption \eqref{eq:crecimientog1} is not required (see \eqref{eq:classicalmaxomum} and  \eqref{eq:contradiction}).$\fin$
\label{eq:directcomparison}
\end{rem}
\par
\medskip
A simple consequence proves that the condition 
\begin{equation}
	\displaystyle\limsup_{x\rightarrow \partial\Omega}\frac{u(x)}{v(x)}=1
\label{eq:unicityonboundary}
\end{equation}
implies uniqueness.
\begin{coro}
Suppose \eqref{eq:crecimientog1}. If $u,v\in \mathcal{C}^2(\Omega)$ are two non negative functions for which
$$
-\trace \cA(x) \rD^{2}u(x)+\beta \big (u(x)\big )=f(x)\quad\hbox{in $\Omega$}
$$
and \eqref{eq:unicityonboundary} hold, then $u=v$ on $\Omega$, provided that $f$ is a non negative continuous function.$\fin$
\label{coro:uniquenessmonotonia}
\end{coro}
\par
So that, a uniqueness condition on the boundary as 
$$
\lim_{x \rightarrow \partial \Omega}\dfrac{u(x)}{v(x)}= 1
$$
plays an important role. Next, we extend Proposition \ref{prop:classicalcomparison} in the Viscosity Solution Theory. 
Hark back to the above reasoning. Let us assume that the upper semi continuous function $u^{*}-v_{*}$ admits some interior point $x_{0}\in \Omega$ such that
$$
(u^{*}-v_{*})(x_{0})=\max_{\overline{\Omega}}~(u^{*}-v_{*})>0.
$$
As in the proof of Proposition \ref{prop:classicalcomparison}, the aim is to work by means of an inequality as \eqref{eq:classicalmaxomum}. Unfortunately $u^{*}$ and $v_{*}$ are not $\cC^{2}$ functions. If we use semijets $(p,\rX)\in \cJ^{2,+}u^{*}(x_{0})$ and
$(q,\rY)\in \cJ^{2,-}v_{*}(x_{0})$  the device goes if
$$
\trace \cA(x_{0})(\rX-\rY)\le 0
$$
holds for adequate $\rX,\rY\in \cS^{\rN}$ but strong dificulties appear due to the semijets $\cJ^{2,+}u^{*}(x_{0})$ and $\cJ^{2,-}v_{*}(x_{0})$  are set valued functions. The reasoning of the Jensen Maximum Principle  is very advantageous (see \cite{J}, \cite{Is} or \cite{CIL} for details). Essentially, it is about working close $x_{0}$ by means of a device that doubles  the number of variables penalizating this doubling (see the works by N.S.  Kruzkov, for instance \cite{Kr}). Indeed, for $\alpha >0$ the upper semi-continuous function
$$
\Psi(x,y)\doteq u^{*}(x)-v_{*}(y)-\dfrac{\alpha}{2}|x-y|^{2},\quad (x,y)\in \Omega\times \Omega
$$
has a finite supremum equals to $\rM_{\alpha}$ and 
admits $(x_{\alpha},y_{\alpha})$ such that
$$
\lim_{\alpha\nearrow \infty}\left (\rM_{\alpha}-\Psi(x_{\alpha},y_{\alpha})\right )=0
$$
and
$$
\left \{
\begin{array}{l}
\disp \lim_{\alpha \rightarrow \infty} \alpha |x_{\alpha}-y_{\alpha}|^{2}=0,\\ [.15cm]
\disp \lim_{\alpha \rightarrow \infty} \rM_{\alpha}=(u^{*}-v_{*})(x_{0})
\end{array}
\right .
$$
(see \cite[Lemma 3.1]{CIL}). Moreover, there exist $\rX,\rY\in\cS^{\rN}$ such that
$$
\big (\alpha(x_{\alpha}-y_{\alpha},\rX\big )\in \overline{\cJ}^{2,+}u^{*}(x_{\alpha})\quad\hbox{and}\quad \big (\alpha(y_{\alpha}-x_{\alpha},\rY\big )\in \overline{\cJ}^{2,-}v_{*}(y_{\alpha})
$$
and
\begin{equation}
-3\alpha \left (
\begin{array}{ll}
\mI & \bcero\\
\bcero & \mI
\end{array}
\right ) \le
\left (
\begin{array}{lr}
	\rX & \bcero\\
	\bcero & -\rY
\end{array}
\right )\le 3\alpha
\left (
\begin{array}{lr}
	\mI & -\mI\\
	\mI & \mI
\end{array}
\right ) 
\label{eq:XYcondition}
\end{equation}
In fact, \eqref{eq:XYcondition} implies $\rX\le \rY$ (see \cite[Theorem 3.2]{CIL}). So that, some suitable adptations near $x_{0}$ can be available in order that the aim can go. For large $\alpha$ one has
$$
(u^{*}-v_{*})(x_{0})\le \rM_{\alpha}\le u^{*}(x_{\alpha})-v_{*}(y_{\alpha})
$$ 
and
$$
\beta \big (u^{*}(x_{\alpha})\big )-\beta \big (v_{*}(y_{\alpha})\big )\le \trace \cA(x_{\alpha})\rX-\trace \cA(y_{\alpha})\rY\le \omega\big (\alpha|x_{\alpha}-y_{\alpha}|^{2}\big )
$$
assuming that there exists a kind of modulus of continuity $\omega : \R_{+}\rightarrow \R_{+}$, with $\omega(0^{+})=0$ such that 
\begin{equation}
\trace \cA(x)\rX-\trace \cA(y)\rY\le \omega\big (\alpha|x-y|^{2}\big )
\label{eq:cA}
\end{equation}
whenever $x,y\in\Omega,~\alpha>0,$ and $\rX,\rY\in\cS^{\rN}$ for which 
\eqref{eq:XYcondition} holds. Again, Lemma \ref{lemma:trsaferingsuoerlinearity} implies 
\begin{equation}
\begin{array}{ll}
\beta \big ((u^{*}-v_{*})(x_{0})\big )&\hspace*{-.2cm }\le \beta \big (u^{*}(x_{\alpha})-v_{*}(y_{\alpha})\big )\le
\beta \big (u^{*}(x_{\alpha})\big )-\beta \big (v_{*}(y_{\alpha})\big )\\ [.15cm]
&\le \omega\big (\alpha|x_{\alpha}-y_{\alpha}|^{2}\big )+f(x_{\alpha})-g(y_{\alpha}) \\ [.15cm]
&\le \omega\big (\alpha|x_{\alpha}-y_{\alpha}|^{2}\big )+\omega_{f}\big (|x_{\alpha}-y_{\alpha}|\big )+\n{(f-g)_{+}}_{\infty}
\end{array}
\label{eq:capitalinequality}
\end{equation}
where $\omega_{f}$ is the modulus of continuity of $f$. Finally, 
 the reasoning of Proposition \ref{prop:classicalcomparison} ends by letting $\alpha\nearrow \infty$. What we have tried is 
\begin{theo}[Weak Maximum  Principle] Assume \eqref{eq:crecimientog1} and \eqref{eq:cA}. Let $u$ be a discontinuous subsolution of 
$$
-\trace \cA(\cdot)\rD^{2}u+\beta (u)=f, \quad x\in \Omega
$$
where $f$ is a uniformly continuous function
and $v$ be a discontinuous supersolution of
$$
-\trace \cA(\cdot)\rD^{2}v+\beta (v)=g, \quad x\in \Omega 
$$
where $g$ is a continuous function.  Then
\begin{equation}
\limsup_{x\rightarrow \partial\Omega}\frac{u^{*}(x)}{v_{*}(x)}\le 1
\label{eq:boundarycomparisonagain}
\end{equation}
implies
$$
u^{*}(x)\le v_{*}(x)+\beta ^{-1}\left (\n{(f-g)_{+}}_{\infty}\right ),\quad x\in\Omega.
$$
Moreover, if  $u$ is a solution of
$$
-\trace \cA(\cdot)\rD^{2}u+\beta (u)=f \quad x\in \Omega 
$$
one has
$$
u^{*}(x)\le u_{*}(x),\quad x\in\Omega,
$$
whence the function $u$ is continuous in $\Omega.\fin$
\label{theo:comparisonprincipleviscosity}
\end{theo}
\begin{rem}\rm We send to Theorems \ref{theo:uniquenessboundary} and  \ref{theo:absoluteboundarycomparison} in order to obtain \eqref{eq:boundarycomparisonagain}.$\fin$
\end{rem}
\begin{rem}\rm In the proof of Theorem \ref{theo:comparisonprincipleviscosity}
the assumption of degenerate ellipticity is not assumed explicitly because \eqref{eq:cA} implies that the matrix function $\cA$ is degenerate elliptic. Indeed, given $\rY\in\cS_{+}^{\rN}$ Cauchy inequality implies
$$
-\pe{\rY\eta}{\eta}\le \varepsilon |\eta|^{2}+\dfrac{1}{\varepsilon}\n{\rY}^{2}\n{\eta}^{2},\quad \varepsilon>0,~\eta \in\R^{\rN},
$$
where $\n{\rY}=\sup_{\n{\eta}\le 1} |\pe{\rY\eta}{\eta}|$. Clearly, we also have the inequalty
$$
-\pe{\rY\eta}{\eta}\le \varepsilon |\eta|^{2}+\left (1+\dfrac{\n{\rY}}{\varepsilon}\right )\n{\rY}\n{\eta}^{2}
$$
that may be rewiten as 
$$
\left (
\begin{array}{lc}
	\bcero & \bcero\\
	\bcero & -(\rY+\varepsilon \rI)
\end{array}
\right )\le \left (1+\dfrac{\n{\rY}}{\varepsilon}\right )\n{\rY}
\left (
\begin{array}{lr}
	\rI & -\rI\\
	\rI & \rI
\end{array}
\right ).
$$	
By choosing $\alpha >\n{\rY}\max \left \{1,1+\dfrac{\n{\rY}}{\varepsilon}\right \}$ and $\varepsilon$ small one obtains \eqref{eq:XYcondition} for the matrix $\rX=\bcero$ and $\rY+\varepsilon\rI$, whence \eqref{eq:cA} implies
$$
-\cA\left (y-\frac{z}{\alpha}\right )(\rY+\varepsilon\rI)=\cA(y)\bcero -\cA\left (y-\frac{z}{\alpha}\right )(\rY+\varepsilon\rI)\le
\omega \left  (\dfrac{\n{z}^{2}}\alpha\right )
$$
for any $z\in\R^{\rN}\setminus\{0\}$. Then letting $\alpha\rightarrow \infty$ and $\varepsilon \searrow 0$ concludes
$$
0\le \cA(y)\xi\otimes\xi,\quad y\in\Omega,~\xi\in\R^{\rN},
$$
provided $\rY=\xi\otimes \xi\in\cS_{+}^{\rN},~\xi\in \R^{\rN}$. See \cite[Remark 3.4]{CIL} for a more general proof. When
$$
\cA(x)=\sigma(x)\sigma(x)^{\tt t}
$$
where $\sigma :~\Omega\rightarrow \R^{\rN}\times\R^{N}$is a Lipschitz continuous function with constant $\rL$ the assumption	\eqref{eq:cA} becomes
\begin{equation}
\trace \cA(x)\rX-\trace \cA(y)\rY\le 3\rL^{2}\alpha |x-y|^{2},
\label{eq:cAcontrol}
\end{equation}
thus $\omega(r)=3\rL^{2}r$ (see \cite[Example 3.6]{CIL}).$\fin$
\end{rem}
\begin{rem}\rm Remark \ref{eq:directcomparison} can be extended to Theorem \ref{theo:comparisonprincipleviscosity}. Indeed, a direct comparison is deduced from  \eqref{eq:capitalinequality} that becomes
$$
\beta \big (u^{*}(x_{\alpha})\big )-\beta \big (v_{*}(y_{\alpha})\big )\le \omega\big (\alpha|x_{\alpha}-y_{\alpha}|^{2}\big )+\omega_{f}\big (|x_{\alpha}-y_{\alpha}|\big )
$$
if we assume $f\le g$ in $\Omega$. Therefore an inequality as
$$
(u^{*}-v_{*})(x_{0})=\max_{\overline{\Omega}}~(u^{*}-v_{*})>0
$$
implies the contradiction 
$$
0<\beta \big (u^{*}(x_{0})\big )-\beta \big (v_{*}(x_{0})\big )\le 0.
$$
Consequently, we have proved a direct Comparison Principle as 
$$
\hbox{$u^{*}\le v_{*}$ on $\partial \Omega$  and $f\le g$ in $\Omega\quad \Rightarrow \quad u^{*}\le v_{*}$ in $\Omega$},
$$ 
under \eqref{eq:cA} and the uniformly continuous of the function $f$. We emphasize that the assumption~\eqref{eq:crecimientog1} is not required.$\fin$
\label{eq:directcomparisonviscosity}
\end{rem}
\par
\medskip
The existence of solution follows from the well known Perron method as, for instance \cite{IsPe} or~\cite{CIL}. For several reasons first we sketch an adaptation of the plan of~\cite{DL} on the {\em elliptic regularization} 
\begin{equation}
\left \{
\begin{array}{ll}
-\varepsilon \Delta u_{\varepsilon,n} -\trace \cA(\cdot)\rD^{2}u_{\varepsilon,n} +\beta (u_{\varepsilon,n}) =\min\{f,n\}&\quad\hbox{in $\Omega$},\\ [.15cm]
u_{\varepsilon,n}\equiv n &\quad\hbox{in $\partial \Omega$},\\
\end{array}
\right .
\label{eq:ellipticregupron}
\end{equation}
where $f$ is a nonnegative uniformly continuous function. From classical argument,  there exists $u_{\varepsilon,n}\in\cC^{2}(\Omega)$ solving \eqref{eq:ellipticregupron} (see \cite{Ev}, \cite{DL} or \cite{Sa}).  In fact, $\{u_{\varepsilon,n}\}_{n}$ is an increasing sequence.  The consideration of the assumption \eqref{eq:KellerOsserman} enables us to obtain universal bounds as it is proved in Section \ref{sec:universalbound} below. More precisely under \eqref{eq:KellerOsserman} we define 
\begin{equation}
\Psi _{\varepsilon}(\zeta)\doteq \dfrac{1}{\rA}\Phi^{-1}\left (\dfrac{\rB}{\sqrt{\Lambda+\varepsilon}}\zeta\right ),\quad 0\le \zeta <\rR,
\label{eq:universalfunction}
\end{equation}
with $\rR<\dfrac{\sqrt{\Lambda+\varepsilon}}{\rB}\Phi(0^{+})\le +\infty$, where
$$
\Phi(t)=\int^{+\infty}_{t}\dfrac{ds}{\sqrt{2\rG(s)}}.
$$
was defined in Introduction (see \eqref{eq:profilephi}).
\begin{rem}\rm
For the power like case $\beta _{m}(t)=t^{m},~m>1,$ one satisfies $
\Phi_{m}(0^{+})=+\infty$.
Also for $\beta (t)=e^{t}$ and for which $\beta (t)=te^{2t}$ one has
$\Phi (0^{+})=+\infty.\fin $
\label{rem:exemplesbeta3}
\end{rem}

A simple consequence of the proof of Theorem \ref{theo:universalestimate} below in Section \ref{sec:universalbound} we prove
\begin{theo} Assume \eqref{eq:ellipticity}, \eqref{eq:KellerOsserman} as well as \eqref{eq:cA}. Then there exists two positive constants $\rA$ and $\rB$, independent on $n$, for which
\begin{equation}
u_{\varepsilon,n}(x)\le \sum_{i=1}^{\rN}\Psi_{\varepsilon}\big (\rR^{2}-|x_{i}|^{2}\big )+\beta^{-1}\big (\n{f}_{\cQ_{\rR}(x_{0})}\big ),\quad x\in \cQ_{\rR}(x_{0})
\label{eq:universalestimateepsilonn}
\end{equation}
holds in any cube $\cQ_{\rR}(x_{0})\doteq \{x\in\R^{\rN}:~|x_{i}-x_{0,i}|<\rR\}
\subset \subset \Omega .\fin$
\label{theo:universalestimateepsilonn}
\end{theo}
So that, the elliptic regularization implies 
\begin{equation}
u_{\varepsilon}(x)=\sup_{n}u_{\varepsilon,n}(x)<+\infty ,\quad x\in \Omega,
\label{eq:minimalsolutionellipticregularozation}
\end{equation}
whence one proves that $u_{\varepsilon}$ is a $\cC^{2}$ large solution of 
\begin{equation}
	-\varepsilon \Delta u_{\varepsilon} -\trace \cA(\cdot)\rD^{2}u_{\varepsilon} +\beta (u_{\varepsilon}) =f\quad\hbox{in $\Omega$}.
	\label{eq:ellipticequ}
\end{equation}
We send to \cite{DL} for details. We note by construction  that the function defined in  \eqref{eq:minimalsolutionellipticregularozation} is the minimal large solution of \eqref{eq:ellipticequ}.
A refinement in the proof of Theorem \ref{theo:universalestimate} leads to a simple consequence  
\begin{theo}
Assume \eqref{eq:ellipticity}, \eqref{eq:KellerOsserman} as well as \eqref{eq:cA}. Then there exists two positive constants $\rA$ and $\rB$, independent on $n$ for which
\begin{equation}
u_{\varepsilon}(x)\le \sum_{i=1}^{\rN}\Psi_{1}\big (\rR^{2}-|x_{i}|^{2}\big )+\beta^{-1}\big (\n{f}_{\cQ_{\rR}(x_{0})}\big ),\quad x\in \cQ_{\rR}(x_{0})
\label{eq:universalestimateepsilon}
\end{equation}
holds in any cube $\cQ_{\rR}(x_{0})\subset \subset \Omega $, where we are assuming $0<\varepsilon <1.\fin$
\label{theo:universalestimateepsilon}
\end{theo}
Now the estimate \eqref{eq:universalestimateepsilon} enables us to  define the functions 
$$
\overline{u}(x)=\limsup_{y\mapsto x\atop \varepsilon\searrow 0}u_{\varepsilon}(y)
$$
and
$$
\underline{u}(x)=\liminf_{y\mapsto x\atop \varepsilon\searrow 0}u_{\varepsilon}(y).
$$
Moreover, the stablity results of the Viscosity Solution Theory (see \cite{CIL}) show that $\overline{u}(x)$ and $\underline{u}(x)$ 
are respectively viscosity sub and supersolution of
$$
-\trace \cA(\cdot)\rD^{2}u +\beta (u) =f\quad\hbox{in $\Omega$}
$$
(see also \cite{BP}, \cite{BB}).  Both functions blow up on the boundary.  Defining the set
$$
\cF_{-}=\big \{u\hbox{ is a (viscosity) subsolution in $\Omega$} \big \}.
$$
one has $\overline{u}\in \cF_{-}$. Therefore the Perron function
\begin{equation}
	\cU(x)\doteq \sup _{u\in\cF_{-}}u^{*}(x),\quad x\in\Omega
	\label{eq:Perronfunction}
\end{equation}
is well defined under the assumptions of Theorem  \ref{theo:universalestimate} and verifies
$$
\overline{u}(x)\le \cU(x)\le +\infty,\quad x\in\Omega
$$
Futhermore, Perron function, $\cU$,  blows up on the boundary. In fact, adapting the reasonings of \cite[Theorem 4.1]{CIL} (see also \cite[Theorem 4]{D}) the Perron function solves 
$$
-\trace \cA(\cdot)\rD^{2}\cU +\beta (\cU) =f\quad\hbox{in $\Omega$}.
$$
A suitable adaptation of the reasoning of the proof of Theorem \ref{theo:comparisonprincipleviscosity} gives
$$
u_{\varepsilon,n}(x)\le v_{*}(x),\quad x\in\Omega
$$
for any large supersolution of \eqref{eq:fullyequation}, the Perron solution is the minimum large solution of \eqref{eq:fullyequation} in the sense
\begin{equation}
\cU^{*}(x)\le v_{*}(x),\quad x\in\Omega.
\label{eq:minimumsolution}
\end{equation}
We collect the result.
\begin{theo}Let us assume  $\Omega $ is a bounded set in $\R^{\rN},~\rN>1,~\partial \Omega\in\cC^{2}$, \eqref{eq:ellipticity}, \eqref{eq:ellipticityneraboundaty}, \eqref{eq:lightellipticity}, \eqref{eq:structureannulusSMPsub}, \eqref{eq:crecimientog} as well as \eqref{eq:cA}. 
Let $f$ be a uniformly continuous function satisfying 
$$
\limsup_{\rd(x)\searrow 0}\dfrac{f(x)}{\beta \big (\phi \left  (\sqrt{\dfrac{1-\ell}{\Lambda}}\rd(x)\right )}\le \ell \le
\liminf_{\rd(x)\searrow 0}\dfrac{f(x)}{\beta \big (\phi \left  (\sqrt{\dfrac{1-\ell}{\lambda}}\rd(x)\right )}
$$
with $0\le \ell <1$.  Then  $\cU\in\cC(\Omega)$  is the unique large solution of the equation \eqref{eq:fullyequation} whose boundary behaviour satisfies the inequality 
$$
\limsup_{\rd(x)\searrow 0}\dfrac{\cU(x)}{\phi \left (\sqrt{\dfrac{1-\ell}{\Lambda}}\rd(x)\right )}\le 1\le \liminf_{\rd(x)\searrow 0}\dfrac{\cU(x)}{\phi \left (\sqrt{\dfrac{1-\ell}{\lambda}}\rd(x)\right )}
$$
{\rm (}see Theorem \ref{theo:absoluteboundarycomparison} {\rm )}.  Here $\phi$ is the function defined in \eqref{eq:profilephi}. 
\label{theo:existenceKOsources}
\end{theo}
\proof Notice that we do not consider the assumption \eqref{eq:crecimientog1}. Let $v$ be an arbitrary continuous large solution of \eqref{eq:fullyequation}. Then 
$$
\cU^{*}(x)\le v(x),\quad x\in\Omega
$$
holds (see \eqref{eq:minimumsolution}). We claim that in fact $\cU_{*}\equiv v$. Indeed,  assume  that there exists $x_{0}\in\Omega$ such that $\cU_{*}(x_{0})\le v(x_{0})$. Since $\Omega$ is bounded we deduce that
$$
\cO ^{\varepsilon}=\{x\in\Omega:~(1+\varepsilon)\cU_{*}(x)<v(x)\}\subset\subset \Omega
$$
is a non empty open subset, provided $\varepsilon>0$ small. The assumption 
\eqref{eq:crecimientog} implies the property
$$
\dfrac{\beta(t)}{t}\quad \hbox{is increasing for $t$ greater than some large $t_{0}$}
$$
and choose $\mu>0$ so small that $\cU_{*}(x)\ge t_{0}$ holds in  $\Omega_{\mu}=\{x\in\Omega:~\rd(x)<\mu\}$. As $\varepsilon$ is small $\cO ^{\varepsilon}_{\mu}=\cO ^{\varepsilon}\cap \Omega_{\mu}$ is a non empty open subset. Moreover
$$
-\trace \cA(\cdot)\rD^{2}(1+\varepsilon)\cU+\beta \big ((1+\varepsilon)\cU\big )
\ge (1+\varepsilon)\big (-\trace \cA(\cdot)\rD^{2}\cU+\beta (\cU)\big )\ge 0,
\quad x	\in \cO ^{\varepsilon}_{\mu}
$$
in the viscosity sense. Then adapting Remark \ref{eq:directcomparisonviscosity} we obtain
$$
v(x)-(1+\varepsilon)\cU_{*}(x)\le \sup_{\partial \cO ^{\varepsilon}_{\mu}}\big (v-(1+\varepsilon)\cU_{*}\big ),\quad x\in \cO ^{\varepsilon}_{\mu}.
$$
By construction the maximum of $v-(1+\varepsilon)\cU_{*}$ can not be achieved in $\partial \cO ^{\varepsilon}$, whence
$$
v(x)-(1+\varepsilon)\cU_{*}(x)\le \sup_{\{y\in \partial \cO ^{\varepsilon},~\rd(y)=\mu\}}\big (v(y)-\cU_{*}(y)\big ),\quad x\in \cO ^{\varepsilon}_{\mu}
$$
and
$$
v(x)-\cU_{*}(x)\le \sup_{y\in \partial \cO ^{\varepsilon},~\rd(y)=\mu\}}\big (v(y)-\cU_{*}(y)\big )\equiv \rC_{\mu},\quad x\in \cO ^{\varepsilon}_{\mu}
$$
hold. On the other hand, by mean of Theorem \ref{theo:withoutcoercivityII} the inequality $\beta (\cU^{*})\le \beta (v)$ in $\Omega$ gives
$$
v(x)-\cU_{*}(x)\le \rC_{\mu},\quad x\in \Omega\setminus\Omega_{\mu},
$$
thus
$$
v(x)-\cU_{*}(x)\le \rC_{\mu},\quad x\in \Omega.
$$
Since $\rC_{\mu}=v(z_{0})-\cU_{*}(z_{0})\ge 0$ for some $z_{0}\in\partial \cO^{\varepsilon}$ and $\rd(z_{0})=\mu>0$, the upper semi-continuous function
$u-\cU_{*}$ attains a maximim value in $\Omega$ at some interior point. Since $\cU\le v$ we may apply The Strong Maximum Theorem without coercive term (Theorem \ref{theo:SMPwithoutcoercivity}) and deduce
$$
v(x)-\cU_{*}(x)\equiv \rC_{\mu},\quad x\in \Omega.
$$
The relative version of \eqref{eq:capitalinequality} for $v$ and $\cU_{*}$ at $z_{0}\in\Omega$ becomes
$$
\beta \big (v(x_{\alpha})\big )-\beta \big (\cU_{*}(y_{\alpha})\big )\le \omega\big (\alpha|x_{\alpha}-y_{\alpha}|^{2}\big )+\omega_{f}\big (|x_{\alpha}-y_{\alpha}|\big )
$$
whence  letting $\alpha\nearrow \infty$ we get
$$
\beta \big (\cU_{*}(z_{0})+\rC_{\mu}\big )=\beta \big (v(z_{0})\big )\le \beta \big (\cU_{*}(z_{0})\big )
$$
implies $\rC_{\mu}=0$. Notice that the inequality
$$
v=\cU_{*}\le \cU^{*}\le v
$$
concludes the continuity of $\cU.\fin$
\par
A similar program can be developed when the source term has the extra Keller-Osserman type growth studied in Section \ref{sec:KOgrowthextra}. Then one deduces
\begin{theo}Let us assume  $\Omega \in\cC^{2},~\beta\in\cC^{2}$,  \eqref{eq:ellipticity}, \eqref{eq:ellipticityneraboundaty}, \eqref{eq:lightellipticity}, \eqref{eq:structureannulusSMPsub}, \eqref{eq:crecimientog} as well as \eqref{eq:cA}. 
Let $f$ be a uniformly continuous function satisfying 
$$
\lim_{\rd(x)\searrow 0} \dfrac{f(x)}{h\big (\rd(x)\big )}=1,
$$
for some $\cC^{2}$ decreasing function $h:~\R_{+}\rightarrow \R_{+} $ such that $h(0)=\infty$ and $h'(0)=0$ satisfying the compatiblity condition \eqref{eq:compatibility} as well as
$$
\left \{
\begin{array}{l}
\rH_{1}(\sigma)\doteq \dfrac{h(\sigma)}{\beta'\big ( \beta^{-1}\big (h(\sigma)\big )}\left (h''(\sigma)-\dfrac{\big (h'(\sigma )\big )^{2}\beta ''\big ( \beta ^{-1}\big (h(\sigma)\big )\big )}{\big (\beta '\big (\beta ^{-1}\big (h(\sigma)\big )\big )^{2}}\right ),\\ [.2cm]
\rH_{2}(\delta)\doteq -\dfrac{h'(\sigma)}{h(\sigma)\beta '\big (\beta^{-1}\big (h(\sigma)\big )\big )}
\end{array}
\right .
$$
are increasing, vanishing  at the origin {\rm (}see \eqref{eq:criticalgrowthextraexpoademas}{\rm )}. Then the equation \eqref{eq:fullyequation} has a unique large continuous solution $\cU\in\cC(\Omega)$ whose boundary behaviour satisfies 
$$
\limsup_{\rd(x)\searrow 0}\dfrac{\cU(x)}{\beta^{-1}\big (h(\rd(x))\big )}= 1.
$$
\fineq
\label{theo:existenceextraKOsources}
\end{theo}

\setcounter{equation}{0}
\section{Second order term in the boundary assymptotic expansion}
\label{sec:secondorder}
As it was pointed out in Introduction, in Sections \ref{sec:KOgrowth} and \ref{sec:KOgrowthextra} we only have obtained the leading term of  the boundary assymptotic expansion of the large solutions.  Here we study the second order term of this expansion. Relative to Keller-Osserman type growth of the sources on the boundary (see \eqref{eq:criticalgrowthd}) we replace the classical large subsolution near the boundary $\rU_{0}$, used in the proof of Proposition \ref{prop:maximalboundaryprofile}, by 
$$
\rU(x)=\rU_{0}(x)+\rU_{1}(x)\doteq \phi (\sigma)+
\rCh(x,\varepsilon,\delta)\Theta\big(\phi (\sigma)\big ),\quad x\in\Omega^{\delta}_{\delta_{1}}\quad (0<\delta<\delta_{1}<\delta_{\Omega})
$$
where $\sigma=\widehat{\sigma}_{\varepsilon}\big (\rd(x)-\delta)$ with
$ \widehat{\sigma}_{\varepsilon}=\sqrt{\dfrac{(1-\varepsilon)(1-\ell_{u})}{\Lambda}}$. The functions $\rCh(\cdot,\varepsilon,\delta):\Omega_{\delta_{1}}\rightarrow \R$ and $\Theta:\R_{+}\rightarrow \R_{+}$ will be determined from the reasoning that follows. In particular, we will suppose $\rCh(\cdot,\varepsilon,\delta)\in\rW^{2,\infty}\big (\Omega^{\delta}_{\delta_{1}}\big )$. On the other hand, let  $\Theta(t)$ be a function satisfying
\begin{equation}
\Theta'(t)=-t\Phi '(t)\quad \hbox{for large $t$}.
\label{eq:phi1}
\end{equation}
See Remark \ref{rem:phi1} below for suitable choices of the function $\Theta$.
Clearly, if $\Theta \big (\phi (0)\big )=\Theta (\infty)<\infty$ one has
\begin{equation}
\lim_{\sigma \rightarrow 0}\dfrac{\Theta\big(\phi (\sigma)\big )}{\phi (\sigma)}
=\lim_{t \rightarrow \infty}\dfrac{\Theta(t)}{t}=0.
\label{eq:boundedsecondorder}
\end{equation}
it implies that $\rCh(x,\varepsilon,\delta)\Theta\big(\phi (\sigma)\big )$ is a bounded second order term in the explosive boundary expansion. When $\Theta \big (\phi (0)\big )=\Theta (\infty)=\infty$ the property \eqref{eq:phi1} implies that $\rCh(x,\varepsilon,\delta)\Theta\big(\phi (\sigma)\big )$ is an explosive second order term in the expansion. Indeed, by using the L'Hôpital rule and \eqref{eq:fromKellerOsserman2} one deduces
\begin{equation}
\lim_{\sigma \rightarrow 0}\dfrac{\Theta\big(\phi (\sigma)\big )}{\phi (\sigma)}
=\lim_{t \rightarrow \infty}\dfrac{\Theta(t)}{t}=\lim_{t \rightarrow \infty}\Theta'(t)=\lim_{t \rightarrow \infty}\dfrac{t}{\sqrt{2\rG(t)}}=0.
\label{eq:explosivesecondorder}
\end{equation}
Let $0<\varepsilon <1$ and $0\le \ell_{u}<1$. Denoting 
$\rCh(x)=\rC(x,\varepsilon,\delta)$,  one has
$$
\begin{array}{ll}
\rD_{ij}\rU_{1}(x)&\hspace*{-.2cm}=\Theta\big(\phi (\sigma)\big )\rD_{ij}\rCh(x)+\sigma_{\varepsilon}\Theta'\big(\phi (\sigma)\big ) \phi '(\sigma)\big (\rD_{i}\rCh(x)\rD_{j}\rd(x)+\rD_{j}\rCh(x)\rD_{i}\rd(x)\big )\\ [.15cm]
&+\widehat{\sigma}_{\varepsilon}^{2}\rCh(x)\big (\Theta ''\big(\phi (\sigma)\big )\big(\phi '(\sigma)\big )^{2}+\Theta '\big(\phi (\sigma)\big )\phi ''(\sigma)\big )\rD_{i}\rd(x)\rD_{j}\rd(x)\\ [.15cm]
&+\widehat{\sigma} _{\varepsilon}\rCh(x) \Theta'\big (\phi (\sigma)\big )\phi '(\sigma)\rD_{ij}\rd(x),\quad x\in\Omega^{\delta}_{\delta_{1}},
\end{array}
$$
whence
$$
\begin{array}{ll}
\trace \cA(x)\rU_{1}(x)&\hspace*{-.2cm}=\Theta\big(\phi (\sigma)\big )\trace \cA(x)\rD^{2}\rCh(x)\\ [.15cm]
&+\sigma_{\varepsilon}\Theta '\big(\phi (\sigma)\big ) \phi '(\sigma)\trace \cA(x)\big (\rD\rCh(x)\otimes\rD\rd(x)+\rD\rd(x)\otimes \rD\rCh(x)\big )\\ [.15cm]
&+\widehat{\sigma}_{\varepsilon}^{2}\rCh(x)\big (\Theta ''\big(\phi (\sigma)\big )\big(\phi '(\sigma)\big )^{2}+\Theta '\big(\phi (\sigma)\big )\phi ''(\sigma)\big )\rE(x)\\ [.15cm]
&+\widehat{\sigma} _{\varepsilon}\rCh(x) \Theta' \big (\phi (\sigma)\big )\phi '(\sigma)\trace \cA(x)\rD^{2}\rd(x),\quad x\in\Omega^{\delta}_{\delta_{1}},
\end{array}
$$
where $\rE(x)\doteq \trace \cA(x)\rD \rd(x)\otimes \rD \rd(x)\in [\lambda,\Lambda]$ is an eikonal-like term.
On the other hand
$$
\begin{array}{ll}
\beta (\rU(x))&\hspace*{-.2cm}=\beta (\rU_{0}(x))+\rU_{1}(x)\beta ' (\rU_{0}(x))+
\dfrac{\rU_{1}(x)^{2}}{2}\beta ''(\xi)\\ [.15cm]
&=\beta \big (\phi(\sigma)\big )+
\rCh(x)\Theta \big (\phi(\sigma)\big )\beta '\big (\phi(\sigma)\big )+
\dfrac{\big (\rCh(x)\Theta \big (\phi(\sigma)\big )\big )^{2}}{2}\beta ''\big (\xi(\sigma)\big ),
\end{array}
$$
for some  $\xi(\sigma)$ between $\phi(\sigma)$ and $\phi(\sigma) +\rCh(x)\Theta \big (\phi(\sigma)\big)$. We are interested in to prove the inequality
$$
-\trace \cA(x)\rD^{2}\rU(x)+\beta \big (\rU(x)\big )-f(x)\ge 0,\quad x\in\Omega^{\delta}_{\delta_{1}}
$$
for $\delta_{1}$ small enough. The computations  are very tedious. So that, we consider suitable groups of terms. Moreover, we will use the change of variable $t=\phi(\sigma)$. As in the proof of Proposition \ref{prop:maximalboundaryprofile} the first group is
$$
\cG_{1}\doteq -\trace \cA(x)\rD^{2}\rU_{0}(x)+\beta \big (\rU_{0}(x)\big )- f(x)+
\sigma _{0}\phi' (\sigma)\trace \cA\rD^{2}\rd\ge  \beta (t) \big (\varepsilon (1-\ell_{u})\big )>0
$$
provided \eqref{eq:criticalgrowthd} (see \eqref{eq:stepfirst}). The second group is
$$
\begin{array}{ll}
\cG_{2}&\hspace*{-.2cm}\doteq -\widehat{\sigma} _{\varepsilon}\phi' (\sigma)\trace \cA\rD^{2}\rd(x)\\
[.15cm]&-\rCh (x)\big [\widehat{\sigma}_{\varepsilon}^{2}\big (\Theta ''\big(\phi (\sigma)\big )\big(\phi '(\sigma)\big )^{2}+\Theta '\big(\phi (\sigma)\big )\phi ''(\sigma)\big )\big )\rE(x)-\Theta \big (\phi(\sigma)\big )\beta '\big (\phi(\sigma)\big )\big ]
\end{array}
$$
for the particular case $t=\phi(\sigma)$
Then 
$$
\begin{array}{ll}
\cG_{2}&\hspace*{-.2cm}= \widehat{\sigma} _{\varepsilon}\trace \cA\rD^{2}\rd\sqrt{2\rG(t)}-\rCh (x)\big [\widehat{\sigma} _{\varepsilon}^{2}\big (2\Theta ''(t)\rG(t)+\Theta '(t)\beta(t)\big )\rE(x)-\Theta (t)\beta '(t)\big ].
\end{array}
$$
From \eqref{eq:phi1} we have
$$
\Theta'(t)=\dfrac{t}{\sqrt{2\rG(t)}}
$$
and
$$
2\Theta ''(t)\rG(t)+\Theta '(t)\beta (t)=\sqrt{2\rG(t)}.
$$
So, the choice
\begin{equation}
\rCh(x)=\widehat{\sigma}_{\varepsilon}\dfrac{\sqrt{2\rG\big (\phi(\sigma)\big )}}{\widehat{\sigma}_{\varepsilon}^{2}\sqrt{2\rG\big (\phi(\sigma)\big )}\rE(x)-\Theta \big (\phi(\sigma)\big )\beta '\big (\phi(\sigma)\big )}
\trace \cA(x)\rD^{2}\rd(x),
\label{eq:Aaprox}
\end{equation}
for $x\in \Omega^{\delta}_{\delta_{1}}$, leads $\cG_{2}=0$. 
\begin{rem}\rm  We may write \eqref{eq:Aaprox} as
\begin{equation}
\begin{array}{ll}
\rCh(x)&\hspace*{-.2cm}=\dfrac{\widehat{\sigma}_{\varepsilon}}{\widehat{\sigma}_{\varepsilon}^{2}\rE(x)-\dfrac{\Theta (t)\beta '(t)}{\sqrt{2\rG(t)}}}
\trace \cA(x)\rD^{2}\rd(x)\\ [1cm]
&\hspace*{-.2cm}=-\widehat{\sigma}_{\varepsilon}\dfrac{\sqrt{2\rG(t)}}{\Theta (t)\beta '(t)}\dfrac{1}{1-\widehat{\sigma}_{\varepsilon}^{2}\rE(x)\dfrac{\sqrt{2\rG(t)}}{\Theta (t)\beta '(t)}}
\trace \cA(x)\rD^{2}\rd(x)\\ [1cm]
&\hspace*{-.2cm}=-\widehat{\sigma}_{\varepsilon}\dfrac{\sqrt{2\rG(t)}}{\Theta (t)\beta '(t)}\disp \sum_{n\ge }\left (\widehat{\sigma}_{\varepsilon}^{2}\rE(x)\dfrac{\sqrt{2\rG(t)}}{\Theta (t)\beta '(t)}\right )^{n}
\trace \cA(x)\rD^{2}\rd(x)\\ [1cm]
&\hspace*{-.2cm}=-\widehat{\sigma}_{\varepsilon}\dfrac{\sqrt{2\rG\big (\phi(\sigma)\big )}}{\Theta \big (\phi(\sigma)\big )\beta '\big (\phi(\sigma)\big )}\trace \cA(x)\rD^{2}\rd(x)\big (1+o(1)\big)
\trace \cA(x)\rD^{2}\rd(x),
\end{array}
\label{eq:restoA}
\end{equation}
provided 
\begin{equation}
\widehat{\sigma}_{\varepsilon}^{2}\rE(x)\dfrac{
\sqrt{2\rG\big (\phi(\sigma)\big )}}{\beta '\big (\phi(\sigma)\big )}\dfrac{1}{\Theta \big (\phi(\sigma)\big )\big )}<1.
\label{eq:condicionrestoA}
\end{equation}
Analogously
\begin{equation}
\begin{array}{ll}
\rCh(x)\Theta \big (\phi(\sigma)\big)&\hspace*{-.2cm}=
\dfrac{\widehat{\sigma}_{\varepsilon}\Theta(t)}{\widehat{\sigma}_{\varepsilon}^{2}\rE(x)-\dfrac{\Theta (t)\beta '(t)}{\sqrt{2\rG(t)}}}
\trace \cA(x)\rD^{2}\rd(x)\\ [1cm]
&=-\dfrac{\widehat{\sigma}_{\varepsilon}\sqrt{2\rG(t)}}{\beta'(t)}\dfrac{1}{1-\widehat{\sigma}_{\varepsilon}^{2}\rE(x)\dfrac{\beta '(t)}{
\sqrt{2\rG(t)}}\left (\Theta \big (\phi(\sigma)\big )\right )^{-1}}\trace \cA(x)\rD^{2}\rd(x)\\ [1cm]
&\disp =-\dfrac{\widehat{\sigma}_{\varepsilon}\sqrt{2\rG(t)}}{\beta'(t)}
\sum_{n\ge 0}\bigg (\widehat{\sigma}_{\varepsilon}^{2}\rE(x)\dfrac{\beta '(t)}{
		\sqrt{2\rG(t)}}\left (\Theta (t)^{-1}\right )^{n}\trace \cA(x)\rD^{2}\rd(x)\\
[.75cm]&=-\dfrac{\widehat{\sigma}_{\varepsilon}\sqrt{2\rG\big (\phi(\sigma)\big )}}{\beta' \big (\phi(\sigma)\big )}\trace \cA(x)\rD^{2}\rd(x)\big (1+o(1)\big ),
\label{eq:resto}
\end{array}
\end{equation}
provided 
\begin{equation}
\widehat{\sigma}_{\varepsilon}^{2}\rE(x)\dfrac{\beta '\big (\phi(\sigma)\big )}{
\sqrt{2\rG\big (\phi(\sigma)\big )}}\dfrac{1}{\Theta \big (\phi(\sigma)\big )\big )}<1.
\label{eq:condicionresto}
\end{equation}
\par
\noindent
We note that the claim implies that if such secondary term exists it is unique term  influenced by the geometry by means of $\trace \cA(x)\rD^{2}\rd(x)$ whenever $\Lambda=\lambda=1.\fin$
\end{rem}
The third group is
$$
\cG_{3}\doteq -\Theta\big (\phi(\sigma)\big)\trace \cA(x)\rD^{2}\rCh(x)=-\Theta(t)\trace \cA(x)\rD^{2}\rCh(x)
$$
for which
$$
\lim_{t \nearrow \infty} \dfrac{\Theta(t)}{\beta (t)}=
\lim_{t \nearrow \infty} \dfrac{\Theta(t)}{t}\frac{t}{\beta (t)}=0
$$
(see \eqref{eq:fromKellerOsserman2} and  \eqref{eq:boundedsecondorder} or \eqref{eq:explosivesecondorder}). The next group to be considered is 
$$
\cG_{4}\doteq \dfrac{1}{2}\big (\rCh(x)\Theta\big (\phi(\sigma)\big )\big )^{2}\beta ''\big (\xi(\sigma)\big ).
$$
Since $\xi (\sigma )\rightarrow \infty$ as $\sigma \rightarrow 0$ and $\disp \lim_{t\nearrow \infty}\dfrac{\sqrt{2\rG(t)}}{\beta (t)}=0$ (see \eqref{eq:fromKellerOsserman}), we will use the property
\begin{equation}
\lim_{\sigma\searrow 0}\dfrac{\big (\rCh(x)\Theta\big (\phi(\sigma)\big )\big )^{2}\beta ''\big (\phi(\sigma)\big )}{\sqrt{2\rG\big (\phi(\sigma)\big )}}
\dfrac{\sqrt{2\rG\big (\phi(\sigma)\big )}}{\beta \big (\phi(\sigma)\big )}=0
\label{eq:techpreviabeta2}
\end{equation}
deduced from the assumption
\begin{equation}
\lim_{\sigma\searrow 0}\dfrac{\big (\rCh(x)\Theta\big (\phi(\sigma)\big )\big )^{2}\beta ''\big (\phi(\sigma)\big )}{\sqrt{2\rG\big (\phi(\sigma)\big )}}=0.
\label{eq:techbeta2}
\end{equation}
The last group is 
$$
\begin{array}{ll}
\cG_{5} &\hspace*{-.2cm}\doteq -\sigma_{\varepsilon}\Theta '\big(\phi (\sigma)\big ) \phi '(\sigma)\left  [\trace \cA(x)\big (\rD\rCh(x)\otimes\rD\rd(x)+\rD\rd(x)\otimes \rD\rCh(x)\big ) \right .\\ [.15cm]
&\quad \left .+\rCh(x)\trace \cA(x)\rD^{2}\rd(x)\right  ]\\  [.15cm]
&\hspace*{-.2cm}=\sigma_{\varepsilon}\Theta '(t) \sqrt{2\rG(t)}\left  [\trace \cA(x)\big (\rD\rCh(x)\otimes\rD\rd(x)+\rD\rd(x)\otimes \rD\rCh(x)\big ) +\rCh(x)\trace \cA(x)\rD^{2}\rd(x)\right  ]\\  [.15cm]
&\hspace*{-.2cm}=-\sigma_{\varepsilon}t\left  [\trace \cA(x)\big (\rD\rCh(x)\otimes\rD\rd(x)+\rD\rd(x)\otimes \rD\rCh(x)\big ) +\rCh(x)\trace \cA(x)\rD^{2}\rd(x)\right  ] 
\end{array}
$$
for which we recall again the property
$$
\lim_{t\nearrow \infty} \dfrac{t}{\beta (t)}=0
$$
(see \eqref{eq:fromKellerOsserman2}). With the above reasonings we have
\begin{theo}[Second order term] Let us assume $\partial \Omega \in \cC^{4},~\beta \in \cC^{2}$ and $\cA(\cdot)$ adequately smooth. Under the assumptions of Proposition \ref{prop:maximalboundaryprofile}  
for any solution of 
$$
-\trace \cA(\cdot)\rD^{2}u+\beta (u)\le f\quad \hbox{in $\Omega$}
$$
we have
\begin{equation}
\limsup_{\rd (x)\searrow 0}\dfrac{u^{*}(x)}{\phi \big (\widehat{\sigma}_{0}\rd(x)\big )+\rC_{u}(x)\Theta\left (\phi \big (\widehat{\sigma}_{0}\rd(x)\big )\right )}\le 1
\label{eq:maximalboundaryprofilesecondKO}
\end{equation}
provided \eqref{eq:criticalgrowthd} and   \eqref{eq:techbeta2}, where 
$\widehat{\sigma}_{0}=\sqrt{\dfrac{1-\ell}{\Lambda}}$ and 
the function $\Theta$ is given by \eqref{eq:phi1} and
$$
\rC_{u}(x)=\dfrac{\Lambda\widehat{\sigma}_{0}}{(1-\ell_{u})\rE(x)-\Lambda\dfrac{\Theta \left (\phi \big (\widehat{\sigma}_{0}\rd(x)\big )\right)}{\sqrt{2\rG\left (\phi \big (\widehat{\sigma}_{0}\rd(x)\big )\right)}}\beta '\big (\phi \big (\widehat{\sigma}_{0}\rd(x)\big )}
\trace \cA(x)\rD^{2}\rd(x),
$$
for $\rd(x)$ small {\rm (}see \eqref{eq:resto}{\rm )},  is assumed in $\rW^{2,\infty}$ and $\rE(x)\doteq\trace \cA(x)\rD\rd(x)\otimes\rD\rd(x)$. Moreover
\begin{equation}
\rC_{u}(x)\Theta \big (\phi(\widehat{\sigma}_{0}\rd(x))\big )=-\dfrac{\widehat{\sigma}_{0}\sqrt{2\rG\left (\phi \big (\widehat{\sigma}_{0}\rd(x)\big )\right)}}{\beta '\big (\phi \big (\widehat{\sigma}_{0}\rd(x)\big )}\trace \cA(x)\rD^{2}\rd(x)\big (1+o(1)\big ),
\label{eq:restou}
\end{equation}
for $\rd(x)$ small {\rm (}see \eqref{eq:resto}{\rm )}. Analogously, under the assumptions of Proposition \ref{prop:minimalboundaryprofile}  
for any nonnegative large solution of 
$$
-\trace \cA(\cdot)\rD^{2}v+\beta (v)\ge f\quad \hbox{in $\Omega$}
$$
we have
\begin{equation}
1\le \liminf_{\rd (x)\searrow 0}\dfrac{v_{*}(x)}{\phi \big (\widetilde{\sigma}_{0}\rd(x)\big )+\rC_{d}(x)\Theta\left (\phi \big (\widetilde{\sigma}_{0}\rd(x)\big )\right )}
\label{eq:minimalboundaryprofilesecondKO}
\end{equation}
provided \eqref{eq:criticalgrowthu} and   \eqref{eq:techbeta2}, where 
$\widetilde{\sigma}_{0}=\sqrt{\dfrac{1-\ell}{\lambda}}$
the function $\Theta$ is given by \eqref{eq:phi1} and 
$$
\rC_{d}(x)=\dfrac{\lambda \widetilde{\sigma}_{0}}{(1-\ell_{d})\rE(x)-\lambda \dfrac{\Theta \left (\phi \big (\widetilde{\sigma}_{0}\rd(x)\big )\right )}{\sqrt{2\rG\left (\phi \big (\widetilde{\sigma}_{0}\rd(x)\big )\right )}}\beta '\big (\phi \big (\widetilde{\sigma}_{0}\rd(x)\big )}
\trace \cA(x)\rD^{2}\rd(x),
$$
for $\rd(x)$ small,  assumed in $\rW^{2,\infty}$. Moreover
\begin{equation}
\rC_{d}(x)\Theta \big (\phi(\widetilde{\sigma}_{0}\rd(x))\big )=-\dfrac{\widehat{\sigma}_{0}\sqrt{2\rG\left (\phi \big (\widetilde{\sigma}_{0}\rd(x)\big )\right)}}{\beta '\big (\phi \big (\widetilde{\sigma}_{0}\rd(x)\big )}\trace \cA(x)\rD^{2}\rd(x)\big (1+o(1)\big ),
\label{eq:restod}
\end{equation}
for $\rd(x)$ small {\rm (}see \eqref{eq:resto}{\rm )}.
\label{theo:secondarytermKOgrowth}
\end{theo}
\proof We only show the result relative to $u^{*}$. The proof  for $v_{*}$ is analogous (see Proposition \ref{prop:minimalboundaryprofile}). Following the above reasonings  the function
$$
\rU(x)=\rU_{0}(x)+\rU_{1}(x)\doteq \phi (\sigma)+
\rCh(x)\Theta\big(\phi (\sigma)\big ),\quad x\in\Omega^{\delta}_{\delta_{1}}\quad (0<\delta<\delta_{1}<\delta_{\Omega})
$$
satisfies
$$
\begin{array}{ll}
-\trace \cA(x)\rD^{2}\rU(x)+\beta \big (\rU(x)\big )-f(x)&\hspace*{-.2cm}\ge \cG_{1}+\cG_{2}+\cG_{3}+\cG_{4}+\cG_{5}\\
&\hspace*{-.2cm}\ge \beta (t) \left (\varepsilon (1-\ell_{u})+\dfrac{\cG_{3}+\cG_{4}+\cG_{5}}{\beta (t)}\right ),\quad x\in \Omega^{\delta}_{\delta_{1}},
\end{array}
$$
for $\Theta(t)$ and $\rCh(x) $ are defined in \eqref{eq:phi1} and \eqref{eq:Aaprox} respectively. Since the assumptions imply
$$
\lim_{t\nearrow \infty}\dfrac{\cG_{3}+\cG_{4}+\cG_{5}}{\beta (t)}=0
$$
and $0<\sigma<\widehat{\sigma}_{0}\delta_{1}$ we deduce
$$
\varepsilon (1-\ell_{u})+\dfrac{\cG_{3}+\cG_{4}+\cG_{5}}{\beta \big (\phi (\sigma)\big )}>0
$$
for $\delta_{1}$ small enough. Therefore
$$
-\trace \cA(x)\rD^{2}\rU(x)+\beta \big (\rU(x)\big )\ge f(x) \quad x\in \Omega^{\delta}_{\delta_{1}}
$$ 
for $\delta_{1}$ small enough.  As in  the proof of Proposition \ref{prop:maximalboundaryprofile} we obtain
$$
u^{*}(x)\le \phi  \left (\sqrt{\dfrac{1-\varepsilon)(1-\ell_{u})}{\Lambda}}\rd(x)\right ) +\rCh(x)\Theta \left ( \phi  \left (\sqrt{\dfrac{1-\varepsilon)(1-\ell_{u})}{\Lambda}}\rd(x)\right )\right )+\rM,\quad x\in \Omega_{\delta_{1}}.
$$
where $\rM\ge \disp \sup_{\rd(x)=\delta_{1}}u^{*}$. Then the estimate \eqref{eq:maximalboundaryprofilesecondKO} follows.$\fin$
\begin{rem}\rm In this paper we omit general assumptions for which $\rC_{u}(x),\rC_{d}(x)\in \rW^{2,\infty}$. The above proof adapts and simplifies results obtained in \cite{ADR} whenever $f\equiv 0$ and $\Lambda =\lambda =1$. See \cite{AP} for other reasonings.$\fin$
\end{rem}
\begin{rem}[Second order boundary assymtotic]\rm From \eqref{eq:resto} we have 
$$
\rC_{u}(x)\Theta \big (\phi(\sigma)\big )=-\widehat{\sigma}_{0}\dfrac{\sqrt{2\rG\big (\phi(\sigma))}}{\beta '\big (\phi(\sigma)\big)}\trace \cA(x)\rD^{2}\rd(x)
\big (1+o(1)\big )
$$
for $\sigma=\widehat{\sigma}_{0}\rd(x)$, with $\widehat{\sigma}_{0}=\sqrt{\dfrac{1-\ell_{u}}{\Lambda}}$, provided \eqref{eq:condicionresto}. Then, denoting
\begin{equation}
\lim_{t\nearrow \infty}\dfrac{\sqrt{2\rG(t)}}{\beta '(t)}=\rL
\label{eq:L}
\end{equation}
we deduce:
\par
\noindent i) If $\rL=\infty$ the second order term $\rC_{u}(x)\Theta \big (\phi(\sigma)\big )$ is explosive on the boundary where the influence of the geometry appears,
\par
\noindent ii) If $0<\rL<\infty$ the second order term $\rC_{u}(x)\Theta \big (\phi(\sigma)\big )$ is bounded and the influence of the geometry appears on the boundary,
\par
\noindent iii) If $\rL=0$ the second order term $\rC_{u}(x)\Theta \big (\phi(\sigma)\big )$ vanishes on the boundary and the influence of the geometry is null on the boundary.
\par
Let us come back to the choice  of function $\Theta$ satisfying \eqref{eq:phi1}:
\par
\medskip
\noindent a) Under 
\begin{equation}
\int^{+\infty}\dfrac{s ds}{\sqrt{2\rG(s)}}=\infty
\label{eq:noKOsecondterm}
\end{equation}
we will consider the  function
\begin{equation}
\Theta(t)=\int^{t}_{t_{0}}\dfrac{s ds}{\sqrt{2\rG(s)}},\quad \hbox{for $t$ large}
\label{eq:phi1explicityexplosive}
\end{equation}
where $t_{0}<\infty$ is arbitrary.  We note that $\Theta (\infty)=\infty$ and $0\le \rL\le \infty$. It is preferable $t_{0}=0$ if 
\begin{equation}
	\int_{0^{+}}\dfrac{s ds}{\sqrt{2\rG(s)}}<+\infty	\label{eq:KOsecondtermatorigin}
\end{equation}
holds. 
\par
\medskip
\noindent b) Under 
\begin{equation}
\int^{+\infty}\dfrac{s ds}{\sqrt{2\rG(s)}}<\infty
\label{eq:noKOsecondterm}
\end{equation}
and
\begin{equation}
\int_{0^{+}}\dfrac{s ds}{\sqrt{2\rG(s)}}=+\infty
\label{eq:noKOsecondtermatorigin}
\end{equation}
we define 
\begin{equation}
\Theta(t)=-\int^{\infty}_{t}\dfrac{s ds}{\sqrt{2\rG(s)}},\quad \hbox{for $t$ large}.
\label{eq:phi1explicitynull}
\end{equation}
We note that $\Theta (\infty)=0$ and $\rL=0$.
\par
\noindent c) Under 
\begin{equation}
\int^{+\infty}\dfrac{s ds}{\sqrt{2\rG(s)}}<\infty
\label{eq:KOsecondterm}
\end{equation}
and
\begin{equation}
\int_{0^{+}}\dfrac{s ds}{\sqrt{2\rG(s)}}<+\infty
\label{eq:KOsecondtermatorigin}
\end{equation}
we define 
\begin{equation}
\Theta(t)=\int^{t}_{0}\dfrac{s ds}{\sqrt{2\rG(s)}},\quad \hbox{for $t$ large}.
\label{eq:phi1explicitybounded}
\end{equation}
We note that $\disp \Theta (\infty)=\int^{\infty}_{0^{+}}\dfrac{s ds}{\sqrt{2\rG(s)}}<+\infty$ and $0<\rL<\infty$.
\par
\noindent
We send to \cite{ADR} for some details on the  explosive boundary expansion of solutions including non explosive second order terms for homogeneous sources $f\equiv 0$ and $\Lambda=\lambda=1.\fin$ 
\label{rem:phi1} 
\end{rem}
\begin{exam}[Power like case]\rm For $\beta_{m}(t)= t^{m},~m>1,$ one has
$$
\sqrt{2\rG_{m}(t)}=\sqrt{\frac{2}{m+1}}t^{\frac{m+1}{2}},\hspace*{.2cm} \Phi_{m}(t)=\dfrac{\sqrt{2(m+1)}}{m-1}t^{-\frac{m-1}{2}}
\hspace*{.2cm} \hbox{and}\hspace*{.2cm}
\phi_{m}(\delta)=\left (\dfrac{\sqrt{2(m+1)}}{m-1}\right )^{\frac{2}{m-1}}\delta^{-\frac{2}{m-1}}
$$
(see Remark \ref{rem:exemplesbeta}).  Let us recall
$$
\Theta_{m}'(t)=-t\Phi_{m}'(t)=\dfrac{t}{\sqrt{2\rG_{m}(t))}}=\sqrt{\frac{m+1}{2}}t^{-\frac{m-1}{2}}
$$
(see \eqref{eq:phi1}) and
\begin{equation}
	\int^{t_{2}}_{t_{1}}\dfrac{sds}{\sqrt{2\rG_{m}(s)}}=\sqrt{\dfrac{m+1}{3}}\dfrac{2}{3-m}\left (t_{2}^{\frac{3-m}{2}}-t_{1}^{\frac{3-m}{2}}\right ).
	\label{eq:primitivesecondorder}
\end{equation}
\par
\noindent {\bf a}) If $1<m<3$ the properties
$$
\int^{\infty}\dfrac{sds}{\sqrt{2\rG_{m}(s)}}=\infty
\quad \hbox{and}\quad \int_{0^{+}}\dfrac{sds}{\sqrt{2\rG_{m}(s)}}<\infty
$$
hold. Then 
$$
\dfrac{\sqrt{2\rG_{m}\big (\phi_{m}(\sigma))}}{\beta _{m}'\big (\phi_{m}(\sigma)\big)}=\sqrt{2}\left (\dfrac{(m+1)^{2-m}}{(m-1)^{3-m}}\right )^{\frac{1}{m-1}}\sigma^{-\frac{3-m}{m-1}}.
$$
On the other hand, we define
\begin{equation}
	\Theta_{m}(t)=\int^{t}_{0}\dfrac{sds}{\sqrt{2\rG_{m}(s)}}=\sqrt{\dfrac{m+1}{2}}\dfrac{2}{3-m}t^{\frac{3-m}{2}}
	\label{eq:phi1power}
\end{equation}
that satisfies $\Theta_{m}(\infty)=\infty$. Since \eqref{eq:resto} and \eqref{eq:condicionresto} hold, we deduce that $\rC_{u,m}(x)\Theta \big (\phi_{m}(\widehat{\sigma}_{0}\rd(x))\big )$ is an explosive second order term and the influence of the geometry appears on the boundary (see Remak \ref{rem:phi1}). More precisely
\begin{equation}
\rC_{u,m}(x)\Theta \big (\phi_{m}(\sigma)\big )=-\widehat{\sigma}_{0}^{\frac{m-2}{m-1}}
\sqrt{2}\left (\dfrac{(m+1)^{2-m}}{(m-1)^{3-m}}\right )^{\frac{1}{m-1}}\rd(x)^{-\frac{3-m}{m-1}}\trace \cA(x)\rD^{2}\rd(x)
\big (1+o(1)\big )
\label{eq:restom}
\end{equation}
for $\sigma=\widehat{\sigma}_{0}\rd(x)$, with $\widehat{\sigma}_{0}=\sqrt{\dfrac{1-\ell_{u}}{\Lambda}}$.  We may detail the explosive boundary expansion. Indeed,  we note that 
\begin{equation}
\Theta_{m}\big (\phi_{m}(\sigma)\big ) =\sqrt{\dfrac{m+1}{2}}\dfrac{2}{3-m}\left (\dfrac{\sqrt{2(m+1)}}{m-1}\right )^{-1+\frac{2}{m-1}}\sigma^{1-\frac{2}{m-1}},
\label{eq:phi1powersigma}
\end{equation}
whence
$$
\begin{array}{l}
\phi_{m}\big (\sigma)+\rC_{u,m}(x)\Theta_{m}\big (\phi_{m}(\sigma)\big )=
\phi_{m}\big (\sigma)\bigg  (1+\\
\hspace*{3cm}\left .+\rC_{u,m}(x)\sqrt{\dfrac{m+1}{2}}\dfrac{2}{3-m}\left (\dfrac{\sqrt{2(m+1)}}{m-1}\right )^{-1}\sigma\trace \cA(x)\rD^{2}\rd(x)\right ).
\end{array}
$$
Moreover, one has
\begin{equation}
\begin{array}{ll}
\rC_{u,m}(x)&\hspace*{-.2cm}=\dfrac{\widehat{\sigma}_{0}}{\widehat{\sigma}_{0}^{2}\rE(x)-\dfrac{\Theta_{m}\big (\phi_{m}(\sigma)\big )}{\sqrt{2\rG_{m}\big (\phi_{m}(\sigma)\big )}}\beta _{m}'\left (\big (\phi_{m}(\sigma)\big )\right )}
\trace \cA(x)\rD^{2}\rd(x)\\ [.8cm]
&=\dfrac{\widehat{\sigma}_{0}(m-3)}{m^{2}+\widehat{\sigma}_{0}^{2}\rE(x)+1\big )m-3\widehat{\sigma}_{0}^{2}\rE(x)}\trace \cA(x)\rD^{2}\rd(x).
\end{array}
\label{eq:Am}
\end{equation}
One proves the property
$$
m^{2}+\widehat{\sigma}_{0}^{2}\rE(x)+1\big )m-3\widehat{\sigma}_{0}^{2}\rE(x)\neq 0
$$
whence $\rC_{u,m}\in\rW^{2,\infty}$, provided  $\partial \Omega \in \cC^{4}$  and $\cA(\cdot)$ is adequately smooth.  
\par
\noindent 
The computation for $\sigma=\widetilde{\sigma}_{0}\rd(x)$, with $\widetilde{\sigma}_{0}=\sqrt{\dfrac{1-\ell_{d}}{\lambda}}$ are analogous.
Finally,  one proves that \eqref{eq:techbeta2} holds and Theorem \ref{theo:secondarytermKOgrowth} applies. We note that the Lowener-Nirenberg choice (see \cite{LN})
$\beta(t)=t^{\frac{\rN+2}{\rN-2}}$ belongs this case if $\rN>4$.
\par
\noindent 
{\bf b)} When $m>3$ the properties 
$$
\int^{\infty}\dfrac{sds}{\sqrt{2\rG_{m}(s)}}<\infty
\quad \hbox{and}\quad \int_{0^{+}}\dfrac{sds}{\sqrt{2\rG_{m}(s)}}=\infty
$$
hold (see \eqref{eq:primitivesecondorder}. Here one defines the  function 
$$
\Theta_{m}(t)=-\int^{\infty}_{t}\dfrac{sds}{\sqrt{2\rG_{m}(s)}}=\sqrt{\dfrac{m+1}{2}}\dfrac{2}{3-m}t^{\frac{3-m}{2}}
$$
that coincides with  \eqref{eq:phi1power}, but now $m>3$ where  the formulae 
\eqref{eq:phi1powersigma} and \eqref{eq:Am} are the same.
Therefore, we deduce that the second order $\rC_{u,m}(x)\Theta _{m}\big (\phi_{m}(\widehat{\sigma}_{0}\rd(x))\big )$ vanishes on the boundary, hence the influence of the geometry is null on the boundary. The Lowener-Nirenberg choice belongs this case if $\rN=3$.
\par
\noindent 
{\bf c)} If $m=3$  the properties 
$$
\int^{\infty}\dfrac{sds}{\sqrt{2\rG_{3}(s)}}<\infty
\quad \hbox{and}\quad \int_{0^{+}}\dfrac{sds}{\sqrt{2\rG_{3}(s)}}=\infty
$$
hold. Then one define
$$
\Theta_{3}(t)=\int^{t}_{t_{0}}\dfrac{s ds}{\sqrt{2\rG_{3}(s)}}=\sqrt{2}\log \bigg (\dfrac{t}{t_{0}}\bigg ),\quad t>t_{0}
$$
where $t_{0}<\infty$ is arbitrary.  Here $\Theta_{3}(\infty)=\infty$, thus $
\Theta_{3}\big (\phi_{3}(0)\big )=+\infty $. Take $t_{0}=1$ by simplicity.
In this case \eqref{eq:restom} implies that the second order $\rA_{u,3}(x)\Theta _{3}\big (\phi_{3}(\widehat{\sigma}_{0}\rd(x))\big )$ is bounded where the influence of the geometry on the boundary appears (see Remak \ref{rem:phi1}).  Then
$$
\Theta _{3}\big (\phi_{3}(\sigma)\big )=\sqrt{2}\log \bigg (\dfrac{\sqrt{2}}{\sigma}\bigg ).
$$
and
$$
\begin{array}{ll}
\rC_{u,3}(x)&\hspace*{-.2cm}=\dfrac{\widehat{\sigma}_{0}}{\widehat{\sigma}_{0}^{2}\rE(x)-\dfrac{\Theta_{3}\big (\phi_{3}(\sigma)\big )}{\sqrt{2\rG_{3}\big (\phi_{m}(\sigma)\big )}}\beta _{3}'\bigg (\big (\phi_{3}(\sigma)\big )\bigg )}
\trace \cA(x)\rD^{2}\rd(x)\\ [1.2cm]
&
\hspace*{-.2cm}=\dfrac{\widehat{\sigma}_{0}}{\widehat{\sigma}_{0}^{2}\rE(x)-6\log \left (\dfrac{\sqrt{2}}{\sigma}\right )}
\trace \cA(x)\rD^{2}\rd(x).
\end{array}
$$
Since 
$$
\widehat{\sigma}_{0}^{2}\rE(x)-6\log \bigg (\dfrac{\sqrt{2}}{\sigma}\bigg )<0\quad \hbox{for $\sigma$ small enough}
$$
one deduces $\rC_{u,3}\in\rW^{2,\infty}$, provided  $\partial \Omega \in \cC^{4}$  and $\cA(\cdot)$ is adequately smooth.  Moreover,
$$
\begin{array}{ll}
\rC_{u,3}(x)\Theta _{3}\big (\phi_{3}(\sigma)\big )&\hspace*{-.2cm}=
\dfrac{\widehat{\sigma}_{0}\sqrt{2}\log \bigg (\dfrac{\sqrt{2}}{\sigma}\bigg )}{
\widehat{\sigma}_{0}^{2}\rE(x)-6\log \bigg (\dfrac{\sqrt{2}}{\sigma}\bigg )}\trace \cA(x)\rD^{2}\rd(x)\\ [1cm]
&=-\dfrac{\widehat{\sigma}_{0}\sqrt{2}}{6}\dfrac{1}{
1-\dfrac{\widehat{\sigma}_{0}^{2}\rE(x)}{6}\left (\log \bigg (\dfrac{\sqrt{2}}{\sigma}\bigg )\right )^{-1}}\trace \cA(x)\rD^{2}\rd(x)\\ [1cm]
\end{array}
$$
(see \eqref{eq:resto}) whence  the full expansion is
$$
\phi_{3}\big (\widehat{\sigma}_{0}\rd(x)\big )-\dfrac{\widehat{\sigma}_{0}\sqrt{2}}{6}\trace \cA(x)\rD^{2}\rd(x)+o(1).
$$
Here the assumption \eqref{eq:techbeta2} also holds. The Lowener-Nirenberg choice belongs this case if $\rN=4$.
\par
\noindent 
For $f\equiv 0,~\ell_{u}=0$ and $\Lambda =\lambda=1$ the relative result was obtained by first time in \cite{DPL}. The case $0\le \ell_{u}<1$ and $\Lambda =\lambda=1$ was obtained in  \cite{ADLR}.$\fin$
\end{exam}
\begin{exam}[Bieberbach choices]\rm One proves that for $\beta(t)=e^{t}$ the conditions
\eqref{eq:KOsecondterm} and  \eqref{eq:KOsecondtermatorigin} hold.  Now
$$
\lim_{t\nearrow \infty}\dfrac{\sqrt{2\rG(t)}}{\beta'(t)}=\lim_{t\nearrow \infty}\sqrt{2}\big (e^{-t}-e^{-2t}\big )=0.
$$
Moreover, one proves
$$
\Theta (t)=\sqrt{2}\left (2-(t+2)e^{-\frac{t}{2}}\right).
$$
Then the second order term $\rC_{u}(x)\Theta \big (\phi(\sigma)\big )$ vanishes on the boundary and the influence of the geometry is null on the boundary
(see Remak \ref{rem:phi1}). The full explosive boundary expasion is 
$$
\phi (\sigma)+o(\sigma).
$$
\par
\noindent  The conditions
\eqref{eq:KOsecondterm} and  \eqref{eq:KOsecondtermatorigin} also hold for the choice $\beta (t)=te^{t^{2}}$.  Here 
$$
\lim_{t\nearrow \infty}\dfrac{\sqrt{2\rG(t)}}{\beta'(t)}=\lim_{t\nearrow \infty}\dfrac{1-e^{-t^{2}}}{1-2t^{2}}=0.
$$
Moreover, one proves
$$
\Theta (t)=\sqrt{2}\left (12-e^{-\frac{t^{2}}{2}}\right).
$$
Therefore the second order term $\rC_{u}(x)\Theta \big (\phi(\sigma)\big )$ vanishes on the boundary and the influence of the geometry is also null on the boundary
(see Remak \ref{rem:phi1}). The full explosive boundary expasion is 
$$
\phi (\sigma)+o(\sigma).
$$
\fineq
\end{exam}
\par
\smallskip 
For general nonlinear function $\beta (t)$ the consideration of the upper terms in the explosive boundary expansion of solution under extra Keller-Osserman type growth of the sources on the boundary (see \eqref{eq:criticalgrowthextra})
is very tangled. The study for the power like choice $\beta_{m}=t^{m},~m>1$, was detailed in \cite{ADLR} when $\lambda=\Lambda=1$. We end this paper by extending the results. By simplicity in the exposition we only detail the reasonings on the equation
$$
-\trace \cA(\cdot)\rD^{2}u^{*}+\big (u^{*}\big )^{m}\le f\quad \hbox{in $\Omega$}
$$
under
$$
\limsup_{\rd(x)\searrow 0} f(x)\rd(x)^{q}\le\ell_{u} \quad (\ell_{u}>0),
$$
for $q>\dfrac{2m}{m-1},~m>1$ (see \eqref{eq:criticalgrowthextrapoweru}).  Here, we replace $\rU_{0}$, used in the proof of Theorem \ref{theo:comparisonboundaryextrapower}, by 
$$
\rU(x)\doteq \rC_{0}\sigma^{-\frac{q}{m}}\bigg ((1+\varepsilon)+\rC_{1}(x,\varepsilon,\delta)\sigma ^{\gamma}\bigg ), \quad x\in \Omega^{\delta}_{\delta_{1}}\quad (0<\delta<\delta_{1}<\delta_{\Omega}),
$$
where $\sigma\doteq \rd(x)-\delta$ and   $\rCh_{1}(\cdot)\doteq\rC_{1}(\cdot,\varepsilon,\delta):\Omega_{\delta_{1}}\rightarrow \R$ are real functions yet to be determined as well as the constant $\rC_{0}$. We will be assumed  that  $\rCh_{1}\in\rW^{2,\infty}\big (\Omega_{\delta_{1}}\big )$. 
Straightforward computations yield
$$
\begin{array}{ll}
\trace \cA(x)\rD^{2}\sigma ^{-\frac{q}{m}}&\hspace*{-.5cm}=\sigma^{-q}\left  [\dfrac{(q+m)q}{m^{2}}\rE(x)\sigma ^{\frac{q(m-1)-2m}{m}} -\dfrac{q}{m}
\trace \cA(x)\rD^{2}\rd(x)\sigma ^{\frac{q(m-1)-m}{m}}\right ],\\ [.275cm]
\trace \cA(x)\rD^{2}\rCh_{1}(x)\sigma ^{-\frac{q}{m}+\gamma}&\hspace*{-.2cm}=\sigma^{-q}
\left [\left (\dfrac{m\gamma -q}{m}\right )\left (\dfrac{m(\gamma-1)-q}{m}\right )\rE(x)\rCh_{1}(x)\sigma^{\frac{q(m-1)+m(\gamma-2)}{m}}\right .\\ [.275cm]
&+\dfrac{m\gamma -q}{m}\rFh_{1}(x)\sigma^{\frac{q(m-1)+m(\gamma-1)}{m}}
+\trace \cA(x)\rD^{2}\rCh_{1}(x)\sigma^{\frac{q(m-1)+m\gamma}{m}}\bigg ],
\end{array}
$$
where $\rFh_{1}(x)=\trace \cA(x)\big (\rD \rd(x)\otimes\rD\rCh_{1}(x)+
\rD \rCh_{1}(x)\otimes \rD \rd(x)\big )+\rCh_{1}(x)\trace\cA(x)\rD^{2}\rd(x)$ is a bounded function. On the other hand, we may write
$$
\big (\rU(x))^{m}=\rC_{0}^{m}\sigma ^{-q}\left ((1+\varepsilon)^{m}+m\rCh_{1}(x)\sigma^{\gamma } +\dfrac{m(m-1)}{2}
\rCh_{1}(x)^{2}\sigma^{2\gamma}\xi(\sigma)^{m-2}\right )
$$
for some $\xi(\sigma)$ in between $1$ and $1+\rC_{1}(x)\sigma^{\gamma}$, determinating a bounded interval.
Therefore, 
$$
\lim_{\sigma\searrow 0}\sigma^{2\gamma}\xi(\sigma)^{m-2}=0.
$$
Let us choose
\begin{equation}
\gamma =\dfrac{q(m-1)-2m}{m}>0\quad \hbox{and}\quad 
\dfrac{(q+m)q}{m^{2}}\rE(x)=m\rC_{0}^{m-1}\rCh_{1}(x).
\label{eq:choosinggamma}
\end{equation}
So that,  the assumption \eqref{eq:criticalgrowthextrapoweru} on the source $f$ implies
$$
-f(x)\ge -\rd(x)^{-q}(\ell_{u}+\varepsilon)\ge (\ell_{u}+\varepsilon)\sigma^{-q},
$$
hence
$$
-\trace \cA(\cdot)\rD^{2}\rU(x)+\big (\rU(x)\big )^{m}-f(x)\ge \sigma^{-q}
\bigg ((1+\varepsilon)^{m}\rC_{0}^{m}-(\ell_{u}+\varepsilon) +\cQ_{1}(x,\sigma)\bigg ),\quad x\in\Omega^{\delta}_{\delta_{1}},
$$
where
\begin{equation}
\begin{array}{ll}
\cQ_{1}(x,\sigma)&\hspace*{-.2cm}\doteq -\dfrac{q}{m}\trace \cA(x)\rD^{2}\rd(x)\rC_{0}\sigma ^{\frac{q(m-1)-m}{m}} \\ [.15cm]
&\hspace*{-.2cm}-\left (\dfrac{m\gamma -q}{m}\right )\left (\dfrac{m(\gamma-1)-q}{m}\right )\rE(x)\rC_{0}\rCh_{1}(x)\sigma^{\frac{q(m-1)+m(\gamma-2)}{m}}\\ [.275cm]
&-\dfrac{m\gamma-q}{m}\rC_{0}\rFh_{1}(x)\sigma^{\frac{q(m-1)+m(\gamma-1)}{m}} -\rC_{0}\trace \cA(x)\rD^{2}\rCh_{1}(x)\sigma^{\frac{q(m-1)+m\gamma}{m}}\\ [.2cm]
&+\dfrac{m(m-1)}{2}
\rCh_{1}(x)^{2}\sigma^{\frac{2\gamma +q(m-2)}{m}}\xi(\sigma)^{m-2},\quad x\in\Omega^{\delta}_{\delta_{1}}.
\end{array}
\label{eq:functionQ}
\end{equation}
The above reasonings lead to

\begin{theo}[Second order term] Let us assume $\partial \Omega \in \cC^{4}$ and $\cA(\cdot)$ adequately smooth. Suppose \eqref{eq:ellipticity} and  \eqref{eq:ellipticityneraboundaty}. Let
\begin{equation}
q>\dfrac{2m}{m-1}\quad \hbox{and}\quad \gamma =\dfrac{q(m-1)-2m}{m}>0,\quad m>1.
\label{eq:generalpowerexponent}
\end{equation}	
For any solution of 
$$
-\trace \cA(\cdot)\rD^{2}u+u^{m}\le f\quad \hbox{in $\Omega$}
$$
we have
\begin{equation}
\limsup_{\rd(x)\searrow 0}\dfrac{u^{*}(x)}{\rd(x)^{-\frac{q}{m}}\big (1+\rC_{1,u}(x)\rd(x)^{\gamma}\big )}\le \ell_{u}^{\frac{1}{m}},
\label{eq:maximalboundaryprofilesecondextragrowth}
\end{equation}
provided
$$
\limsup_{\rd(x)\searrow 0} f(x)\rd(x)^{q}\le\ell_{u} \quad (\ell_{u}>0),
$$
where
\begin{equation}
\rC_{1,u}(x)=\dfrac{(q+m)q}{m^{3}}\ell_{u}^{\frac{m}{m-1}}\rE(x).
\label{eq:choosinggamma}
\end{equation}
Analogously, for any nonnegative large solution of 
$$
-\trace \cA(\cdot)\rD^{2}v+v^{m}\ge f\quad \hbox{in $\Omega$}
$$
we have
\begin{equation}
\liminf_{\rd(x)\searrow 0}\dfrac{v_{*}(x)}{\rd(x)^{-\frac{q}{m}}\big (1+\rC_{1,d}(x)\rd(x)^{\gamma}\big )}\ge \ell_{d}^{\frac{1}{m}},
\label{eq:miniimalboundaryprofilesecondextragrowth}
\end{equation}
provided
\begin{equation}
\limsup_{\rd(x)\searrow 0} f(x)\rd(x)^{q}\ge\ell_{d}>0
\label{eq:criticalgrowthextraKOs}
\end{equation}
and
$$
\rC_{1,d}(x)=\dfrac{(q+m)q}{m^{3}}\ell_{d}^{\frac{m}{m-1}}\rE(x).
$$
\label{theo:secondarytermextraKOgrowth}
\end{theo}
\proof We only show the result relative to $u^{*}$. The proof  for $v_{*}$ is analogous (see Theorem \ref{theo:comparisonboundaryextrapower}).
With the above notation one proves
that the function $\cQ_{1}(x,\sigma)$ (see \eqref{eq:functionQ}) satisfies
$$
\lim_{\sigma \searrow 0}\cQ_{1}(x,\sigma )=0
$$
uniformly in $\Omega ^{\delta}_{\delta_{1}}$, for $\delta_{1}$ small enough. Then by choosing
$$
\rC_{0}>\dfrac{(\ell_{u}+\varepsilon)^{\frac{1}{m}}}{1+\varepsilon}
$$
we deduce
$$
-\trace \cA(\cdot)\rD^{2}\rU(x)+\big (\rU(x)\big )^{m}\ge f(x),\quad x\in\Omega^{\delta}_{\delta_{1}},
$$
for $\delta_{1}$ small. As in the proof of Theorem \ref{theo:comparisonboundaryextrapower} we deduce
$$
\limsup_{\rd(x)\nearrow 0}\dfrac{u^{*}(x)}{\rd(x)^{-\frac{q}{m}}\big ( 1+\rCh_{1}(x)\rd(x)^{\gamma}\big )}\le \rC_{0},
$$
whence letting $\rC_{0}\searrow \ell_{u}^{\frac{1}{m}}$ one concludes the proof.$\fin$
\begin{rem}\rm From Theorem \ref{theo:secondarytermextraKOgrowth} one deduces that for extra Keller-Osserman type growth of the sources on the boundary the second order in the explosive boundary expansion of large solutions is independent on the geometry whenever $\lambda=\Lambda =1$, as it was obtained in \cite{ADLR}.$\fin$
\end{rem}
\begin{rem}\rm Since
$$
-\dfrac{q}{m}+\gamma =\dfrac{q(m-2)-2m}{m}
$$
we deduce that if $\dfrac{2m}{m-1}<q<\dfrac{2m}{m-2}$ the second term of the expansion is explosive on the boundary. If $q=\dfrac{2m}{m-2}>0$ the second term is bounded on the boundary while $0<\dfrac{2m}{m-2}<q$ implies that these second term of the expansion vanishes on the boundary.$\fin$
\end{rem}

\setcounter{equation}{0}
\section{Universal bounds}
\label{sec:universalbound}

Here we adapt the proof of Lemma 1 and Theorem 1 of \cite{D} (see also the pioneer work \cite{V}) for the solution of the equation
$$
-\trace \cA(x)\rD^{2}u+\beta (u)\le f\quad \hbox{in $\Omega\subseteq\R^{\rN}$}
$$
where $\beta $ is a function satisfying \eqref{eq:KellerOsserman} and  \eqref{eq:ellipticity} holds. Let us consider the decreasing function
$$
\Phi(t)=\int^{+\infty}_{t}\dfrac{ds}{\sqrt{2\rG(s)}}\quad \hbox{where $\rG'(s)=\beta (s)$}
$$
defined in Introduction (see \eqref{eq:profilephi}). We note that $\Phi(a)< +\infty$. For each couple of positive constants $\rA$ and $\rB$ we consider the $\cC^{2}$ convex function 
$$
\Psi (\zeta)\doteq \dfrac{1}{\rA}\Phi^{-1}\left (\dfrac{\rB}{\sqrt{\Lambda}}\zeta\right ),\quad 0\le \zeta <\rR,
$$
where $\rR<\dfrac{\sqrt{\Lambda}}{\rB}\Phi(a)$ (see \eqref{eq:universalfunction}). Straightforward computations imply
$$
\Psi'(\zeta)=-\dfrac{\rB}{\rA\sqrt{\Lambda}}\sqrt{2\rG\big (\rA\Psi(\zeta)\big )}
\quad \hbox{and}\quad
\Psi''(\zeta)=\dfrac{\rB^{2}}{\rA\Lambda}\beta \big (\rA\Psi(\zeta)\big ).
$$ 
Next, we define the $\cC^{2}$ convex real function
$$
\rW(x)=\Psi\big (\zeta (x)\big ),\quad \zeta(x)=\rR^{2}-|x|^{2},\quad |x|<\rR,
$$
for which 
$$
\begin{array}{ll}
\rW''(x)&\hspace*{-.2cm}=\Psi_{\zeta \zeta}\big (\zeta(x)\big )\big (\zeta '(x)\big )^{2}+
\Psi_{\zeta}\big (\zeta(x)\big )\zeta ''(x)\\ [.2cm] 
&\hspace*{-.2cm}=4x^{2}\Psi_{\zeta \zeta}\big (\zeta(x)\big )+2\dfrac{\rB}{\rA\sqrt{\Lambda}}\sqrt{2\rG\big (\rA\Psi(\zeta(x))\big )}\\ [.2cm]
&\hspace*{-.2cm}\le 4\rR^{2}\dfrac{\rB^{2}}{\rA\Lambda}\beta \big (\Psi(\zeta(x))\big )+2\dfrac{\rB}{\rA\sqrt{\Lambda}}\sqrt{2\rA\Psi\big (\zeta(x)\big )\beta \big (\Psi\big (\zeta(x)\big )\big )},\quad |x|<\rR.
\end{array}
$$
On the other hand, if $0<\rA<1$ the inequality
$$
\dfrac{\rB}{\sqrt{\Lambda}}\rR^{2}\ge \dfrac{\rB}{\sqrt{\Lambda}} \zeta(x)\ge
\int^{\Psi(\zeta(x))}_{\rA\Psi (\zeta(x))}\dfrac{ds}{\sqrt{2\rG(s)}}\ge \dfrac{(1-\rA)\Psi\big (\zeta(x)\big )}{\sqrt{2\rG\big (\Psi\big (\zeta(x)\big )\big )}}\ge 
\dfrac{(1-\rA)\Psi\big (\zeta(x)\big )}{\sqrt{2\Psi\big (\zeta(x)\big )\beta \big (\Psi\big (\zeta(x)\big )\big )}}
$$
implies
$$
\sqrt{2\rA\Psi\big (\zeta(x)\big )\beta \big (\rA\Psi\big (\zeta(x)\big )\big )}\le \dfrac{2\rA\rB\rR^{2}}{(1-\rA)}\sqrt{\dfrac{\rA}{\Lambda}}\beta \big (\rA\Psi(\zeta(x))\big ),
$$
whence
$$
\Lambda \rW''(x)\le 4\rB^{2}\rR^{2}\left (\dfrac{1}{\rA}+\dfrac{1}{(1-\rA)\sqrt{\rA}}\right )\beta \big (\rW(x)\big ),\quad |x|<\rR.
$$
Denote 
$$
c(\rA)\doteq \dfrac{1}{\sqrt{\dfrac{1}{\rA}+\dfrac{1}{(1-\rA)\sqrt{\rA}}}},\quad 0<\rA<1.
$$
Since $c(0)=c(1)=0$ the continuous function $c(\rA)$ attains a positive maximum at some $\rA_{\infty}$ in the interval $(0,1)$. Then for
$$
\rB=\dfrac{c(\rA_{\infty})}{2\rR^{2}}
$$
one has
\begin{equation}
\Lambda\rW''(x)\le \left (\dfrac{2\rB\rR}{c(\rA_{\infty})}\right )^{2}\beta \big (\rW(x)\big ),\quad |x|<\rR.
\label{eq:supersolutionuniversal}
\end{equation}
\begin{theo} Assume \eqref{eq:ellipticity} and \eqref{eq:KellerOsserman}. Then there exists two positive constants $\rA_{\infty}$ and $\rB$, satisfying \eqref{eq:ajusteuniversalbound} below, for which
\begin{equation}
	u^{*}(x)\le \sum_{i=1}^{\rN}\Psi\big (\rR^{2}-|x_{i}|^{2}\big )+\beta^{-1}\big (\n{f}_{\cQ_{\rR}(x_{0})}\big ),\quad x\in \cQ_{\rR}(x_{0})
	\label{eq:universalestimate}
\end{equation}
holds in any cube $\cQ_{\rR}(x_{0})\subset \subset \Omega $, for any subsolution $u$ of
$$
-\trace \cA(\cdot )\rD^{2}u+\beta (u)\le f\quad \hbox{in $\Omega$}.
$$ 
\label{theo:universalestimate}
\end{theo}
\proof Assume with no loss of generality $x_{0}=0$. The above reasoning proves that the function
$$
\cW(x)=\sum_{i=1}^{\rN}\Psi\big (\zeta(x_{i})\big ),\quad \zeta(x_{i})=\rR^{2}-|x_{i}|^{2},\quad |x_{i}|<\rR,
$$
satisfies in the cube $\cQ_{\rR}(0)$ the inequality
$$
\trace \cA(x)\rD^{2}\cW(x)\le
\Lambda \sum_{i=1}^{\rN}\rW'' (x_{i})
\le   \left (\dfrac{2\rB\rR}{c(\rA_{\infty})}\right )^{2}\sum_{i=1}^{\rN}\beta \big (\Psi (\zeta (x_{i}))\big )
\le  \rN\left (\dfrac{2\rB\rR}{c(\rA_{\infty})}\right )^{2}\beta \big (\cW(x)\big ),
$$
thus
\begin{equation}
-\trace \cA(x)\rD^{2}\cW(x)+\beta \big (\cW(x)\big )\ge 0,\quad x\in\cQ_{\rR}(0),
\label{eq:supersolutionuniversal}
\end{equation}
provided
\begin{equation}
\rN\left (\dfrac{2\rB\rR}{c(\rA_{\infty})}\right )^{2}=1.
\label{eq:ajusteuniversalbound}
\end{equation}
Since $\cW(x)=+\infty $, Remark \ref{rem:supverification} concludes the estimate.$\fin$
\begin{rem}\rm We emphasize that no uniform ellipticity assumption is required in the above proof. On the other hand, the choice \eqref{eq:ajusteuniversalbound} is independent on the constant $\Lambda.\fin$
\label{rem:universalestimate}
\end{rem} 
\begin{rem}\rm Relative to the equation
$$
-\varepsilon \Delta u_{\varepsilon} -\trace \cA(\cdot)\rD^{2}u_{\varepsilon} +	\beta (u_{\varepsilon}) =f\quad\hbox{in $\Omega$}
$$
with $0<\varepsilon <1$, the estimate \eqref{eq:universalestimateepsilonn} follows from \eqref{eq:universalestimate} by replacing $\Lambda$ by $\Lambda +\varepsilon$ in the definition~\eqref{eq:universalfunction} of the function $\Phi(\zeta)$. The estimate \eqref{eq:universalestimateepsilon} also follows from \eqref{eq:universalestimate} by  replacing $\Lambda$ by $\Lambda +1$.
Indeed \eqref{eq:supersolutionuniversal} becomes
$$	
(\varepsilon +\Lambda )\rW''(x)\le (1+\Lambda )\rW''(x)\le \left (\dfrac{2\rB\rR}{c(\rA_{\infty})}\right )^{2}\beta \big (\rW(x)\big ),\quad |x|<\rR,
$$
whence one deduces
$$
-\varepsilon \Delta \cW(x)-\trace \cA(x)\rD^{2}\cW(x)+\beta \big (\cW(x)\big )\ge 0,\quad x\in\cQ_{\rR}(0).
$$
\fineqnum
\label{rem:universalestimateestimate}
\end{rem}

\setcounter{equation}{0}
\section{Strong Maximum Principle  without coercivity term} 
\label{sec:SMP}
\par
Here we group some results used in above reasoning. So, in the proof of Theorem  \ref{theo:absoluteboundarycomparison} we require a version of Weak Maximum Principle without coercivity term.
\begin{theo}[Weak Maximum Principle without coercivity term (I)] 
Assume \eqref{eq:cA}. Let $u$ be a discontinuous subsolution of 
$$
-\trace \cA(\cdot)\rD^{2}u=f, \quad x\in \Omega
$$
where $f$ is a uniformly continuous function and $v$ be a discontinuous supersolution of
$$
-\trace \cA(\cdot)\rD^{2}v=g, \quad x\in \Omega 
$$
where $g$ is a continuous function. Then 
$$
u^{*}(x)-v_{*}(x)\le 
\sup_{\partial \Omega}~(u^{*}-v_{*})+\dfrac{\rC_{\Omega}}{2\lambda_{\cA}}\n{(f-g)_{+}}_{\rL_{\infty}},\quad x\in\Omega,
$$
where $\disp \rC_{\Omega}=\max_{x\in \Omega}|x|^{2}$. See the Introduction for  the positive constant $\lambda_{\cA}$.
\label{theo:withoutcoercivityI}
\end{theo}
\proof It is clear that the smooth function $\rW(x)=-|x|^{2},~x\in \Omega$, satisfies 
$$
-\trace \cA(x)\rD^{2}\rW(x)\ge 2\lambda_{\cA},\quad x\in\Omega.
$$
So,  for any positive constant $a$ we have
$$
-\trace \cA(\cdot)\rD^{2}(v_{*}+a\rW)\ge g+2a\lambda_{\cA}, \quad x\in \Omega.
$$
Arguing on an eventual interior maximum of $u^{*}-(v_{*}+a\rW)$ at some $x_{0}\in\Omega$, as in the Theorem~ \ref{theo:comparisonprincipleviscosity}, we get 
$$
0\le \omega\big (\alpha|x_{\alpha}-y_{\alpha}|^{2}\big )+\omega_{f}\big (|x_{\alpha}-y_{\alpha}|\big )+f(y_{\alpha})-g(y_{\alpha})-2a\lambda_{\cA}
$$
(see \eqref{eq:capitalinequality}). Then  by letting $\alpha\nearrow \infty$ one derives the contradiction
$$
0\le f(x_{0})-g(x_{0})-2a\lambda_{\cA}<0
$$
provided $a=\dfrac{\n{(f-g)_{+}}_{\infty}}{2\lambda_{\cA}}+\varepsilon$ with $\varepsilon>0$. So that, we have obtained the inequality
$$
\big (u^{*}-v_{*}\big )(x)\le \big (u^{*}-v_{*}-a\rW\big )(x)\le \sup_{\partial \Omega}~\big (u^{*}-v_{*}-a\rW\big ) \le 
\sup_{\partial \Omega}~(u^{*}-v_{*})+\left (\dfrac{\n{(f-g)_{+}}_{\infty}}{2\lambda_{\cA}}+\varepsilon\right )\rC_{ \Omega},
$$
for $x\in\Omega$. The result follows sending $\varepsilon\searrow 0.\fin$
\par
\medskip
Next we detail the application of certain results and reasonings of \cite{BBus} in what follows. The classical device by E. Hopf (see \cite{GT}) enables us to obtain a Strong Maximum Principle as follows. It was used in the proof of Theorem \ref{theo:existenceKOsources}.

\begin{lemma}[Hopf Lemma] Assume \eqref{eq:cA} and \eqref{eq:structureannulusSMPsub}.
Let $u$ be a discontinuous subsolution of 
$$
-\trace \cA(\cdot)\rD^{2}u=f, \quad x\in \Omega
$$
and $v$ be a discontinuous supersolution of
$$
-\trace \cA(\cdot)\rD^{2}v=f, \quad x\in \Omega
$$
where $f$ is a uniformly continuous function. Suppose that $\Omega$ is an open set such that there exists $x_{0}\in\partial \Omega$ for which
$$
\bB_{\rR}(y)\subset \Omega\quad\hbox{and}\quad \partial \bB_{\rR}(y)\cap\partial \Omega=\{x_{0}\} \quad\hbox{{\rm (}interior sphere condition at $x_{0}\in\partial \Omega${\rm )}}
$$
and
$$
\big (u^{*}-v_{*}\big )(x)\le \big (u^{*}-v_{*}\big )(x_{0}),\quad x\in\bB_{\rR}(y)
$$
hold. Then
$$
\liminf_{x\rightarrow x_{0}}\dfrac{\big (u^{*}-v_{*}\big )(x_{0})-\big (u^{*}-v_{*}\big )(x)}{|x-x_{0}|}>0.
$$
\label{lemma:Hopf}
\end{lemma}
\proof The smooth and positive function 
$$
\rW(x)=e^{-\alpha|x-y|^{2}}-e^{-\alpha \rR^{2}},\quad x\in\bB_{\rR}(y)\subset \Omega
$$
satisfies $\rD \rW(x)=-2\alpha e^{-\alpha |x-y|^{2}}(x-y) $ whence
$$
\left \{
\begin{array}{l}
\alpha \rR e^{-\alpha \rR^{2}}\le |\rD \rW(x)|=2\alpha e^{-\alpha|x-y|^{2}}|x-y|\le 2\rR\alpha e^{-\alpha\frac{\rR^{2}}{4}},\\ [.15cm]
\begin{array}{ll}
	\rD^{2}\rW(x)&\hspace*{-.2cm}= -2\alpha e^{-\alpha |x-y|^{2}}\rI+4\alpha ^{2}e^{-\alpha |x-y|^{2}}\big (x-y\big )\otimes \big (x-y\big )\\ [.15cm]
	&\ge -\dfrac{2}{\rR}|\rD \rW(x)|\rI+\alpha\rR\dfrac{\rD \rW(x)\otimes \rD \rW(x)}{|\rD \rW(x)|}
\end{array}
\end{array}
\right . \qquad x\in \bB_{\rR}(y)\setminus\bB_{\frac{\rR}{2}}(y).
$$
Then from \eqref{eq:structureannulusSMPsub} we have 
$$
\trace \cA(x)\rD^{2}\rW(x)\ge 0,\quad x\in\bB_{\rR}(y)\setminus \bB_{\frac{\rR}{2}}(y),
$$
for some provided $\alpha$. Therefore
$$
-\trace \cA(\cdot)\rD^{2}(u^{*}+\varepsilon \rW)\le f\quad \hbox{in $\bB_{\rR}(y)\setminus\bB_{\frac{\rR}{2}}(y)$},
$$
for any $\varepsilon>0$. Moreover
$$
(u^{*}-v_{*}+\varepsilon \rW)(x)\le (u^{*}-v_{*})(x_{0}),\quad x\in\partial \bB_{\frac{\rR}{2}}(y)
$$
if
$$
\varepsilon=\dfrac{(u^{*}-v_{*})(x_{0})-\sup_{x\in\partial \bB_{\frac{\rR}{2}}(y)}(u^{*}-v_{*})(x)}{e^{-\alpha\frac{\rR^{2}}{4}}-e^{-\alpha \rR^{2}}}>0.
$$
So that, the Weak Maximum Principle (Theorem \ref{theo:withoutcoercivityI}) applied to the set $\bB_{\rR}(y)\setminus \overline{\bB}_{\frac{\rR}{2}}(y)$ leads to
$$
(u^{*}-v_{*})(x_{0})-(u^{*}-v_{*})(x)\ge \rW(x),\quad \bB_{\rR}(y)\setminus \overline{\bB}_{\frac{\rR}{2}}(y).
$$
Sinces $\rW(x_{0})=0$ one concludes
\begin{equation}
\liminf_{x\rightarrow x_{0}}\dfrac{(u^{*}-v_{*})(x_{0})-(u^{*}-v_{*})(x)}{|x-x_{0}|}\ge\varepsilon \pe{\rD \rW(x_{0})}{\bn(\rR)}=\varepsilon\alpha e^{-\alpha \rR^{2}}>0,
\label{eq:Hopf}
\end{equation}
where $\bn(\rR)$ is the outer normal vector to $\partial \Omega$ at $x_{0}.\fin$
\begin{rem}\rm In fact, straightforward computations show that there exists $\delta >0$ for which 
$$
\liminf_{x\rightarrow x_{0}\atop \angle (x_{0}-x,\bn(\rR))<\frac{\pi}{2}-\delta
}\dfrac{\big (u^{*}-v_{*}\big )(x_{0})-\big (u^{*}-v_{*}\big )(x)}{|x-x_{0}|}>0.
$$
\fineqnum
\end{rem}
\begin{rem}\rm Even in the broader case in which $\alpha$ is large enough, assumption \eqref{eq:structureannulusSMPsub} is more restrictive to \eqref{eq:ellipticity}.$\fin$
\end{rem}

\begin{theo}[Strong Maximum Principle without coercivity term] Assume \eqref{eq:lightellipticity}
and \eqref{eq:structureannulusSMPsub}. Let $u$ be a discontinuous subsolution of 
$$
-\trace \cA(\cdot)\rD^{2}u=0, \quad x\in \Omega
$$
and $v$ be a discontinuous supersolution of
$$
-\trace \cA(\cdot)\rD^{2}v=0, \quad x\in \Omega,
$$
Then if $u^{*}-v_{*}$ achieves its maximum in the interior of $\Omega$, then 
$u^{*}-v_{*}$ is constant in $\Omega$. 
\label{theo:SMPwithoutcoercivity}
\end{theo}

Before the proof of Theorem \ref{theo:SMPwithoutcoercivity} we collect some useful results. First we recall a property of any convex function, $\psi $, defined in $\cO$: it attains a local maximum at $z_{0}\in \cO$ then $\rD \psi(z_{0})=\bcero$.  More precisely
\begin{lemma}[DM] Let $\psi$ be a function achieving a local maximum at some $z_{0}\in\cO$. Assume that there exists a function  $\widehat{\psi}$ defined in $\cO$ such that $\widehat{\psi}(z_{0})=0,~\Psi=\psi+\widehat{\psi}$ is convex on $\cO$ and
$$
\widehat{\psi}(x)\le \rK|x-z_{0}|^{2},\quad x\in\cO~\hbox{with $|x-z_{0}|$ small,}
$$
for some constant $\rK>0$. Then the function $\psi$ is differentiable at $z_{0}$ and $\rD \psi(z_{0})={\bf 0}$.
\label{lemma:DM}
\end{lemma}
\proof  By simplicity we can take $z_{0}=0\in\cO$. By applying the  Convex Separation Theorem there exists $\bp\in\R^{\rN}$ such that
$$
\Psi(x)\ge \Psi(0)+\langle \bp,x\rangle=\psi(0)+\langle \bp,x\rangle, \quad x\in\cO,~\hbox{with $|x|$ small}.
$$
Then we have
\begin{equation}
\begin{array}{ll}
	\psi(x)&\hskip-.2cm=\Psi(x)-\widehat{\psi}(x)\ge \psi(0)+\langle \bp,x\rangle-\rK|x|^{2}\\ [.15cm]
	&\hskip-.2cm \ge \psi (x)+\langle \bp,x\rangle -\rK|x|^{2},\quad x\in \cO~\hbox{with $|x|$ small}
\end{array}
\label{eq:semiauxi2}
\end{equation}
whence
$$
\langle \bp,x\rangle\le \rK|x|^{2},\quad x\in\cO\hbox{ with $|x|$ small}.
$$
For $\tau >0$ small enough we can choose $x=\tau \bp\in \cO$ and $\tau \rK <1$, for which
$$
\tau |\bp|^{2}\le \rK \tau^{2} |\bp|^{2}.
$$
Therefore $\bp=\bcero $. Finally, (\ref{eq:semiauxi2}) leads to
$$
0\ge \psi(x)-\psi(0)\ge -\rK|x|^{2},\quad x\in \cO~\hbox{ with $|x|$ small},
$$
and the result follows. $\fin$
\begin{rem}\rm As it was pointed out the result is immediate if $\psi$ is convex for which we can choose $\widehat{\psi}\equiv~0.\fin$ 
\end{rem}

As is in the uniform ellipticity case, the proof of Theorem \ref{theo:SMPwithoutcoercivity} rests on a contradiction generated on the differentiability at a local interior maximum. It is more direct whenever (DM) property holds. We make it by regularizing $u^{*}$ and $v_{*}$ by sup- and inf- convolution
$$
u_{\varepsilon}(x)\doteq \sup_{y\in \Omega} \left \{u^{*}(y)-\dfrac{|y-x|^{2}}{2\varepsilon^{2}}\right \}\quad \hbox{and}\quad 
v^{\varepsilon}(x)\doteq \inf_{y\in \Omega} \left \{v_{*}(y)+\dfrac{|y-x|^{2}}{2\varepsilon^{2}}\right \}.
$$
By construction, $u_{\varepsilon}$ and $v^{\varepsilon}$ are semiconvex and semiconcave function, respectively, in a slightly smaller domain (still denote by $\Omega$ in the sequel for simplicity) where we argue. Moreover, the magical properties of the regularization implies that they are sub and supersolution of 
$$
-\trace \cA(\cdot)\rD^{2}u=0, \quad x\in \Omega.
$$
respectively (see \cite{CIL}). In fact, each of them satisfies the Strong Maximum Principle as follows from 
\begin{coro} Assume \eqref{eq:structureannulusSMPsub}. Let $w$ be a subsolution (respectively supersolution) of 
$$
-\trace \cA(\cdot)\rD^{2}u=0, \quad x\in \Omega.
$$
Assume  that $w$ achieves a local maximum (resp. a local minimum) at $z_{0}\in\Omega$  then $w$ is constant in a neighborhood of $z_{0}$.
\label{coro:SMPwithoutcoercivitysemiconvex}
\end{coro}
\proof For simplicity we only treat the case of a subsolution. First of all
$$
w_{\varepsilon}(z_{0})\ge w^{*}(z_{0})\ge u^{*}(y)-\dfrac{|y-x|^{2}}{2\varepsilon ^{2}}\quad \hbox{for all $y,~x$}
$$
implies $w_{\varepsilon}(z_{0})\ge w_{\varepsilon}(x)$, thus $w_{\varepsilon}$ also attains a maximum at $z_{0}$. By means of geometrical construction we may assume that $z_{0}\in\partial \bB$ for an adequate ball $\bB\subset \Omega$. The reasoning of the proof of Lemma \ref{lemma:Hopf} applied to $w_{\varepsilon}$ derives the contradiction
$$
\bcero \not = \rD w_{\varepsilon}(z_{0})=\bcero.
$$
So that, $w_{\varepsilon}$ must be constant in a ball $\bB(z_{0})$. Then denoting $\widehat{x}$ the point such that 
$$
w_{\varepsilon}(x)=w^{*}(\widehat{x})-\dfrac{|\widehat{x}-x|^{2}}{2\varepsilon^{2}}
$$ 
we have
$$
w^{*}(\widehat{x})-\dfrac{|\widehat{x}-x|^{2}}{2\varepsilon^{2}}=w_{\varepsilon}(x)=
w_{\varepsilon}(z_{0})=w^{*}(z_{0})\ge w^{*}(\widehat{x}),\quad x\in \bB(z_{0}).
$$
Therefore $\widehat{x}=x$ and $w^{*}(z_{0})= w^{*}(x),~x\in\bB(z_{0}).\fin$
\par
\medskip
\noindent
{\sc Sketch of Proof of Theorem \ref{theo:SMPwithoutcoercivity}}  Let us assume 
$$
\max_{\overline{\Omega}}(u^{*}-v_{*})=(u^{*}-v_{*})(x_{0})>0
$$
at some $x_{0}\in\Omega$. Regularizing $u^{*}$ and $v_{*}$ by sup- and inf-convolution respectively we may assume 
$$
\rM(0)\doteq \max_{\overline{\Omega}}(u_{\varepsilon}-v^{\varepsilon})=(u_{\varepsilon}-v^{\varepsilon})(x_{0})>0
$$ 
for $\varepsilon>0$ small enough. Next, we introduce
\begin{equation}
\rM(\bh)\doteq \max_{x\in \overline{\Omega}_{|\bh|}}\big (u_{\varepsilon}(x+\bh)-v^{\varepsilon}(x)\big )=u_{\varepsilon}(x_{\bh}+\bh)-v^{\varepsilon}(x_{\bh})>0,
\label{eq:maximumMh}
\end{equation}
for some $x_{\bh}\in \Omega_{|\bh|}\doteq \{x\in\Omega:~\rd(x)>|\bh|\}$. If $|\bh|$ is small enough we can assume $x_{0}\in\Omega_{|\bh|}$. Notice that $\rM(\bh)$ is convex in a neighborhood of $\bcero\in\R^{\rN}$ due to it is the maximum over a family of semiconvex functions.
\par
From the property (DM) (see Lemma \ref{lemma:DM}) the function
$u_{\varepsilon}(\cdot+\bh)-v^{\varepsilon}(\cdot)$ is differentiable at $x_{\bh}$ and 
$$
\rD u_{\varepsilon}(x_{\bh}+\bh)=\rD v_{\varepsilon}(x_{\bh}).
$$
We assume the key stone 
\begin{equation}
\rD v_{\varepsilon} (x_{\bh})=\bcero
\label{eq:keystoineSMP}
\end{equation}
proved by the sharp reasonings of the points 4,5 and 6 of the proof of Theorem 1 of \cite{BBus} (see Lemma~\ref{lemma:nullgradient} below). Since
$\rD u_{\varepsilon}(x_{\bh}+\bh)=\rD v_{\varepsilon}(x_{\bh})=\bcero$ for any $\bh'$ in a neighborhood of $\bh$ we have
$$
\rM(\bh')\ge u_{\varepsilon}(x_{\bh}+\bh')-v^{\varepsilon}(x_{\bh})\ge u_{\varepsilon}(x_{\bh}+\bh)-v^{\varepsilon}(x_{\bh})-\rC|\bh-\\bh'|^{2}
$$
whence
$$
\rM(\bh')\ge \rM(\bh)-\rC|\bh-\bh'|^{2}.
$$
Therefore $\bcero \in\partial \rM(\bh)$ (the subdifferential of $\rM$ at $\bh$) for any $\bh$ in a neighborhood of $\bcero$ in which $\rM(\bh)$ must be constant. So that
$$
u_{\varepsilon}(x_{0})-v^{\varepsilon}(x_{0})=\rM(0)=\rM(\bh)\ge u_{\varepsilon}(x_{0}+\bh)-v^{\varepsilon}(x_{0})
$$
implies that $x_{0}$ is also a point of local maximum for $u_{\varepsilon}$.
\par
\noindent 
Finally, from Corollary \ref{coro:SMPwithoutcoercivitysemiconvex} the function $u_{\varepsilon}$ is constant in a neighborhood of $x_{0}$ in which
$$
v^{\varepsilon}(x)=v^{\varepsilon}(x)-u_{\varepsilon}(x)+u_{\varepsilon}(x)\ge 
v^{\varepsilon}(x_{0})-u_{\varepsilon}(x_{0})+u_{\varepsilon}(x)=v^{\varepsilon}(x_{0}).
$$
Since $x_{0}$ is a point of local minimum, the function  $v^{\varepsilon}$ is constant in a neighborhood of $x_{0}$ (see again Corollary \ref{coro:SMPwithoutcoercivitysemiconvex}) and consequently $u_{\varepsilon}-v^{\varepsilon}$ is constant in a neighborhood of $x_{0}.\fin$
\par
\medskip
A  direct consequence of Theorem \ref{theo:SMPwithoutcoercivity} is a version of  Theorem \ref{theo:withoutcoercivityI}.
\begin{theo}[Weak Maximum Principle without coercivity term (II)] Assume \eqref{eq:lightellipticity}
and \eqref{eq:structureannulusSMPsub}. Let $u$ be a discontinuous subsolution of 
$$
-\trace \cA(\cdot)\rD^{2}u=0, \quad x\in \Omega
$$
and $v$ be a discontinuous supersolution of
$$
-\trace \cA(\cdot)\rD^{2}v=0, \quad x\in \Omega .
$$
Then 
$$
u^{*}(x)-v_{*}(x)\le  \sup_{\partial \Omega}~(u^{*}-v_{*}),\quad x\in\Omega.
$$
\label{theo:withoutcoercivityII}
\end{theo}
\proof First of all, with no loss of generality we may assume $\sup_{\partial \Omega}~(u^{*}-v_{*})\le 0$ by replacing $v^{*}$ by $v^{*} +\sup_{\partial \Omega}~(u^{*}-v_{*})$. So that we will prove that a condition as
$$
\max_{\overline{\Omega}}(u_{*}-v_{*})>0
$$
contradicts the fact $u^{*}-v_{*}\le 0$ on $\partial \Omega$.  As in the proof of Theorem \ref{theo:SMPwithoutcoercivity} it is enough to work with the approximations $u_{\varepsilon}$ and $v^{\varepsilon}$ and to argue in a slightly smaller domain, also denote by $\Omega$, in which
$$
\rM(0)=\max_{\overline{\Omega}}(u^{\varepsilon}-v_{\varepsilon})>0
$$ 
and $u_{\varepsilon}-v^{\varepsilon}< 0$ on $\partial \Omega$. From the consequence of Theorem \ref{theo:SMPwithoutcoercivity} we deduce that $\rM(0)$ is achieved in an open set that by construction 
it is also closed. So that $\rM(0)>0$ is achieved on the whole $\overline{\Omega}$ that contradicts the fact  $u^{*}-v_{*}\le 0$ on $\partial \Omega.\fin$

\begin{lemma} Assumed  \eqref{eq:lightellipticity} the function $v^{\varepsilon}$ defined in  the proof of Theorem \ref{theo:SMPwithoutcoercivity} safisfies \eqref{eq:keystoineSMP}
\label{lemma:nullgradient}
\end{lemma}
\proof The lack of a coercive term can be supplied by means of a kind of Kruzkov change of variable (see \cite{BBus} or \cite{BDD}). More precisely, let $u\in\cC^{2}$ and define
$$
\rU=\psi(u)
$$
where $\psi$ is a smooth function close to the identity map. Since
\begin{equation}
\rD u=\varphi '(\rU)\rD \rU\quad \hbox{and}\quad \rD^{2}u=\varphi' (\rU)\rD^{2}\rU+\varphi''(\rU)\rD \rU\otimes\rD \rU, 
\label{eq:derivatives}
\end{equation}
for $\varphi =\psi^{-1}$, the classical equation
$$
-\trace \cA(x)\rD^{2}u(x)=0,\quad x\in\Omega
$$
becomes
$$
-\trace \widehat{\cA}_{\varphi}\big (x,\rU(x),\rD \rU(x),\rD ^{2}\rU(x)\big )=0,\quad x\in\Omega,
$$
where
\begin{equation}
\widehat{\cA}_{\varphi}(x,t,\bq,\rY)\doteq \dfrac{1}{\varphi '(t)}\cA(x)\rX=\cA(x)\left (\rY+\dfrac{\varphi ''(t)}{\varphi '(t)}\bq \otimes \bq\right ),\quad (x,t,\bq,\rY)\in\Omega\times\R\times\R^{\rN}\times \cS^{\rN}_{+}.
\label{eq:Aphi}
\end{equation}
From \eqref{eq:derivatives} we note that the new variable $\bq$ correspond with 
the old one $\bp=\varphi '(t)\bq$.
Then,
$$
\dfrac{\partial }{\partial t}\widehat{\cA}_{\varphi}(x,t,\bq,\rY)=\left (\dfrac{\varphi ''(t)}{\varphi '(t)}\right )'\cA (x)\bq\otimes \bq \bigg |_{\bq=\frac{1}{\varphi '(t)}\bp}=\left (\dfrac{\varphi'''(t)}{(\varphi '(t))^{3}}-\dfrac{{(\varphi ''(t))^{2}}}{(\varphi '(t))^{4}}\right )\cA(x)\bp\otimes \bp.
$$
As in \cite{BBus} an adequate approximation of the identity map  can be constructed by means of the relationship given by $\omega (t) = \left (\dfrac{\varphi _{\eta}''\big (\psi _{\eta}(t)\big )}{\varphi ' \big (\psi _{\eta}(t)\big )}\right )'$ with $
\omega (t)\doteq \exp \left (-\eta ^{-1}\big (t+\eta^{-1}\big )\right )$ that implies
\begin{equation}
\dfrac{\partial }{\partial t}\trace \widehat{\cA}_{\varphi}(x,t,\bq,\rY)=\omega '(t)\trace \cA (x)\bp\otimes \bp<\omega '(t)\lambda (\delta)<0
\label{eq:coercivity}
\end{equation}
provided $\bp\neq\bcero$. This choice gives 
$$
\varphi _{\eta}'(t)=\exp \left (\int^{t}_{0}\omega (s)ds\right )
$$
that satisfies $\varphi _{\eta}(t)\rightarrow t,~\varphi _{\eta}'(t)\rightarrow 1$ and $\varphi _{\eta}''(t)\rightarrow 0$ locally uniformly in $\R$ as $\eta \rightarrow 0$.
\par
\noindent So that, we return to the proof of Theorem \ref{theo:SMPwithoutcoercivity}. Let us introduce the functions
$$
\rU_{\varepsilon,\eta}=\psi_{\eta}\big (u_{\varepsilon}\big)\quad \hbox{and}\quad
\rV^{\varepsilon,\eta}=\psi_{\eta}\big (v^{\varepsilon}\big)
$$
sub and supersolution, respectively, of
$$
-\trace \widehat{\cA}_{\eta}(x,\rU,\rD\rU,\rD^{2}\rU)=0\quad \hbox{in $\Omega$},
$$ 
where $\widehat{\cA}_{\eta}$ is defined by replacing $\varphi$ by $\varphi_{\eta}$  in \eqref{eq:Aphi} (see \cite{CIL} for viscosity solution after change of variables).
\par
\noindent 
Les us assume that there exists a sequence $\{\bh_{n}\}_{n}\rightarrow \bcero$ such that
\begin{equation}
\rD u_{\varepsilon}(x_{\bh_{n}}+\bh_{n})=\rD v^{\varepsilon}(x_{\bh_{n}})\neq \bcero,
\label{eq:hipcomtradiction}
\end{equation}
where $x_{\bh_{n}}$ is the maximum point introduced in \eqref{eq:maximumMh}. The convexity properties of $u_{\varepsilon}$ and $v^{\varepsilon}$ enable us  to use the {\em partial continuity of the gradient} (PCG) and to prove
\begin{equation}
|\rD u_{\varepsilon}(x_{\bh_{n}}+\bh_{n})|=|\rD v^{\varepsilon}(x_{\bh_{n}})|\ge c(n)>0
\label{eq:hipcomtradiction2}
\end{equation}
for some constant $c(n)$. On the other hand, by construction the function $ \rU_{\varepsilon,\eta}(\cdot+\bh_{n})-\rV^{\varepsilon,\eta}(\cdot)$ achives a positive maximum point at some $x_{\eta}\in\Omega$, for $\eta $ small enough. Using (PCG) we have
$$
|\rD \rU_{\varepsilon,\eta}(x_{\eta}+\bh_{n})|=|\rD \rV^{\varepsilon,\eta}(x_{\eta})|.
$$
Extracting a subsequence $\eta_{k}\rightarrow 0$ such that $x_{\eta_{k}}\rightarrow \overline{x}$, a maximum point of $u_{\varepsilon}(\cdot+\bh_{n})- v^{\varepsilon}(\cdot)$ for which \eqref{eq:hipcomtradiction2}  implies
$$
|\rD \rU_{\varepsilon,\eta}(\overline{x}+\bh_{n})|=|\rD \rV^{\varepsilon,\eta}(\overline{x})|\ge \dfrac{c(n)}{2}>0
$$
for $\eta$ small enough. By means of a sharp argument one proves in the point 6 of the proof of Theorem~1 of \cite{BBus} that $\rD \rV^{\varepsilon,\eta}(z)\rightarrow \widehat{\bp}\in\R^{\rN},~|\widehat{\bp}|>\dfrac{c(n)}{4}>0$ and $\rD^{2} \rV^{\varepsilon,\eta}(z)\rightarrow \widehat{\rX}\in\cS^{\rN}_{+}$ where $z$ is a maximum point of a suitable approximation of $\rU_{\varepsilon,\eta}(\cdot+\bh_{n})- \rV^{\varepsilon,\eta}(\cdot)$ such that $z\rightarrow x_{\eta}$. Moreover, since $
\rU_{\varepsilon,\eta}=\psi_{\eta}\big (u_{\varepsilon}\big)$ and $
\rV^{\varepsilon,\eta}=\psi_{\eta}\big (v^{\varepsilon}\big)$
are sub and supersolution, taking limit in these approximation one proves in \cite{BBus} 
$$
-\trace \widehat{\cA}\big (x_{\varepsilon},\rU_{\varepsilon,\eta}(x_{\eta}+\bh_{n}),\widehat{\bp},\widehat{\rX}
\big )\le 0\le -\trace \widehat{\cA}\big (x_{\varepsilon},\rV^{\varepsilon,\eta}(x_{\eta}),\widehat{\bp},\widehat{\rX}
\big ),
$$
whence \eqref{eq:coercivity} leads to the contradiction
\begin{equation}
\trace \widehat{\cA}\big (x_{\varepsilon},\rU_{\varepsilon,\eta}(x_{\eta}+\bh_{n}),\widehat{\bp},\widehat{\rX}
\big )>\trace \widehat{\cA}\big (x_{\varepsilon},\rU_{\varepsilon,\eta}(x_{\eta}+\bh_{n}),\widehat{\bp},\widehat{\rX}
\big ). 
\label{eq:contradiction}
\end{equation}
\fineqnum
\begin{rem}\rm Notice that \eqref{eq:contradiction} is deduced from  \eqref{eq:coercivity} that  requires $\widehat{\bp}\neq \bcero$. It is proved from the contradiction assumption \eqref{eq:hipcomtradiction}.$\fin$
\end{rem}

\begin{center}
{\bf \sffamily Acknowledgements}	
\end{center}
The research of G. D\'{\i}az was partially supported
by the project ref. PID2020-112517GB-I00 of the Ministerio de Ciencia e Innovaci\'{o}n
(Spain).

{\footnotesize
\begin{tabular}{l}
Gregorio Díaz \\
Dpto. Análisis Matemático y Matem\'atica Aplicada\\
U. Complutense de Madrid \\
28040 Madrid, Spain \\
{\tt gdiaz@ucm.es}, {\tt gregoriodiazdiaz@gmail.com} 
\end{tabular}
}

\end{document}